\documentclass[11pt]{article}
\usepackage[latin1]{inputenc}
\usepackage[english]{babel}
\newtheorem{theorem}{Theorem}

\newtheorem{definition}{Definition}

\usepackage{stmaryrd}
\usepackage{enumerate}
\usepackage[noend]{algpseudocode}
\usepackage{calrsfs}
\usepackage{graphicx}
\usepackage{amssymb}
\usepackage{array}
\usepackage{amsmath}
\usepackage{epstopdf}
\usepackage{epsfig}
\usepackage{subcaption}
\usepackage{pstricks}
\usepackage{fancyheadings}
\usepackage{pgfplots}
\usepackage{physics}
\usepackage{stackengine}
\usepackage[margin=1.0in]{geometry}

\usepackage[framemethod=tikz]{mdframed}
\usepackage{lipsum}

\usepackage{url}
\usepackage{xcolor}
\usepackage{framed,multirow}

\usepackage{latexsym}

\usepackage{amsfonts}
\usepackage{mathtools}
\usepackage{bm}
\usepackage{float}

\definecolor{mycolor}{rgb}{0.122, 0.435, 0.698}
\definecolor{blue-green}{rgb}{0.0, 0.87, 0.87}
\definecolor{royalblue}{rgb}{0.01, 0.28, 1.0}
\definecolor{bluet}{rgb}{0.1, 0.1, 0.8}

\newmdenv[innerlinewidth=1pt, roundcorner=0.5pt,linecolor=mycolor,innerleftmargin=0.5pt,
innerrightmargin=0.5pt,innertopmargin=0.4pt,innerbottommargin=0.4pt]{mybox}

\DeclareMathAlphabet{\pazocal}{OMS}{zplm}{m}{n}

\newtheorem{assumption}{Assumption}
\newtheorem{corollary}{Corollary}
\newtheorem{remark}{Remark}

\newcommand{\lJump}{[\![}
\newcommand{\rJump}{]\!]}

\makeatletter
\def\BState{\State\hskip-\ALG@thistlm}
\makeatother

\begin{document}
\title{The perfectly matched layer (PML) for hyperbolic wave propagation problems: A review}
\author{Kenneth Duru\thanks{Corresponding author. Mathematical Sciences Institute, The Australian National University, Australia.} \and Gunilla Kreiss \thanks{Division of Scientific Computing, Department of Information Technology, Uppsala University, Sweden.} }
\pagenumbering{arabic}
\maketitle
\begin{abstract}
 It is well-known that reliable and efficient domain truncation is crucial to accurate numerical solution of most wave propagation problems. The perfectly matched layer (PML) is a method which, when
stable, can provide a domain truncation scheme which is convergent with increasing layer width/damping. The difficulties in using the PML are primarily associated with stability, which can be present at the continuous level or be triggered by numerical approximations.  The mathematical and numerical analysis of the PML for hyperbolic wave propagation problems has been an area of active research. It is now possible to construct stable and high order accurate numerical wave solvers by augmenting wave equations with the PML and approximating the equations  using  summation-by-parts finite difference methods, continuous and discontinuous Galerkin finite element methods. In this review we summarise the progress made, from mathematical, numerical and practical perspectives, point out some open problems  and set the stage for future work. We also present numerical experiments of model problems corroborating the theoretical analysis, and numerical simulations of real-world wave propagation demonstrating impact. Stable and parallel implementations of the PML in the high performance computing  software packages WaveQLab3D and ExaHyPE  allow  to sufficiently limit the computational domain of seismological problems with only a few grid points/elements around the computational boundaries where the PML is active, thus saving as much as $96\%$ of  the required computational resources for a three space dimensional seismological benchmark problem.
\end{abstract}

\section{Introduction}
\paragraph{Background.} 
Real-world wave propagation problems are often formulated in unbounded or very large domains. Because of limited computational resources numerical simulations are performed in truncated (smaller) computational domains by introducing artificial boundaries.
However, one of the main  features of propagating waves is that they can travel long distances relative to their characteristic length-scale, the wavelength.
As an example, a strong ground shaking resulting from an earthquake or a nuclear explosion that occurred in one continent can be recorded far away in another continent sometime after it has occurred.
This innate capability of waves to travel long distances while carrying important information of their sources and the medium of propagation  allows us to use them to probe the world around us and develop modern multi-media technologies.
In practice, computer simulations of such wave phenomena are restricted to smaller computational domains and areas of interest, for example, where the effects of the strong ground motion is significant. This is typical of numerical simulations of most real-world wave propagation problems including earthquake seismology, seismic imaging, aero-acoustics, sonar applications, wireless communication, ground penetrating radar technologies, to name only a few.

For numerical simulations, it is precisely this essential feature of waves, the radiation of waves to far field, that leads to the greatest difficulties. In truncated computational domains outgoing wave radiations will result to spurious reflections at artificial boundaries, which will travel back into the simulation domain and destroy the accuracy of numerical simulations everywhere. Therefore, in order to ensure the accuracy of numerical simulations, artificial boundaries introduced to limit the computational domain must be closed with reliable and accurate boundary conditions.

\paragraph{Absorbing boundary conditions and absorbing layers.} 
The endeavour to design efficient domain truncation strategies for waves began over forty five  years ago \cite{Engquist77}, and has yielded two powerful and equally appealing  approaches for realising effective domain truncation schemes.
The first approach is the absorbing boundary condition (ABC) \cite{Engquist77,ALPERT2002270,10.2307/2008166,https://doi.org/10.1002/cpa.3160320303,collino1993high,LysmerKuhlemeyer1969,Baffe2012,Hagstrom2009,Hagstrom2003,GIVOLI2004319,hagstrom2010radiation,Baffe2012,hagstrom_1999,https://doi.org/10.1002/cpa.3160330603,doi:10.1137/S0036139993269266,RABINOVICH2017190,10.1785/BSSA0670061529}, which is a boundary condition  enforced at an artificial boundary such that unwanted reflections from the boundary are significantly minimised. The second approach corresponds to extending the domain to an absorbing layer of finite thickness where the underlying equations are transformed such that waves decay rapidly in the layer \cite{LAVELLE2008270,APPELO20094200,berenger1994,Be_etAl,SkAdCr,AppeloKreiss2006,KDuru2016,DuKr,petersson_sjogreen_2014,KREISS20161,Renaut_1996,grote2010efficient,ABARBANEL1999266,kaltenbacher2013modified,doi:10.1190/segam2014-0860.1,doi:10.1080/10618560410001673524}. In order for an absorbing layer to be effective the equations must be perfectly matched \cite{Duruthesis2012,HalpernPetit-BergezRauch2011,berenger1994,Be_etAl}.
The perfectly matched layer (PML), \cite{berenger1994,Be_etAl,SkAdCr,AppeloKreiss2006,KDuru2016,DuKr,Duruthesis2012,doi:10.1137/S1064827596301406}, are absorbing layers with the desirable properties that all waves enter the layer and are absorbed without any reflections, regardless of frequency and angle of incidence. The PML  was initially introduced over twenty five years ago for electromagnetic waves \cite{berenger1994,Roden_and_Gedney_2000,Chew1994,Kuzuoglu96}.  However, the approach has since then been extended to other wave propagation problems such as elastic wave propagation, see for example \cite{AppeloKreiss2006,ElasticDG_PML2019,DuKrSIAM,Pled2021,Komatitsch10.1046/j.1365-246X.2003.01950.x,Duruthesis2012,cmes.2010.056.017,cmes.2008.037.274}, and acoustic wave propagation, for examples \cite{grote2010efficient,ABARBANEL1999266,kaltenbacher2013modified, Duru2012JSC, DURU2014757}. The popularity of the  PML is mainly due to its effective absorption properties, versatility, simplicity, ease of derivation and implementation using standard numerical methods.  For example, using the so-called complex coordinate stretching technique \cite{Chew1994, ZhaoCangellaris491508, PETROPOULOS1998184} the derivation of the PML for many linear hyperbolic systems is simplified and takes the following three straightforward steps.
Consider the Cartesian spatial coordinates $(x,y,z)$ and a PML in $\xi$-direction with $\xi \in \{x, y, z\}$.
The steps for a standard PML model are as follows:\\
$\bullet$ Take the Laplace transform in time of the underlying equations, $\partial /\partial t \to s$, \\
$\bullet$  Transform the partial derivatives using the PML complex metric $S_\xi = 1 + d_\xi(\xi)/s$,  $\partial/{\partial \xi}\to \left(1/S_\xi\right) \partial/{\partial \xi}$,  \\
$\bullet$  Introduce local auxiliary variables and  invert the Laplace transforms.\\
Here, $s$ is the dual time variable and $d_\xi\ge 0$ is the absorption function. Note that the PML metric $S_\xi = 1 + d_\xi(\xi)/s$ is standard, and can be enhanced by including more parameters \cite{Chew1994, ZhaoCangellaris491508, PETROPOULOS1998184, Duruthesis2012, doi:10.1137/S1064827596301406,Danieldoi:10.1137/050639107} (see also Section 3 of this review) such that the PML is more robust.
The PML transformation  has some important mathematical and numerical consequences which are critical for effective use in many practical computations. For the past twenty years,  the mathematical analysis of well-posedness and stability, and numerical analysis of error and convergence of the PML has  attracted substantial attention in the literature, see for examples \cite{Duruthesis2012,Be_etAl,SkAdCr,AppeloKreiss2006,DuruGabKreiss2019,HalpernPetit-BergezRauch2011,Becache1296843,becache_joly_2002,Danieldoi:10.1137/050639107,ABARBANEL1998331,ABARBANEL1997357,HESTHAVEN1998129,HU1996201,HU2001455,DURU201434,TAM1998213,Lions2002,ABARBANEL1999266,Baffet2019,BecaheKachanovska2021,SkAdCr,doi:10.1137/110835268,doi:10.1137/040604315,refId0,BECACHE20121639,Becache2004}. 


\paragraph{Aim of the paper.} 
In this paper, we shall review the state-of-the-art results for the PML, both when it comes to theoretical results and practical issues. We have collected some useful results concerning how to discretely handle PML boundaries and interfaces for several types of wave propagation problems. We will consider PML initial boundary value problems (IBVPs) for first order hyperbolic systems and mention some extensions to second order hyperbolic systems, on domains partially or completely surrounded by a PML. Throughout the review we shall highlight important results and point out some open problems.

Mathematical well-posedness and stability are fundamental features of any useful model for accurate numerical absorption of waves.  
Analysis of well-posedness and stability of the PML has  been of significant interests  \cite{Duruthesis2012,Be_etAl,SkAdCr,AppeloKreiss2006,DuruGabKreiss2019,HalpernPetit-BergezRauch2011,Becache1296843,becache_joly_2002,ABARBANEL1998331,Sun_Zhang_Zhan_Liu_2018,ModaveLambrechtsGeuzaine2016,ABARBANEL1997357,HESTHAVEN1998129,HU1996201,HU2001455,doi:10.1137/040604315,refId0,BECACHE20121639,Becache2004,refId0,BecacheAndKachanovska}. We will try to summarise the results without going into too much details of the proofs. For initial value problems (IVPs) the stability of the PML can be predicted by the so--called geometric stability condition \cite{Be_etAl} using classical Fourier analysis. This is corroborated by the analysis in \cite{Danieldoi:10.1137/050639107} by deriving decaying energy densities in the Fourier space and in \cite{HalpernPetit-BergezRauch2011} using geometric optics, and verified through several numerical experiments published in the literature. 
We note that for some simple models such as the acoustic wave equation with constant damping and constant medium parameters, analytical solution of the PML can be derived,  for instance using the  Cagniard-De  Hoop  method \cite{BecaheKachanovska2021,DIAZ20063820}.
However, several numerical experiments presented in the literature suggest that, even when the corresponding IVP is stable,  the PML can suffer from exponential energy growth, in particular   when the PML interacts with certain boundary conditions. Furthermore, even when there are no physical boundaries, in practice the PML must be implemented in a bounded computational domain, as an IBVP.  Therefore we must prove well-posedness and stability  for PML IBVPs and derive stable and convergent numerical approximations of the PML IBVPs.

The theory and numerical methods to solve IBVPs are more elaborate and complicated than those of a corresponding IVP; for a well-posed partial differential equation (PDE) which is stable in the absence of boundaries can support unstable solutions, or become ill-posed, when boundaries are introduced \cite{Duruthesis2012,GustafssonKreissOliger1995,SkAdCr,KDuru2016,doi:10.1137/1.9780898719130,DuKrSIAM}. 
Indeed, the analysis of the IVP is a necessary (first) step towards the analysis of the IBVP \cite{GustafssonKreissOliger1995,KDuru2016,doi:10.1137/1.9780898719130}.
In numerical approximations of IBVPs, most of the difficulties arise from the boundaries. Similarly, a numerical method which is stable in the absence of boundaries can support growth when boundary conditions are imposed \cite{GustafssonKreissOliger1995,BGustafsson,doi:10.1137/060654943,BGustafsson,BGustafsson98}.

We note that the mode analysis for PMLs has been extended to PML IBVPs, to investigate the effects of boundary conditions on the well-posedness and stability of PML IBVPs, see \cite{Duruthesis2012,DuKrSIAM,KDuru2016,DuruKozdonKreiss2016}.
The overarching result is that the PML IBVP will not support exponentially growing solutions as long as neither the underlying undamped IBVP nor the PML IVPs support exponentially growing solutions.   However these results are  { too technical to be extended to the analysis of numerical approximations}. An approach that has been demonstrated to provide a cure for this problem is the extension of the theoretical stability results using the energy method in the Laplace space \cite{KDuru2016,DuruKozdonKreiss2016,ElasticDG_PML2019,DuruGabKreiss2019}. 

One aim of this review is to demonstrate the importance and strength of the theoretical results in the form of energy estimates. 
We will discuss the results in both continuous and discrete settings, and demonstrate that the discrete energy estimates are essential for the usefulness of PMLs. We will focus here on summation-by-part (SBP) finite difference methods and DG finite element methods, where there are well developed techniques to impose boundary and interface conditions weakly with penalty terms. Such a technique allows for mimicking continuous energy estimates. 

For IBVPs, a standard approach to achieve numerical stability is to first derive an energy estimate for the continuous problem. The next step is to derive a discrete energy estimate, which typically is done by mimicking the procedure in the continuous setting. If such an estimate is possible, numerical stability is ensured. For PMLs this approach has not worked so well. The main issue is that the PML is generally asymmetric. It becomes extremely difficult to derive energy estimates that can be used to design stable numerical methods for the PML in truncated domains.  In many settings, a straight forward approach would yield an exponentially growing estimate in physical space. This would indicate well-posedness, but exponential growth is not optimal for an absorbing model. A well-functioning PML should not support exponential growth. We note that at constant coefficient an energy estimate for the acoustic wave equation in second order form was recently derived  for a specific time-domain formulation of the PML \cite{doi:10.1137/110835268,Baffet2019}. However, we will review the use of energy methods in the Laplace space \cite{KDuru2016,DuruKozdonKreiss2016,ElasticDG_PML2019,DuruGabKreiss2019}, in general media, which is applicable to most PML models and useful for developing  stable and high order accurate numerical methods.

Most modern numerical methods such as the discontinuous Galerkin (DG) finite/spectral element methods and multi-block SBP finite difference methods require interface conditions such that locally adjacent elements/blocks can be coupled together. Interface conditions connecting elements together are often implemented through numerical fluxes. For symmetric or symmetrisable linear hyperlolic systems, numerical fluxes are often designed such that numerical scheme obeys a discrete energy estimate, thus ensuring robustness of numerical approximations. However, when the PML is present, the development of  accurate and stable numerical fluxes has proven to be a nightmare for practitioners. Exponential and/or linear growth is often seen in numerical simulations using modern numerical methods.
As above, numerical procedures  based on the energy method in the Laplace space has proven to yield robust numerical methods and provide a possible cure for this problem \cite{KDuru2016,DuruKozdonKreiss2016,ElasticDG_PML2019,DuruGabKreiss2019}.  We will  summarise the  energy-based method in the Laplace  \cite{ElasticDG_PML2019,DuruGabKreiss2019} for stable DG implementations of the PML in acoustics and linear elastodynamics.

Another aim is to demonstrate impact on practical real-world simulations.
Our algorithms and the PML have been implemented in two different freely open source  high performance computing (HPC) software packages, WaveQLab3D \cite{DuruandDunham2016,DuruFungWilliams2020} and ExaHyPE \cite{ElasticDG_PML2019,Duru_exhype_2_2019}, for large-scale simulation of seismic waves in geometrically complex 3D Earth models. The software package WaveQLab3D \cite{DuruandDunham2016,DuruFungWilliams2020} is  a  high order accurate SBP-SAT finite difference solver. ExaHyPE is a DG solver of arbitrary accuracy for large-scale numerical simulation of hyperbolic wave propagation problems on dynamically adaptive curvilinear meshes. Stable and parallel implementations of the PML in the HPC software packages WaveQLab3D and ExaHyPE \cite{DuruandDunham2016,DuruFungWilliams2020,DuruGabKreiss2019,ElasticDG_PML2019,Duru_exhype_2_2019} allow  us to sufficiently limit the computational domain of seismological problems with only a few grid points/elements around the computational boundaries where the PML is active. Thus saving as much as $96\%$ of  the required computational resources for a 3D seismological benchmark problem \cite{Seismowine}. We also present  a real-world wave propagation propagation problem which involves the simulation of 3D seismic waves in a section of European Alpine region, with strong non-planar free-surface topography \cite{Copernicus,Duru_exhype_2_2019}.

\paragraph{Structure of the paper}
The remaining parts of the review will proceed as follows. In section 2 we introduce a general model of linear hyperbolic systems in first and second order forms. We perform dispersion relation analysis, introduce the {\it geometric stability condition} and formulate well-posed boundary and interface conditions. In section 3 we review the derivation of the PML using complex coordinate stretching. Mathematical analysis of the PML at the continuous level is reviewed in section 4. In section 5 we summarise theoretical results on the numerical analysis of the PML discrestised with the DG finite element method and SBP finite difference methods, respectively. In section 6 we present some numerical experiments verifying the theoretical results and demonstrating the importance of the continuous and discrete energy estimates. In section 7 we present the use  of PML in practical seismological applications problems, and demonstrate impact. We draw conclusions in section 8 and speculate on the directions for future work.

\section{Linear hyperbolic  partial differential equations}
\label{sec:s2}
In this section, we present model equations that are representative of  linear hyperbolic  wave propagation problems encountered in different application areas, such as acoustics, seismology, electrodynamics and elastodynamics.
We will consider the equations in both first order form and second order formulations.
Plane waves will be used to discuss the wave propagation properties of the medium at constant coefficients.
We will also introduce boundary conditions, comment on well-posedness and derived energy estimates.
\subsection{First order systems}
Consider the linear first order hyperbolic partial differential equation (PDE)  in Cartesian coordinates and a source free heterogeneous medium:
{
\normalsize
 \begin{equation}\label{eq:first_order_linear_hyp_pde}
\begin{split}
\mathbf{P}^{-1} \frac{\partial{\mathbf{U}}}{\partial t} = \sum_{\xi = x, y, z}\mathbf{A}_{\xi}\frac{\partial{\mathbf{U}}}{\partial \xi}, \quad (x, y, z) \in \Omega \subset \mathbb{R}^d,\quad d = 3, \quad t \ge  0,
\end{split}
\end{equation}
}
subject to the initial condition
\begin{align*}
 \mathbf{U}(x,y,z,0) = \mathbf{U}_0(x,y,z) \in \mathrm{L}^2\left(\Omega\right).
\end{align*}
Here, $\mathbf{U} \in \mathbb{R}^m$, $m \ge 1$ is the unknown vector field, the coefficient matrices are symmetric, $\mathbf{A}_{\xi} = \mathbf{A}_{\xi}^T\in \mathbb{R}^{m\times m}$ and $\mathbf{P}= \mathbf{P}^T\in \mathbb{R}^{m\times m}$,  with $ \mathbf{U}^T\mathbf{P}\mathbf{U} > 0$. In general the matrix $\mathbf{P}$ depends on the spatial coordinates $x, y,z$, and encodes the material parameters of the underlying medium. The constant coefficients and non-dimensional  matrices $\mathbf{A}_{\xi}$ encapsulate the underlying linear conservation law and the corresponding linear constitutive relation. 

As will be shown below, depending on the coefficient matrices $\mathbf{P}, \mathbf{A}_{\xi}$ the  system \eqref{eq:first_order_linear_hyp_pde} describes first order linear hyperbolic wave propagation problems encountered in different application areas such acoustics, elastodynamics, electromagnetics, etc. 
For example in an acoustic medium  we denote $p$  the pressure and $\mathbf{v}= (v_x,v_y,v_z)^T$  the particle velocity, $\rho$ the density, and $\kappa  > 0$  the bulk modulus of the medium, we have 
{ 
\footnotesize
\begin{align*}
 \mathbf{U} = \begin{bmatrix}
p\\
\mathbf{v}
 \end{bmatrix},
 \quad 
\mathbf{P} = 
\begin{pmatrix}
\kappa  & \mathbf{0}^T\\
\mathbf{0}  & \rho^{-1} \mathbf{I}_3
\end{pmatrix},
\quad
\mathbf{A}_{\xi} = 
\begin{pmatrix}
{0} & -\mathbf{e}_\xi^T\\
-\mathbf{e}_\xi & \mathbf{0}_3
\end{pmatrix},
  \quad
\mathbf{0} = \left(0, 0, 0\right)^T,
\quad
\mathbf{e}_{\xi} = \left({e}_{\xi x}, {e}_{\xi y}, {e}_{\xi z}\right)^T, 
\end{align*}
where $\mathbf{I}_3$ is a 3-by-3 identity matrix,  $\mathbf{0}_3$ is a 3-by-3 zero matrix, and 
}
\begin{equation*}
\begin{split}
{e}_{\xi\eta} = \left \{
\begin{array}{rl}
0 \quad {}  \quad \text{if} \quad \xi \ne \eta,\\
1 \quad {}  \quad  \text{if} \quad \xi = \eta.
\end{array} \right.
\end{split}
\end{equation*}

To describe wave propagation in elastic solids, we introduce the unknown wave fields 
\begin{align}\label{eq:velocity_stress}
\mathbf{U}\left(x,y,z,t\right) = \begin{bmatrix}
\mathbf{v}(x,y, z,t) \\
\boldsymbol{\sigma}(x,y,z,t)
 \end{bmatrix},
\end{align}
with the particle velocity vector is $\mathbf{v} = \left[ v_x, v_y, v_z \right]^T$, and the stress vector is given by\\ $\boldsymbol{\sigma} = \left[ \sigma_{xx}, \sigma_{yy}, \sigma_{zz}, \sigma_{xy},  \sigma_{xz},  \sigma_{yz}\right]^T$.  
The symmetric constant coefficient matrices $\mathbf{A}_\xi $ describing the conservation of momentum and the constitutive relation, defined by Hooke's law, are given by
%
{
\begin{align*}
&\mathbf{A}_{\xi} = 
\begin{pmatrix}
\mathbf{0}_3 & \mathbf{a}_{\xi}\\
\mathbf{a}_{\xi}^T & \mathbf{0}_6
\end{pmatrix},
\quad
\mathbf{a}_x = 
\begin{pmatrix}
1& 0& 0& 0&0& 0\\
0& 0& 0& 1&0& 0\\
0& 0& 0& 0&1& 0
\end{pmatrix},
\\
&\mathbf{a}_y = 
\begin{pmatrix}
0& 0& 0& 1&0& 0\\
0& 1& 0& 0&0& 0\\
0& 0& 0& 0&0& 1
\end{pmatrix},
\quad
\mathbf{a}_z = 
\begin{pmatrix}
0& 0& 0& 0&1& 0\\
0& 0& 0& 0&0& 1\\
0& 0& 1& 0&0& 0\\
\end{pmatrix},
\end{align*}
}
where $\mathbf{0}_3$ and $\mathbf{0}_6$ are  the $3$-by-$3$ and $6$-by-$6$ zero  matrices.
The symmetric positive definite material parameter matrix $\mathbf{P}$  is defined by
{
\footnotesize
\begin{align}\label{eq:material_coeff}
\mathbf{P} = 
\begin{pmatrix}
\rho^{-1} \mathbf{I}_3  & \mathbf{0}\\
 \mathbf{0}^T  & \mathbf{C}
\end{pmatrix}
  ,
  \quad
  \mathbf{0} =   \begin{pmatrix}
  0 & 0 & 0& 0 & 0 & 0 \\
  0 & 0 & 0& 0 & 0 & 0 \\
  0 & 0 & 0& 0 & 0 & 0 
  \end{pmatrix}
  ,
\end{align}
}
where  $\rho(x,y,z) > 0$ is the  density of the medium, and  $\mathbf{C} = \mathbf{C}^{T} > 0$ is the symmetric positive definite  matrix of elastic constants.

For real functions, we introduce the weighted $L^2$ inner product and the corresponding energy norm
\begin{align}\label{eq:weighted_scalar_product}
\left(\mathbf{U}, \mathbf{F}\right)_P = \int_{\Omega}{\frac{1}{2}[\mathbf{U}^T\mathbf{P}^{-1}\mathbf{F}] dxdydz}, \quad \|\mathbf{U}\left(\cdot, \cdot, \cdot, t\right)\|_P^2 = \left(\mathbf{U}, \mathbf{U}\right)_P 
\end{align}
%
The weighted $L^2$-norm $\|\mathbf{U}\left(\cdot, \cdot, \cdot, t\right)\|_P^2$ is the physical energy in the medium.

For the IVP \eqref{eq:first_order_linear_hyp_pde}, $\Omega = \mathbb{R}^3$, the energy (conservation) equation follows
\begin{align}\label{eq:cauchy_estimate_first}
\frac{d}{dt}\|\mathbf{U}\left(\cdot, \cdot, \cdot, t\right)\|_P^2 = 0.
\end{align}
The conservation of energy \eqref{eq:cauchy_estimate_first} indicates that the IVP \eqref{eq:first_order_linear_hyp_pde} is  well-posed and strongly stable. 

%
%
%
\subsection{Second order systems}
Now consider the linear second order hyperbolic PDE in Cartesian coordinates in a source free heterogeneous medium
\begin{equation}\label{eq:second_order_linear_hyp_pde}
\begin{split}
\boldsymbol{\Phi}^{-1} \frac{\partial^2{\mathbf{u}}}{\partial t^2} = \sum_{\xi = x, y, z}\frac{\partial }{\partial \xi}\left(\sum_{\eta = x, y, z}\mathbf{A}_{\xi \eta}\frac{\partial{\mathbf{u}}}{\partial \eta}\right), \quad (x, y, z) \in \Omega \subset \mathbb{R}^3, \quad t \ge  0,
\end{split}
\end{equation}
%
subject to the initial condition
\begin{align*}
\mathbf{u}(x,y,z,0) = \mathbf{u}_0(x,y,z)\in  \mathrm{H}^1\left(\Omega\right), \quad  \frac{\partial }{\partial t}\mathbf{u}(x,y,z,0) = \mathbf{v}_0(x,y,z)\in \mathrm{L}^2\left(\Omega\right).
\end{align*}
Here, $ \mathbf{u}(x,y,z,t) \in \mathbb{R}^{n}$ is the unknown vector field, $\boldsymbol{\Phi} \in \mathbb{R}^{n\times n}$ is a symmetric and positive definite matrix, and the matrices $\mathbf{A}_{\xi\eta} \in \mathbb{R}^{n\times n}$ are such that the potential energy matrix $\mathcal{P} \in \mathbb{R}^{dn\times dn}$ defined below is symmetric semi-positive definite, that is
$$
\mathcal{P} 
=
\begin{pmatrix}
\mathbf{A}_{xx} & \mathbf{A}_{xy} & \mathbf{A}_{x z}\\
\mathbf{A}_{xy}^T & \mathbf{A}_{yy} & \mathbf{A}_{y z}\\
\mathbf{A}_{xz}^T & \mathbf{A}_{yz}^T & \mathbf{A}_{z z}
\end{pmatrix}, 
\quad
\mathbf{U}^T \mathcal{P}\mathbf{U} \ge 0, \quad \forall \mathbf{U}\in \mathbb{R}^{dn}, \quad d = 3.
$$
Many linear hyperbolic PDEs originally arise as second order form \eqref{eq:second_order_linear_hyp_pde} but can be reduced to a first order system \eqref{eq:first_order_linear_hyp_pde}, by introducing new variables, see eg. the classical acoustic wave equation and the elastic wave equation. Some models also initially appear as  first order systems \eqref{eq:first_order_linear_hyp_pde} but they can be rewritten as second order systems \eqref{eq:second_order_linear_hyp_pde}, eg linearised Euler equations of fluid dynamics, Maxwell's equation, and linearised MHD equations, to name a few. 

For example, consider an acoustic medium where $p$ denote the pressure, $ \rho, \kappa  > 0$ denote the density and the bulk modulus of the medium with the acoustic waves peed $c = \sqrt{\kappa/\rho}$, we have 
\normalsize
\begin{align*}
\mathbf{u} = p,
\quad 
\boldsymbol{\Phi} = {\kappa}  ,
\quad
\mathbf{A}_{\xi\eta} = \left \{
\begin{array}{rl}
0 \quad {}  \quad \text{if} \quad \xi \ne \eta,\\
\frac{1}{\rho}\quad {}  \quad  \text{if} \quad \xi = \eta.
\end{array} \right.
\end{align*}
Define the energy in the medium
\begin{align}\label{eq:energy_second}
E(t) =  \frac{1}{2}\int_{\Omega}{\left(\frac{\partial{\mathbf{u}}}{\partial t}^T \boldsymbol{\Phi}^{-1} \frac{\partial{\mathbf{u}}}{\partial t} + 
\begin{pmatrix}
\frac{\partial{\mathbf{u}}}{\partial x} \\
\frac{\partial{\mathbf{u}}}{\partial y}\\
\frac{\partial{\mathbf{u}}}{\partial z} 
\end{pmatrix}^T
\mathcal{P}
\begin{pmatrix}
\frac{\partial{\mathbf{u}}}{\partial x} \\
\frac{\partial{\mathbf{u}}}{\partial y}\\
\frac{\partial{\mathbf{u}}}{\partial z} 
\end{pmatrix}\right)}dxdydz > 0.
\end{align}
The energy $E(t)$ given in \eqref{eq:energy_second} defines a semi-norm.
As above, for Cauchy problems, $\Omega = \mathbb{R}^3$, the energy equation follows
{
\begin{align}\label{eq:cauchy_estimate_second}
\frac{d}{dt} E(t) = 0.
\end{align}
}
Thus the energy is conserved, $E(t)  = E(0) $, for all $t\ge 0$. 
Again, this analysis indicates that the Cauchy problem for \eqref{eq:second_order_linear_hyp_pde} is  well-posed and asymptotically stable. 

\subsection{Dispersion relation}
To understand the wave propagation properties of the models \eqref{eq:first_order_linear_hyp_pde} and \eqref{eq:second_order_linear_hyp_pde}  it is useful to
consider wave-like solutions
\begin{align}\label{eq:plane_wave}
\mathbf{U}\left(x, y, z, t\right) = \mathbf{U}_0 e^{st + i\left( k_x x  + k_y y + k_z z \right)}, \quad \left(k_x, k_y, k_z\right) \in \mathbb{R}^3.
\end{align}
In \eqref{eq:plane_wave}, $\mathbf{k}=\left(k_x, k_y, k_z\right) $ is the wave vector, and $ \mathbf{U}_0$ is a vector of constant amplitude
called the polarization vector.  By inserting \eqref{eq:plane_wave} in \eqref{eq:first_order_linear_hyp_pde}  we have the solvability condition, often called the dispersion relation
{
\begin{align}\label{eq:dispersion_relation_f}
&F\left(s,  \mathbf{k}\right):= \Big|{P}^{-1} s -  \sum_{\xi = x,y, z}{ik_{\xi}}\mathbf{A}_{\xi}\Big| = 0.
\end{align}
}
The indeterminate $s$ is related to the temporal frequency and will be determined by the dispersion relation \eqref{eq:dispersion_relation_f}.
Since \eqref{eq:first_order_linear_hyp_pde}  conserve energy the roots $s$ must be zero or purely imaginary,   that is $s\in \mathbb{C}$ with $\Re{s} =0$.
Otherwise if the roots $s$ have non-zero real parts the energy will grow or decay, contradicting the energy equations \eqref{eq:cauchy_estimate_first}.

Similarly, inserting the simple wave solution \eqref{eq:plane_wave} in the second order system \eqref{eq:second_order_linear_hyp_pde}, we have
{
\begin{align*}
G\left(s,  \mathbf{k}\right):= \Big|\boldsymbol{\Phi}^{-1} s^2 +  \sum_{\xi = x,y, z}\sum_{\eta = x, y,z}{k_{\xi}}k_\eta\mathbf{A}_{\xi,\eta}\Big| = 0.
\end{align*}
}
If  \eqref{eq:first_order_linear_hyp_pde} and \eqref{eq:second_order_linear_hyp_pde} model the same physical phenomena, for some natural number  $r\in \mathbb{N}$, we have 
\begin{align*}
F\left(s,  \mathbf{k}\right) = s^{r}G\left(s,  \mathbf{k}\right), \quad {r} \ge 0,
\end{align*}
where
\begin{align*}
s  =\{s_0\in \mathbb{C}: \quad G\left(s_0,  \mathbf{k}\right)=0, \quad |s_0| \ne 0\}.
\end{align*}
Note that $F\left(s,  \mathbf{k}\right)$ will have $r$ zero roots $s =0$ and $m-r$ nonzero roots $s =\{s_0\}$.
The nonzero roots $s$ of $G\left(s,  \mathbf{k}\right)$
correspond to propagating physical modes. 
Specifically in acoustics with the wave speed $c = \sqrt{\kappa/\rho}$ we have
$$
G\left(s,  \mathbf{k}\right) = \frac{s^2}{c^2} + |\mathbf{k}|^2, \quad F\left(s,  \mathbf{k}\right) = s^2G\left(s,  \mathbf{k}\right),
$$
where  $
|\mathbf{k}| = \sqrt{k_x^2+k_y^2 +k_z^2}$ and the nonzero roots are $s = \pm i c |\mathbf{k}|$.
In isotropic linear elastic medium we have
$$
G\left(s,  \mathbf{k}\right) = \left(\frac{s^2}{c_s^2} + |\mathbf{k}|^2\right) \left(\frac{s^2}{c_s^2} + |\mathbf{k}|^2\right)\left(\frac{s^2}{c_p^2} + |\mathbf{k}|^2\right), \quad F\left(s,  \mathbf{k}\right) = s^3G\left(s,  \mathbf{k}\right),
$$
where $c_s$, $c_p$ are the S-wave speed and P-wave speed, respectively, of the medium, and the nonzero roots are are given by the simple roots $s = \pm i c_s |\mathbf{k}|$,  $s = \pm i c_s |\mathbf{k}|$, and $s = \pm i c_p |\mathbf{k}|$.
In anisotropic media the dispersion relation are complicated nonlinear algebraic expressions, and the roots are difficult to determine explicitly.  However, as above, one can show that the roots are purely imaginary,  that is $s\in \mathbb{C}$ with $\Re{s} =0$.

We write $s = i\omega$ where $\omega \in \mathbb{R}$ is called the temporal frequency, and introduce
\begin{align}
&\mathbf{K} = \left(\frac{k_x}{|\mathbf{k}|}, \frac{k_y}{|\mathbf{k}|}, \frac{k_z}{|\mathbf{k}|}\right), \quad \text{normalized propagation direction},\\
&\mathbf{V}_p = \left(\frac{\omega}{k_x}, \frac{\omega}{k_y}, \frac{\omega}{k_z}\right), \quad \text{phase velocity},\\
&\mathbf{S} = \left(\frac{k_x}{\omega}, \frac{k_y}{\omega}, \frac{k_z}{\omega}\right), \quad \text{slowness vector},\\
&\mathbf{V}_g = \left(\frac{\partial \omega}{\partial k_x}, \frac{\partial \omega}{\partial k_y}, \frac{\partial \omega}{\partial k_z}\right), \quad \text{group velocity}.
\end{align}

For a constant coefficient medium with no boundaries, or with periodic boundary
conditions, the dispersion relation $F\left(i\omega,  k_x, k_y, k_z\right)=0$   and the quantities $\mathbf{K}$, $\mathbf{V}_p $, $\mathbf{S}$, $\mathbf{V}_g$, defined
above, provide a rather comprehensive description of the wave propagation properties of the medium.
When physical boundaries are introduced the additional presence of boundary wave modes, such as Rayleigh waves, introduces more complex and interesting wave phenomena.
More general media, such as complex elastic media, can also support guided waves, such as Love waves, Stoneley waves, Scholte waves, conical waves, and other interface waves. 
\begin{remark}
We note however that the analysis of the PML in the presence of interface waves is still an open problem and will not be considered in this review. 
\end{remark}

Since \eqref{eq:first_order_linear_hyp_pde}  contains only first derivative terms and no lower order terms then  $G\left(s,  \mathbf{k}\right)$ is homogeneous in $\omega, k_x,k_y, k_z$, and we can rewrite $G\left(s,  \mathbf{k}\right)=0$  as $G\left(i,  \mathbf{S}\right) = 0$.

\begin{definition}[Geometric stability condition]\label{def:geometric_condition}
At any point on the slowness surface $\mathbf{S}=(S_x, S_y, S_z)$ with $G\left(i,  \mathbf{S}\right) = 0$ , the geometric stability condition in the $\xi$-direction (for $\xi = x, y, z$) is defined by
$$
{V}_{p \xi} {V}_{g \xi} \ge 0.
$$
\end{definition}
As was shown in \cite{Duruthesis2012,Be_etAl,AppeloKreiss2006} the Geometric stability condition is necessary for the stability of the PML IVP. This will be discussed in more detail in Section \ref{sec:math_analysis} where the mathematical analysis of the PML is reviewed.

\subsection{Boundary conditions}
Effective treatments of boundary conditions are critical for the  stability and efficacy of the PML \cite{Duruthesis2012,DuKrSIAM,KDuru2016,DuruKozdonKreiss2016,ElasticDG_PML2019,DuruGabKreiss2019}. 
In this section, we formulate well-posed and energy stable boundary  conditions for the first order system \eqref{eq:first_order_linear_hyp_pde} and the second order system \eqref{eq:second_order_linear_hyp_pde}. The boundary conditions will be extended when the PML is introduced. 

Consider the 3D cuboidal domain
{
\begin{equation}\label{eq:physical_domain}
\Omega = \{(x,y,z): -1 \le x \le 1, \quad -1 \le y\le 1, \quad -1 \le z\le 1\}.
\end{equation}
}
 Stable and well-posed boundary conditions are needed to close the rectangular surfaces  of the boundaries of the cuboidal domain.
Define the reference boundary surface
 \begin{align}
 \widetilde{\Gamma} = [-1, 1]\times[-1, 1].
 \end{align}
\subsubsection{Boundary condition for first order systems}
The application of the energy method to \eqref{eq:first_order_linear_hyp_pde}
 yields 
\begin{align}\label{eq:cauchy_estimate_first_bt}
\frac{d}{dt}\|\mathbf{U}\left(\cdot, \cdot, \cdot, t\right)\|_P^2 = \mathrm{BT}\left(\mathbf{U}\right),
\end{align}
where the boundary term is given by
\begin{align}\label{eq:product_BT_first_order}
  \mathrm{BT}\left(\mathbf{U}\right) &\equiv \frac{1}{2}\oint_{\partial \Omega} \left[\mathbf{U}^T \left(\sum_{\xi = x, y, z}n_{\xi}\mathbf{A}_{\xi}\right)\mathbf{U}\right] dS \nonumber \\
   &= \sum_{\xi = x, y, z}\int_{\widetilde{\Gamma}}  \frac{1}{2}\left[ \mathbf{U}^T\mathbf{A}_{\xi}\mathbf{U}  \Big|_{\xi =1}- \mathbf{U}^T\mathbf{A}_{\xi}\mathbf{U}  \Big|_{\xi =-1}  \right]  \frac{dxdydz}{d\xi}.
\end{align}
For the acoustic wave equation we have
$$
\frac{1}{2}\mathbf{U}^T\mathbf{A}_{\xi}\mathbf{U} = -v_\xi p.
$$
For linear elastodynamics, at the boundaries $\xi = \pm 1$ we introduce with the tractions and  particle velocity  on the boundary, 
$$\mathbf{T}: =\mathbf{a}_\xi\boldsymbol{\sigma}= (T_x, T_y, T_z)^T,
 \quad \mathbf{v} = (v_x, v_y, v_z)^T.
$$
 We have
$$
\frac{1}{2}\mathbf{U}^T\mathbf{A}_{\xi}\mathbf{U} = \mathbf{v}^T\mathbf{T} = \sum_{\xi = x, y, z}v_\xi T_\xi.
$$
See \cite{Duru_exhype_2_2019,ElasticDG_PML2019,DuruGabKreiss2019} for details.

Well-posed boundary conditions are imposed such that the boundary term   $\mathrm{BT}\left(\mathbf{U}\right)$ defined in \eqref{eq:product_BT_first_order} is negative semi-definite. 
We consider linear boundary operators 
\begin{align}\label{eq:bc_first_order}
\mathbf{L}\left(\mathbf{U}\right) = 0, \quad \xi = \pm 1,
\end{align}
such that 
$$
\mathbf{U}^T \mathbf{A}_\xi \mathbf{U} \ge 0, \quad \text{at $\xi = -1$},
\quad
\mathbf{U}^T \mathbf{A}_\xi \mathbf{U} \le 0 , \quad \text{at $\xi =  1$}.
$$
Such boundary conditions are formulated  by prescribing the amplitude of the incoming characteristic variables, see  \cite{GustafssonKreissOliger1995,DuruGabrielIgel2017,Duru_exhype_2_2019,ElasticDG_PML2019} for more details. 
Thus the number of boundary conditions must be the same as the number of ingoing characteristics at the boundary. 

Specifically for acoustics we introduce the impedance $Z = \rho c$, the reflection coefficients  $r_{\xi}$, with $|r_{\xi}|\le 1$ and the well-posed boundary condition
{
\begin{equation}\label{eq:boundary_condition_acoustics}
\mathbf{L}\left(\mathbf{U}\right):=\frac{1-r_{\xi}}{2}Zv_{\xi} \mp \frac{1+r_{\xi}}{2}p = 0, \quad |r_{\xi}|\le 1, \quad \text{at} \quad \xi = \pm 1.
 \end{equation}
}
For linear elasticity, on a given boundary $\xi = \pm 1$, we introduce  the impedance $Z_\eta >0$ (for  $\eta=x,y,z$) and the well-posed boundary conditions
{
\begin{equation}\label{eq:boundary_condition_elastic}
[\mathbf{L}\left(\mathbf{U}\right)]_{\eta}:=\frac{1-r_{\xi}}{2}Z_{\eta}v_{\eta} \pm \frac{1+r_{\xi}}{2}T_{\eta} = 0, \quad |r_{\xi}|\le 1, \quad \text{at} \quad \xi = \pm 1.
 \end{equation}
}
The reflection coefficient $r_{\xi} = -1$ will correspond to boundary conditions for velocities, that is, soft-wall boundary conditions ($v_{\xi} =0$) for acoustics  and clamped-wall boundary conditions ($v_{\eta} =0$) for elastodynamics. With  $r_{\xi} = 1$ we will have boundary conditions for pressure and traction, that is, hard-wall boundary conditions ($p =0$) for acoustics  and free-surface boundary conditions ($\mathbf{T} = 0$) for elastodynamics. The zero reflection coefficients, $r_{\xi} = 0$, yield the classical first order ABC
\cite{Engquist77,LysmerKuhlemeyer1969,Baffe2012,Hagstrom2009,Hagstrom2003,GIVOLI2004319,hagstrom2010radiation,Baffe2012,hagstrom_1999,10.1785/BSSA0670061529}.

The boundary conditions \eqref{eq:boundary_condition_acoustics} and \eqref{eq:boundary_condition_elastic}  with  $|r_{\xi}|\le 1$
satisfy $\mathrm{BT}\left(\mathbf{U}\right) \le 0$, and ensure bounded solutions through
\begin{align}\label{eq:cauchy_estimate_first_0}
\frac{d}{dt}\|\mathbf{U}\left(\cdot, \cdot, \cdot, t\right)\|_P^2 = \mathrm{BT}\left(\mathbf{U}\right) \le 0.
\end{align}
This energy loss through the boundary is what a stable numerical method should emulate.
\subsubsection{Boundary condition for second order systems}
Consider the second order system in the bounded domain $(x,y,z)\in [-1, 1]^{3}$, the energy method gives
\begin{align*}
{ \frac{d }{ dt}{E}(t) = \mathrm{BT}\left(\mathbf{u}\right)},
\end{align*}
where the boundary term is given by
{\small
\begin{equation}\label{eq:product_BT_second_order}
\begin{split}
 &\mathrm{BT}\left(\mathbf{u}\right) \equiv \frac{1}{2}\oint_{\partial \Omega} \left[\frac{\partial{\mathbf{u}}}{\partial t}^T \left(  \sum_{\eta = x, y, z}\mathbf{A}_{\xi \eta}\frac{\partial{ \mathbf{u}}}{\partial \eta} \right)\right] dS  \\
  &= \sum_{\xi = x, y, z}\int_{\widetilde{\Gamma}}  \frac{1}{2}\left[ \frac{\partial{\mathbf{u}}}{\partial t}^T \left(  \sum_{\eta = x, y, z}\mathbf{A}_{\xi \eta}\frac{\partial{ \mathbf{u}}}{\partial \eta} \right)  \Big|_{\xi =1}- \frac{\partial{\mathbf{u}}}{\partial t}^T \left(  \sum_{\eta = x, y, z}\mathbf{A}_{\xi \eta}\frac{\partial{ \mathbf{u}}}{\partial \eta} \right) \Big|_{\xi =-1}  \right]  \frac{dxdydz}{d\xi}.
  \end{split}
\end{equation}
}
We consider linear boundary operators 
\begin{align}\label{eq:bc_second_order}
\mathbf{\mathcal{B}}\mathbf{u} = 0, 
\end{align}
such that the boundary term is never positive $ \mathrm{BT}\left(\mathbf{u}\right) \le 0$.
As above, we will consider specifically the well-posed linear operator
\begin{align}\label{eq:bc_second_order_operator}
\mathbf{\mathcal{B}} = \frac{1-r_{\xi}}{2}\mathbf{Z}\frac{\partial{}}{\partial t} \pm \frac{1+r_{\xi}}{2}\sum_{\eta = x, y, z}\mathbf{A}_{\xi \eta}\frac{\partial{}}{\partial \eta}, \quad \mathbf{Z} =\mathbf{Z} ^T\ge 0, \quad \xi = \pm 1,
\end{align}
with reflection coefficients $|r_{\xi}| \le 1$.
For the acoustic wave equation we have $\mathbf{Z} = \rho c >0$, and for linear elasticity the impedance is a $3$-by-$3$ diagonal matrix with diagonal entries
\begin{equation*}
\begin{split}
&Z_{\eta, \eta} = \left \{
\begin{array}{rl}
Z_p =  \rho c_p>0, \quad {}  \quad {}& \text{if} \quad \eta = \xi,\\
Z_s =  \rho c_s>0, \quad {}  \quad {}& \text{if} \quad \eta \ne \xi.
\end{array} \right.
\end{split}
\end{equation*}
Again, the reflection coefficient $r_{\xi} = -1$ will correspond to Dirichlet boundary conditions 
\begin{align}\label{eq:bc_second_order_operator_dirichlet}
\mathbf{u} = 0 \implies \frac{\partial{\mathbf{u}}}{\partial t} = 0, \quad \xi = \pm 1,
\end{align}
equivalent to the soft-wall boundary conditions for acoustics  and clamped-wall boundary conditions for elastodynamics. With  $r_{\xi} = 1$ we  have the Neumann boundary condition
\begin{align}\label{eq:bc_second_order_operator_neumann}
 \sum_{\eta = x, y, z}\mathbf{A}_{\xi \eta}\frac{\partial{\mathbf{u}}}{\partial \eta} = 0, \quad \xi = \pm 1,
\end{align}
equivalent to hard-wall boundary conditions for acoustics  and free-surface boundary conditions for elastodynamics. The zero reflection coefficients, $r_{\xi} = 0$, yield the classical first order ABC,
\begin{align}\label{eq:bc_second_order_operator_abc}
 \frac{1}{2}\mathbf{Z}\frac{\partial{\mathbf{u}}}{\partial t} \pm \frac{1}{2}\sum_{\eta = x, y, z}\mathbf{A}_{\xi \eta}\frac{\partial{\mathbf{u}}}{\partial \eta}=0,  \quad \xi = \pm 1.
\end{align}
Note that for $|r_{\xi}| = 1$ we have 
$$
\mathrm{BT}\left(\mathbf{u}\right) =0,
$$
and if $|r_{\xi}| < 1$
$$
 \mathrm{BT}\left(\mathbf{u}\right) =  -\sum_{\xi = x, y, z}\int_{\widetilde{\Gamma}}  \frac{1}{2}\left[ \frac{1-r_{\xi}}{1+r_{\xi}}\frac{\partial{\mathbf{u}}}{\partial t}^T \mathbf{Z}\frac{\partial{\mathbf{u}}}{\partial t}  \Big|_{\xi =1}+ \frac{1-r_{\xi}}{1+r_{\xi}}\frac{\partial{\mathbf{u}}}{\partial t}^T \mathbf{Z}\frac{\partial{\mathbf{u}}}{\partial t}  \Big|_{\xi =-1}   \right]  \frac{dxdydz}{d\xi} \le 0.
$$
We have
\begin{align}\label{eq:cauchy_estimate_second_0}
 \frac{d }{ dt}{E}(t) = \mathrm{BT}\left(\mathbf{u}\right) \le 0.
\end{align}
\subsection{Interface conditions}
Some numerical methods such as the DG and multi-block SBP finite difference methods require interface conditions such that locally adjacent elements/blocks can be coupled together.
This is also the situation in many settings where material parameters are discontinuous.
We introduce physical interface conditions  which we will use to couple local   elements to the global domain.
Later we will summarise how the interface conditions will be perturbed by the PML. 

Consider the IVP with a planar interface at $\xi = 0$ and  denote the corresponding fields and material parameters in the positive/negative sides of the interface with the superscripts $+/-$.
We define the jumps in scalar or vector valued fields by $\lJump{{a} \rJump} =  a^{+} - a^{-}$.
\subsubsection{Interface condition for the first order system}
 As before, using the energy method we have
 \begin{align}\label{eq:energy_estimate_fault_30_first}
\frac{d}{dt}\left(\|\mathbf{U}^{-}\left(\cdot, \cdot, \cdot, t\right)\|_P^2+ \|\mathbf{U}^{+}\left(\cdot, \cdot, \cdot, t\right)\|_P^2\right)   =   \mathrm{IT}\left( \mathbf{U}^T\mathbf{A}_{\xi}\mathbf{U}  \right).
 \end{align}
Here, the interface term is defined by
\begin{align}\label{eq:interface_term}
  \mathrm{IT}\left(\mathbf{a} \right) = - \frac{1}{2}\int_{\widetilde{\Gamma}} \lJump \mathbf{a}  \rJump   \frac{dxdydz}{d\xi}.
\end{align}
At interfaces the physics determines the interface conditions. In a standard model with conservation of the 
physical energy, interface conditions 
\begin{align}\label{eq:interface_condition}
\mathbb{I}\mathbb{T}\left(\mathbf{A}_{\xi}, \mathbf{U} \right) =0,
\end{align}
must guarantee that the interface term vanishes identically  $\mathrm{IT}\left( \mathbf{U}^T\mathbf{A}_{\xi}\mathbf{U}  \right) =0$.
For well-posedness, the number of interface conditions must correspond to the number of in-going characteristics at the interface.

In acoustics we prescribe the continuity of the pressure and the normal velocity, having 
\begin{equation}\label{eq:elastic_acoustic_interface_first}
 \mathbb{I}\mathbb{T}\left(\mathbf{A}_{\xi}, \mathbf{U} \right):=\{ \lJump{{p}} \rJump =0, \quad \lJump {{v}}_\xi \rJump = 0, \quad \text{at} \quad \xi =0 \}.
\end{equation}
In linear elastodynamics, we impose force balance, no slip and no opening conditions 
\begin{equation}\label{eq:elastic_elastic_interface_first}
  \mathbb{I}\mathbb{T}\left(\mathbf{A}_{\xi}, \mathbf{U} \right):=\{ \lJump{\mathbf{T}}\rJump =0, \quad \lJump {\mathbf{v}} \rJump = 0, \quad \text{at} \quad \xi =0\}.
\end{equation}
Note that the interface conditions \eqref{eq:elastic_acoustic_interface_first}--\eqref{eq:elastic_elastic_interface_first} 
ensure
$\mathrm{IT}\left( \mathbf{U}^T\mathbf{A}_{\xi}\mathbf{U}  \right) =0$
and conserve energy
 \begin{align*}
\frac{d}{dt}\left(\|\mathbf{U}^{-}\left(\cdot, \cdot, \cdot, t\right)\|_P^2+ \|\mathbf{U}^{+}\left(\cdot, \cdot, \cdot, t\right)\|_P^2\right)   =0.
 \end{align*}
An accurate and stable numerical method should as much as possible mimic this energy estimate.

\subsubsection{Interface condition for the second order system}
As before, using the energy method we have
\begin{align}\label{eq:energy_estimate_fault_30_second}
\frac{d}{dt}\left({E}^{-}(t)+ {E}^+(t)\right)   =   \mathrm{IT}\left( \frac{\partial{\mathbf{u}}}{\partial t}^T \left(  \sum_{\eta = x, y, z}\mathbf{A}_{\xi \eta}\frac{\partial{ \mathbf{u}}}{\partial \eta} \right) \right),
\end{align}
where ${E}^{-}(t), {E}^+(t)$ are the  energy in the negative and positive sub-domains, and the interface term $\mathrm{IT}$ is the surface integral defined by \eqref{eq:interface_term}.
Here the the well-posed interface condition is
\begin{equation}\label{eq:elastic_acoustic_interface_second}
\lJump \sum_{\eta = x, y, z}\mathbf{A}_{\xi \eta}\frac{\partial{ \mathbf{u}}}{\partial \eta}  \rJump = 0, \quad \lJump \frac{\partial{\mathbf{u}}}{\partial t} \rJump = 0.
\end{equation}
The interface condition \eqref{eq:elastic_acoustic_interface_second} imposes the continuity of the flux and the time-derivatives of the wave fields,  and ensures that  the interface term $\mathrm{IT}$ vanishes
\begin{align*}
\mathrm{IT}\left( \frac{\partial{\mathbf{u}}}{\partial t}^T \left(  \sum_{\eta = x, y, z}\mathbf{A}_{\xi \eta}\frac{\partial{ \mathbf{u}}}{\partial \eta} \right) \right) =0.
\end{align*}

We again have
\begin{align*}
\frac{d}{dt}\left({E}^-(t)+ {E}^+(t)\right)   =   0,
\end{align*}
 and the total  energy is conserved.

\section{Perfectly matched layers}\label{sec:s3}
%
%

There are  two standard approaches, the complex coordinate stretching technique \cite{Chew1994, ZhaoCangellaris491508, PETROPOULOS1998184} and B{\'e}renger's splitting method \cite{berenger1994,Be_etAl}. The two approaches are mathematically analogous but yield two different formulations of the PML, the unsplit PML and the split-field PML, respectively. The difference lies in the choice of variables in the time domain. As  shown in \cite{KDuru2016}, using standard PML metric, by appropriate choice of variables the split-field PML \cite{berenger1994} can be formulated as the unsplit PML \cite{ABARBANEL1998331, KDuru2016, DURU201434}. This shows that, for linear problems,  the general solutions of the two PML formulations are identical. However, the complex coordinate stretching technique \cite{Chew1994} simplifies the construction of the PML for several hyperbolic PDEs \cite{KDuru2016,DuruKozdonKreiss2016,AppeloKreiss2006} 

Here, we will first demonstrate the construction of the PML for the  first order system \eqref{eq:first_order_linear_hyp_pde} by using the complex coordinate stretching technique \cite{Chew1994}. The techniques extends to linear second order hyperbolic systems in a straightforward manner \cite{DuKrSIAM,DuKr,Duruthesis2012,Komatitsch10.1093/gji/ggu219,MODAVE2017684,DIAZ20063820,KALTENBACHER2013407,Komatitsch10.1046/j.1365-246X.2003.01950.x}.  
 
Let the Laplace transform, in time, of $\mathbf{U}\left(x,y, z, t\right)$ be defined by
{
\begin{equation}
\widetilde{\mathbf{U}}(x,y,z,s)  = \int_0^{\infty}e^{-st}{\mathbf{U}}\left(x,y,z,t\right)\text{dt},  \quad s = a + ib, \quad \Re{s} = a > 0.
\end{equation}
}
We consider  a setup where the PML is included in all spatial coordinates.  Take the Laplace transform,  in time,  of    equation  \eqref{eq:first_order_linear_hyp_pde}. 
The PML can be constructed  in each coordinate,  $\xi = x, y, z $, using 
$\partial/\partial{\xi} \to 1/S_{\xi} \partial/\partial{\xi} $. 
Here
{
\begin{align}\label{eq:PML_metric}
\widetilde{\xi} = \int_{0}^{\xi}S_{\eta} d\eta,  \quad S_{\xi}= \gamma_\xi\left(\xi\right)\left(1 +\frac{d_{\xi}\left(\xi\right)}{s + \alpha_{\xi}(\xi)}\right),
\end{align}
}
 are the  complex cordinates and the PML metric, with $s$ denoting the Laplace dual time variable.
The PML damping functions   $d_{\xi}(\xi) \ge 0$  are  zero in the interior, where we solve the PDE \eqref{eq:first_order_linear_hyp_pde}, and take positive values in the PML.
 The grid stretching/compression parameter $\gamma_\xi>0$ \cite{DuKr,Duruthesis2012} are unity $\gamma_\xi \equiv 1$ in the interior  where $d_\xi = 0$ but increases/decreases inside the PML.
 The nonnegative  function $\alpha_\xi \ge 0$ is called the complex frequency shift (CFS) \cite{Kuzuoglu96}. The standard PML metric parameters will correspond to $\gamma_\xi =1$ and $\alpha_\xi =0$. It is possible to use even a more general  and complicated PML metric with more parameters, see for example \cite{Danieldoi:10.1137/050639107}.

 \subsection{First order systems}
 Now, take the Laplace transform in time of \eqref{eq:first_order_linear_hyp_pde},  \eqref{eq:bc_first_order} and \eqref{eq:interface_condition}, and introduce the complex change of variable $\partial/\partial{\xi} \to 1/S_{\xi} \partial/\partial{\xi} $. We have the PML in the Laplace space
 {
\begin{equation}\label{eq:linear_wave_PML_Laplace}
\begin{split}
\mathbf{P}^{-1} s \widetilde{\mathbf{U}}(x,y,z,s) &= \sum_{\xi = x, y, z}\mathbf{A}_{\xi}\frac{\partial{\widetilde{\mathbf{U}}(x,y,z,s)}}{\partial \widetilde{\xi}} =\sum_{\xi = x, y, z}\mathbf{A}_{\xi}\frac{1}{S_{\xi} }\frac{\partial{\widetilde{\mathbf{U}}(x,y,z,s)}}{\partial \xi},
  \end{split}
  \end{equation}
}
subject to
\begin{align}\label{eq:bc_first_order_Laplace}
\mathbf{L}\left( \widetilde{\mathbf{U}}\right) = 0, \quad \xi = \pm 1,
\end{align}
and
\begin{align}\label{eq:interface_condition_Laplace}
\mathbb{I}\mathbb{T}\left(\mathbf{A}_{\xi}, \widetilde{\mathbf{U}} \right) =0.
\end{align}
Note that the for the first order system the boundary condition \eqref{eq:bc_first_order_Laplace} and the interface condition \eqref{eq:interface_condition_Laplace} are not modified by the PML.
 Also note  the identity
 \begin{align}\label{eq:PML_metric_inverse}
\frac{1}{S_{\xi}} = \frac{1}{\gamma_\xi}-\frac{1}{S_\xi} \frac{d_\xi}{s+\alpha_\xi}.
\end{align}
In order to localize the PML in time, we introduce the auxiliary variable
\begin{equation}\label{eq:auxiliary_PML_Laplace}
\begin{split}
 \widetilde{\mathbf{w}}_{\xi}= \frac{1}{\left(s + \alpha_\xi\right)S_{\xi} }\mathbf{A}_{\xi}\frac{\partial{\widetilde{\mathbf{U}}}}{\partial \xi}, 
  \end{split}
  \end{equation}
  having
\begin{equation}\label{eq:linear_wave_PML_Laplace_simplify}
  \mathbf{P}^{-1} s \widetilde{\mathbf{U}} = \sum_{\xi = x, y, z}\mathbf{A}_{\xi}\frac{1}{S_{\xi} }\frac{\partial{\widetilde{\mathbf{U}}}}{\partial \xi} =  \sum_{\xi = x, y, z}\left[\frac{1}{\gamma_\xi}\mathbf{A}_{\xi}\frac{\partial{\widetilde{\mathbf{U}}}}{\partial \xi} - d_\xi  \widetilde{\mathbf{w}}_{\xi}\right].
   \end{equation}
  We then invert the Laplace transforms in \eqref{eq:linear_wave_PML_Laplace_simplify} and \eqref{eq:auxiliary_PML_Laplace},
  and we have the time-dependent PML 
\begin{equation}\label{eq:elastic_pml_1A}
\begin{split}
\mathbf{P}^{-1}\frac{\partial{\mathbf{U}}}{\partial t} = \sum_{\xi = x, y, z}\left[\frac{1}{\gamma_\xi}\mathbf{A}_{\xi}\frac{\partial{\mathbf{U}}}{\partial \xi}  - d_{\xi}\mathbf{w}_{\xi}\right],  \quad (x, y, z) \in \Omega, \quad t\ge 0,
  \end{split}
  \end{equation}
\begin{equation}\label{eq:elastic_pml_2A}
\begin{split}
\frac{\partial{\mathbf{w}_{\xi}}}{\partial t} 
 = \frac{1}{\gamma_\xi}\mathbf{A}_{\xi}\frac{\partial{\mathbf{U}}}{\partial \xi} -   \left(\alpha_\xi + d_{\xi}\right)\mathbf{w}_{\xi},
  \end{split}
  \end{equation}
with the boundary condition \eqref{eq:bc_first_order} and the interface conditions \eqref{eq:interface_condition}.
 In the interior, where the PML is not activated, we have $\gamma_\xi=1,d_{\xi} =0,$ which recovers the original equation \eqref{eq:first_order_linear_hyp_pde}.
 \begin{remark}
 Note that for the first order systems the boundary/interface conditions, \eqref{eq:bc_first_order} and \eqref{eq:interface_condition}, are not modified by the PML. However, the numerical treatments of the boundary/interface conditions will be critical for numerical stability when the PML is active.
 \end{remark}


\subsection{Second order systems}
For second order systems the ideas are similar. 
As above, take the Laplace transform of \eqref{eq:second_order_linear_hyp_pde} in time, and introduce the complex change of variable $\partial/\partial{\xi} \to 1/S_{\xi} \partial/\partial{\xi} $, $\partial/\partial{\eta} \to 1/S_{\eta} \partial/\partial{\eta} $. We have the PML in the Laplace space
{
\normalsize
\begin{equation}\label{eq:linear_2nd_wave_PML_Laplace}
\begin{split}
\boldsymbol{\Phi}^{-1} s^2 \widetilde{\mathbf{u}} &= \sum_{\widetilde{\xi} = \widetilde{x}, \widetilde{y}, \widetilde{z}}\frac{\partial }{\partial \widetilde{\xi}}\sum_{\widetilde{\eta} = \widetilde{x}, \widetilde{y}, \widetilde{z}}\mathbf{A}_{\widetilde{\xi} \widetilde{\eta}}\frac{\partial{\widetilde{\mathbf{u}}}}{\partial \widetilde{\eta}},\\
&= \sum_{\xi = x, y, z}\frac{1}{S_{\xi} }\frac{\partial }{\partial \xi}\sum_{\eta = x, y, z}\frac{1}{S_{\eta} }\mathbf{A}_{\xi \eta}\frac{\partial{\widetilde{\mathbf{u}} }}{\partial \eta},
\end{split}
\end{equation}
}
subject to the boundary condition
\begin{equation}\label{eq:linear_2nd_wave_PML_Laplace_BC}
\begin{split}
\widetilde{\mathbf{\mathcal{B}}}\widetilde{\mathbf{u}} = 0, \quad \widetilde{\mathbf{\mathcal{B}}} = \frac{1-r_{\xi}}{2}\mathbf{Z}\mathbf{I}s \pm \frac{1+r_{\xi}}{2}\sum_{\eta = x, y, z}\frac{1}{S_\eta}\mathbf{A}_{\xi \eta}\frac{\partial{}}{\partial \eta}, \quad \mathbf{Z} =\mathbf{Z} ^T\ge 0, \quad \xi = \pm 1.
\end{split}
\end{equation}
and the interface conditions,
\begin{equation}\label{eq:elastic_acoustic_interface_second_laplace_IC}
\lJump \sum_{\eta = x, y, z}\frac{1}{S_{\eta} }\mathbf{A}_{\xi \eta}\frac{\partial{\widetilde{\mathbf{u}} }}{\partial \eta}  \rJump = 0, \quad \lJump s{\widetilde{\mathbf{u}} }\rJump = 0,
\end{equation}
which ensure that  the interface term vanishes 
$$
\mathrm{IT}\left( \left(s{\widetilde{\mathbf{u}} }\right)^\dagger \left(  \sum_{\eta = x, y, z}\frac{1}{S_{\eta} }\mathbf{A}_{\xi \eta}\frac{\partial{\widetilde{\mathbf{u}} }}{\partial \eta}  \right) \right) =0.
$$

For second order systems, if the underlying boundary condition involves spatial derivatives the PML will transform the boundary conditions.

Localising the PML in time for the second order system may involve significant algebraic manipulations. To do this, first, we multiply \eqref{eq:linear_2nd_wave_PML_Laplace} with $S_xS_yS_z$,  and  \eqref{eq:linear_2nd_wave_PML_Laplace_BC} and \eqref{eq:elastic_acoustic_interface_second_laplace_IC} with $S_xS_yS_z/S_\xi$.  We have
{
\normalsize
\begin{equation*}
\begin{split}
\boldsymbol{\Phi}^{-1} s^2S_xS_yS_z \widetilde{\mathbf{u}}(x,y,z,s) &= \sum_{\xi = x, y, z}\frac{\partial }{\partial \xi}\sum_{\eta = x, y, z}\frac{S_xS_yS_z}{S_{\xi} S_{\eta} }\mathbf{A}_{\xi \eta}\frac{\partial{\mathbf{u}}}{\partial \eta}
\end{split}
\end{equation*}
}
with
$$
\widetilde{\mathbf{\mathcal{B}}}\widetilde{\mathbf{u}} = 0, \quad  \widetilde{\mathbf{\mathcal{B}}} = \frac{1-r_{\xi}}{2}\mathbf{Z}\mathbf{I}s\frac{S_xS_yS_z}{S_\xi} \pm \frac{1+r_{\xi}}{2}\sum_{\eta = x, y, z}\frac{S_xS_yS_z}{S_{\xi}S_\eta}\mathbf{A}_{\xi \eta}\frac{\partial{}}{\partial \eta}, \quad \xi = \pm 1,
$$
and
$$
\lJump \sum_{\eta = x, y, z}\frac{S_xS_yS_z}{S_{\xi}S_\eta}\mathbf{A}_{\xi \eta}\frac{\partial{\widetilde{\mathbf{u}} }}{\partial \eta}  \rJump = 0, \quad \lJump s{\widetilde{\mathbf{u}} }\rJump = 0.
$$

Second, introduce the auxiliary variables
$$
\widetilde{\boldsymbol{w}}_\eta = \left(\frac{\gamma_\eta S_xS_yS_z}{S_{\eta}\gamma_x\gamma_y\gamma_z }-1 \right)\widetilde{\mathbf{u}}(x,y,z,s), \quad \widetilde{\boldsymbol{v}}_{\eta \xi} = \frac{1}{ \left(\alpha + s\right) S_\xi} \frac{\partial }{\partial \eta}\left(\widetilde{\mathbf{u}}(x,y,z,s) + \widetilde{\boldsymbol{w}}_\eta \right),
$$

$$
\widetilde{\boldsymbol{\psi}} = \frac{1}{s + \alpha}\widetilde{\mathbf{u}}(x,y,z,s), 
\quad 
\widetilde{\boldsymbol{\phi}} = \frac{s}{(s + \alpha)^2}\widetilde{\mathbf{u}}(x,y,z,s), 
\quad 
\widetilde{\boldsymbol{\theta}} = \frac{1}{s + \alpha}\widetilde{\boldsymbol{\phi}}.
$$
And invert the Laplace transforms, we have
{
\footnotesize
\begin{equation}\label{eq:linear_2nd_wave_PML_1A}
\begin{split}
&\boldsymbol{\Phi}^{-1} \left(\frac{\partial^2 \mathbf{u}}{\partial t^2} +  \sum_{\xi =x,y,z  }d_\xi \frac{\partial  }{\partial t}\left(\mathbf{u} - \alpha{\boldsymbol{\psi}} \right) +   \prod_{\eta=x,y,z  }d_\eta\sum_{\xi=x,y,z  }\frac{1}{d_\xi}\left(\mathbf{u}-\alpha\left({\boldsymbol{\psi}}+ \boldsymbol{\phi}\right) \right)+ \prod_{\xi  }d_\xi \left(\boldsymbol{\phi}  - \alpha  \boldsymbol{\theta}\right) \right)\\
&= \sum_{\xi = x, y, z}\frac{1}{\gamma_\xi}\frac{\partial }{\partial \xi}\sum_{\eta = x, y, z}\left(\frac{1}{\gamma_\eta} \mathbf{A}_{\xi \eta}\left(\frac{\partial}{\partial \eta}\left({\mathbf{u}} + {\boldsymbol{w}}_\eta \right) - d_{\xi} \boldsymbol{v}_{\eta \xi}\right)\right),
\end{split}
\end{equation}
}
\begin{equation}\label{eq:linear_2nd_wave_PML_2A}
\begin{split}
\frac{\partial}{\partial t}\boldsymbol{\psi} =  \mathbf{u} - \alpha \boldsymbol{\psi}, \quad \frac{\partial}{\partial t}\boldsymbol{\phi} = \mathbf{u} - \alpha \left(\boldsymbol{\psi}+ \boldsymbol{\phi}\right), \quad \frac{\partial}{\partial t}\boldsymbol{\theta} = \boldsymbol{\phi}  - \alpha \boldsymbol{\theta},
\end{split}
\end{equation}
{
\normalsize
\begin{equation}\label{eq:linear_2nd_wave_PML_3A}
\begin{split}
\frac{\partial}{\partial t}\boldsymbol{w}_{\eta} =  \sum_{\xi \ne \eta }d_\xi\mathbf{u} +  \prod_{\xi \ne \eta }d_\xi  \boldsymbol{\psi} - \alpha \boldsymbol{w}_\eta, \quad \frac{\partial}{\partial t}\boldsymbol{v}_{\eta \xi} = \frac{\partial}{\partial \eta}\left({\mathbf{u}} + {\boldsymbol{w}}_\eta \right) - \left(d_{\xi} +\alpha\right)\boldsymbol{v}_{\eta \xi},
\end{split}
\end{equation}
}
subject to the boundary conditions,
{
\begin{equation}\label{eq:bc_second_order_bc_PML}
\begin{split}
&\frac{1-r_{\xi}}{2}\mathbf{Z}\left(\frac{\partial \mathbf{u}}{\partial t} +  \sum_{\eta \ne \xi  }d_\eta \left(\mathbf{u} - \alpha{\boldsymbol{\psi}} \right) +   \prod_{\eta \ne \xi  }d_\eta \boldsymbol{\phi} \right) \\
&\pm\\
&\frac{1+r_{\xi}}{2} \sum_{\eta = x, y, z}\left(\frac{1}{\gamma_\eta} \mathbf{A}_{\xi \eta}\left(\frac{\partial}{\partial \eta}\left({\mathbf{u}} + {\boldsymbol{w}}_\eta \right) - d_{\xi} \boldsymbol{v}_{\eta \xi}\right) \right) = 0,
\end{split}
\end{equation}
}
and the interface conditions 
\begin{equation}\label{eq:elastic_acoustic_interface_second_PML}
\lJump \sum_{\eta = x, y, z} \frac{1}{\gamma_\eta} \mathbf{A}_{\xi \eta}\left(\frac{\partial}{\partial \eta}\left({\mathbf{u}} + {\boldsymbol{w}}_\eta \right) - d_{\xi} \boldsymbol{v}_{\eta \xi}\right)  \rJump = 0, \quad \lJump \frac{\partial{\mathbf{u}}}{\partial t} \rJump = 0.
\end{equation}
The boundary/interface conditions, \eqref{eq:bc_second_order}--\eqref{eq:bc_second_order_operator} and \eqref{eq:elastic_acoustic_interface_second} for the second order wave equation are transformed to \eqref{eq:bc_second_order_bc_PML} and \eqref{eq:elastic_acoustic_interface_second_PML} for the second order PML equations \eqref{eq:linear_2nd_wave_PML_1A}--\eqref{eq:linear_2nd_wave_PML_3A}.
\begin{remark}
 Note that for the second order systems the boundary/interface conditions, \eqref{eq:bc_second_order}--\eqref{eq:bc_second_order_operator} and \eqref{eq:elastic_acoustic_interface_second}, are are transformed in a nontrivial manner to the boundary conditions  \eqref{eq:bc_second_order_bc_PML} and \eqref{eq:elastic_acoustic_interface_second_PML}  for the PML equations. This is in contrast to the PML for first order systems where the boundary conditions remain unchanged for the PML.
 \end{remark}

Since, we have used the well known PML metric \eqref{eq:PML_metric}, the PML models \eqref{eq:elastic_pml_1A}-\eqref{eq:elastic_pml_2A} and \eqref{eq:linear_2nd_wave_PML_1A}--\eqref{eq:linear_2nd_wave_PML_3A} can be shown to be  analogous to other  PML models such as  \cite{KDuru2016,DuruKozdonKreiss2016,AppeloKreiss2006,Be_etAl,DuKr,XieKomatitschMartinMatzen,SkAdCr,DuKrSIAM}.     We will initialize the PML with zero initial data and  terminate the PML  \eqref{eq:elastic_pml_1A}-\eqref{eq:elastic_pml_2A}  with the boundary conditions  \eqref{eq:bc_first_order} and the PML \eqref{eq:linear_2nd_wave_PML_1A}--\eqref{eq:linear_2nd_wave_PML_3A} with the boundary condition \eqref{eq:bc_second_order_operator}.
We again note  that the PML absorption functions  and auxiliary  functions  vanish almost everywhere except in the layers defining the PML where $d_\xi \ge 0$.

\section{Mathematical analysis of PMLs at the continuous PDE level}\label{sec:math_analysis}
The PML transformations \eqref{eq:PML_metric} and \eqref{eq:linear_wave_PML_Laplace}  have some important mathematical and physical properties which can be revealed through rigorous mathematical analysis of PMLs. For a first order system  and a second order system modelling the same physical phenomena,  the second order PML   \eqref{eq:linear_2nd_wave_PML_1A}--\eqref{eq:linear_2nd_wave_PML_3A} can be rewritten as the first order PML \eqref{eq:elastic_pml_1A}-\eqref{eq:elastic_pml_2A} model by introducing suitable variables. This implies that the mathematical properties of the systems \eqref{eq:elastic_pml_1A}-\eqref{eq:elastic_pml_2A} and \eqref{eq:linear_2nd_wave_PML_1A}--\eqref{eq:linear_2nd_wave_PML_3A} are equivalent. In this review, we will focus on the first order system, and discuss \emph{perfect matching},  \emph{well-posedness} and \emph{stability}.  Analysis of well-posedness and stability of the PML has  been an area of active research \cite{Duruthesis2012,Be_etAl,SkAdCr,AppeloKreiss2006,DuruGabKreiss2019,HalpernPetit-BergezRauch2011,Becache1296843,becache_joly_2002,Danieldoi:10.1137/050639107}. We will try to summarise the results without going into too much details of the proofs. To simplify the discussions we will often consider the PML only in the $x$-direction, with  $d_x(x) \ge 0$ and $\gamma_\xi = 1$, $d_\xi = 0$ for $\xi = y, z$.

\subsection{Perfect matching} 
We will now consider a whole space PML problem, with $\Omega = \mathbb{R}^3$, and  the PML in the $x$-direction. In the interior $x \le 0$ we have $d_x(x) = 0$, $\gamma_x =1 $ and in the PML $x > 0$, $d_x(x) > 0$, $\gamma_x(x) >0$.
Probably, the most important property of the PML \eqref{eq:elastic_pml_1A}-\eqref{eq:elastic_pml_2A} is that the equations are perfectly matched to the underlying hyperbolic PDEs. This means that the restriction of the general solutions to \eqref{eq:elastic_pml_1A}-\eqref{eq:elastic_pml_2A} in the interior of the domain, $x \le 0$, coincides with the general solutions to \eqref{eq:first_order_linear_hyp_pde}. Therefore, there are no reflections as waves propagate from the interior, $x\le 0$ into the layer $x> 0$. 

We will formalise the discussion below.
\begin{definition}
Let $\mathbf{U}$ denote the general solution of  the wave equation \eqref{eq:first_order_linear_hyp_pde} in the unbounded domain and $\mathbf{U}_{\mathrm{pml}}$ denote the general solution of the PML  \eqref{eq:elastic_pml_1A}-\eqref{eq:elastic_pml_2A}. The equations \eqref{eq:elastic_pml_1A}-\eqref{eq:elastic_pml_2A} are perfectly matched to the wave equation \eqref{eq:first_order_linear_hyp_pde} if 
\begin{align*}
\mathbf{U}_{\mathrm{pml}} = \mathbf{U}, \quad x\le 0, \quad \forall t \ge 0.
\end{align*}
\end{definition}

There are two standard methods \cite{Berenger999615, berenger1994, Danieldoi:10.1137/050639107,Duruthesis2012} that have been used to study the perfect matching property of the PML. The approach \cite{Berenger999615,berenger1994} uses plane wave analysis and only accounts for propagating modes. The technique \cite{Danieldoi:10.1137/050639107,Duruthesis2012}, which is rooted in the construction of the general solutions to the wave equations in the Laplace--Fourier space, is more general since it includes both the propagating mode regime and the evanescent mode regime. 

Consider the wave equation \eqref{eq:first_order_linear_hyp_pde}, and take the Laplace transform in time and Fourier transforms in the tangential directions
 {
\begin{equation}\label{eq:linear_wave_Laplace_Fourier}
\begin{split}
\mathbf{P}^{-1} s \widehat{\widetilde{\mathbf{U}}}\left(x,ik_y,ik_z,s\right) &= \mathbf{A}_{x}\frac{d }{d {x}}{\widehat{\widetilde{\mathbf{U}} }\left(x,ik_y,ik_z,s\right)} - \sum_{\xi = y, z}ik_{\xi}\mathbf{A}_{\xi}\widehat{\widetilde{\mathbf{U}}}\left(x,ik_y,ik_z,s\right)\\
  \end{split}
  \end{equation}
}
For \eqref{eq:linear_wave_Laplace_Fourier} we can construct modal solutions 
\begin{align*}
\widehat{\widetilde{\mathbf{U}}}\left(x,ik_y,ik_z,s\right)  =  \widehat{\widetilde{\mathbf{U}}}_0\left(ik_y,ik_z,s\right) e^{\lambda x}
\end{align*}
where $\lambda$ satisfies
\begin{align}\label{eq:charac_eqn}
 \left|\Upsilon\left(\lambda,ik_y,ik_z,s\right) \right| =0, \quad \Upsilon\left(\lambda,ik_y,ik_z,s\right) = s\mathbf{P}^{-1}  + \sum_{\xi = y, z}ik_{\xi}\mathbf{A}_{\xi}- \lambda\mathbf{A}_{x}.
\end{align}
and $\widehat{\widetilde{\mathbf{U}}}_0\left(ik_y,ik_z,s\right)$ is the eigenfunction of $\Upsilon\left(\lambda,ik_y,ik_z,s\right)$.
Similarly, for the PML \eqref{eq:linear_wave_PML_Laplace} we can construct modal solutions 
\begin{align*}
\widehat{\widetilde{\mathbf{U}}}_{\mathrm{pml}}\left({x},ik_y,ik_z,s\right)  =  \widehat{\widetilde{\mathbf{U}}}_0\left(ik_y,ik_z,s\right) e^{\lambda \widetilde{x}},
\end{align*}
with $\lambda$ satisfying \eqref{eq:charac_eqn}.
Note that in the interior $x\le 0$ we have $d_x = 0$, $\gamma_x = 1$, yielding $x = \widetilde{x}$, and we have
\begin{align*}
\widehat{\widetilde{\mathbf{U}}}_{\mathrm{pml}}\left({x},ik_y,ik_z,s\right)  =  \widehat{\widetilde{\mathbf{U}}}\left(x,ik_y,ik_z,s\right), \quad \forall x \le 0.
\end{align*}
The solutions $\widehat{\widetilde{\mathbf{U}}}\left({x},ik_y,ik_z,s\right)$ and $\widehat{\widetilde{\mathbf{U}}}_{\mathrm{pml}}\left({x},ik_y,ik_z,s\right)=\widehat{\widetilde{\mathbf{U}}}\left(\widetilde{x},ik_y,ik_z,s\right)$ are identical in the interior $x \le 0$. The solutions are perfectly matched by construction.

\begin{remark}
Note the similarity between \eqref{eq:linear_wave_PML_Laplace} and \eqref{eq:linear_wave_Laplace_Fourier}. Since the eigenfunction $\widehat{\widetilde{\mathbf{U}}}_0\left(ik_y,ik_z,s\right)$ is unchanged the PML \eqref{eq:elastic_pml_1A}-\eqref{eq:elastic_pml_2A} corresponds to \eqref{eq:linear_wave_PML_Laplace}, a complex change of coordinates in the Laplace space. Any wave propagating to the right in the interior corresponds to a wave propagating to the right in the PML. It then follows by construction that \eqref{eq:elastic_pml_1A}-\eqref{eq:elastic_pml_2A}  is perfectly matched to the underlying hyperbolic PDE \eqref{eq:first_order_linear_hyp_pde}, see \cite{Danieldoi:10.1137/050639107,Duruthesis2012}.
\end{remark}

\subsection{Well-posedness}
 We begin by introducing the notions of well-posedness and stability. By a well-posed problem, we mean that there is a unique solution which depends continuously on the data of the problem. To be precise, consider the Cauchy problem
\begin{align}\label{eq:PML_IVP}
\frac{\partial \mathbf{Q}}{\partial t} = \mathcal{D}\left({\partial }/{\partial \xi}, d_\xi, \gamma_\xi, \alpha\right) \mathbf{Q}, \quad \mathbf{Q}(\bar{\xi}, 0) = \mathbf{Q}_0\left(\bar{\xi}\right), \quad \bar{\xi} = (x, y, z)
\end{align}
Here, $\mathbf{Q} = \left(\mathbf{U}, \mathbf{w}_\xi\right)^T$,  the differential symbol $\mathcal{D}\left({\partial }/{\partial \xi}, d_\xi, \gamma_\xi, \alpha\right) $ denotes the spatial operator, and for the PML \eqref{eq:elastic_pml_1A}-\eqref{eq:elastic_pml_2A}  it is given by
\begin{equation}
\begin{split}
\mathcal{D}\left({\partial }/{\partial \xi}, d_\xi, \gamma_\xi, \alpha\right)  
= 
\underbrace{\begin{pmatrix}
\sum_{\xi = x, y, z}\frac{1}{\gamma_\xi}\mathbf{P}\mathbf{A}_{\xi}\frac{\partial{}}{\partial \xi}, &  \mathbf{0}\\
\frac{1}{\gamma_\xi}\mathbf{A}_{\xi}\frac{\partial{}}{\partial \xi}, &  \mathbf{0}
\end{pmatrix}}_{\text{principal part}}
-
\begin{pmatrix}
\mathbf{0}, &  \sum_{\xi = x, y, z}d_{\xi}\mathbf{P}\\
\mathbf{0}, &  \left(d_{\xi}+\alpha\right)\mathbf{I}
\end{pmatrix}.
  \end{split}
  \end{equation}
\begin{definition}\label{def:wellposedness}
The Cauchy problem \eqref{eq:PML_IVP} is weakly (resp. strongly) well-posed if for every $t \ge 0$
$$
\| \mathbf{Q}\| \le K e^{\kappa t}\| \mathbf{Q}_0\|_{H^{\omega}},
$$
for $\mathbf{Q}_0$ giving in the Sobolev space $H^{\omega}$, ${\omega}>0$ (resp. ${\omega}=0$), where $\kappa$ and $K$ are constants independent of $\mathbf{Q}_0$.
\end{definition}

For the PML \eqref{eq:elastic_pml_1A}-\eqref{eq:elastic_pml_2A}, we consider the  principal part and the coefficient matrices
$$
\widehat{\mathbf{A}}_{\xi} = \begin{pmatrix}
\frac{1}{\gamma_\xi}\mathbf{P}\mathbf{A}_{\xi}, &  \mathbf{0}\\
\frac{1}{\gamma_\xi}\mathbf{A}_{\xi}, &  \mathbf{0}
\end{pmatrix} \sim  \begin{pmatrix}
\frac{1}{\gamma_\xi}\mathbf{P}\mathbf{A}_{\xi}, &  \mathbf{0}\\
\frac{1}{\gamma_\xi}\mathbf{P}\mathbf{A}_{\xi}, &  \mathbf{0}
\end{pmatrix},
$$
where $\sim$ denotes the similarity of matrices. Since the underlying system \eqref{eq:first_order_linear_hyp_pde}, without the PML, is strongly hyperbolic the coefficient matrices $\mathbf{P}\mathbf{A}_{\xi}$ have real eigenvalues and a complete system of eigenvectors. Thus, by inspection the matrices $\widehat{\mathbf{A}}_{\xi}$ have real eigenvalues with multiple eigenvectors. Therefore, by standard definitions \cite{KreissOliger1972}, the PML \eqref{eq:elastic_pml_1A}-\eqref{eq:elastic_pml_2A} is weakly hyperbolic, and thus weakly well-posed.
For some systems such as the acoustic wave equation \cite{KDuru2016} or the Maxwell's equations \cite{ABARBANEL1998331,KDuru2016}, it is possible to rewrite the equations, by introducing new variables such that the PML is strongly hyperbolic, and thus strongly well-posed.
However, the resulting strongly hyperbolic system is linearly equivalent to the original weakly hyperbolic systems, and their general solutions are identical \cite{ABARBANEL1998331, KDuru2016, DURU201434}.

For time-dependent problems, it is not sufficient that the PML is well-posed, it must also be stable in the sense that the solution remains bounded after a sufficiently long time. By Definition \ref{def:wellposedness}, a well-posed problem can  support exponentially growing solutions. This is of course undesirable of an absorbing model. Any growth in the layer can propagate into the computational domain and pollute the solutions everywhere. In order for the PML to be useful it must therefore be stable.
\subsection{Stability}
The stability  of the PML has attracted considerable attention in the literature, 
\cite{Duruthesis2012,Be_etAl,SkAdCr,AppeloKreiss2006,DuruGabKreiss2019}.   The perfectly matched layer is indeed a variable coefficient problem, but many classical techniques require constant coefficients, and therefore  one often considers corresponding frozen coefficient problems.  Further simplifications involve considering sub-problems  where only one PML damping is nonzero,  eg. the $x$-dependent PML strip with  $d_x > 0$ and $d_{\xi} =0$ for $\xi \ne x$.
%
%
\subsubsection{Fourier analysis for the PML Cauchy problem}
\label{sec:FoureierAnalysis}
~\
Consider the Cauchy problem for Equations \eqref{eq:elastic_pml_1A}-\eqref{eq:elastic_pml_2A}, with constant coefficients, $d_x=d>0$, $d_y=d_z=0$,   $\alpha_x=\alpha\ge0$, $\alpha_y=\alpha_z=0$, and $\gamma_x = \gamma> 0$.
The main idea is to  look for wave-like solutions of the form
\begin{align}\label{eq:planewave}
\mathbf{Q} = \mathbf{Q}_0e^{ik_x x + ik_y y + ik_z z + s t} , \quad s = \Re{s} + i\Im{s},
\end{align}
where $s\in \mathbb{C}$ is the unknown eigenvalue and will be determined. If there is an eigenvalue with a positive real part the PML is unstable, since this will correspond to a plane wave solution with exponentially growing amplitude.
We introduce $|\mathbf{k}| = \sqrt{\left(k_x/\gamma_x\right)^2 + k_y^2 + k_z^2}$ and the scaling
$$
\lambda = \frac{s}{|\mathbf{k}|}, \quad k_1 = \frac{k_x}{\gamma_x|\mathbf{k}|},  \quad k_2 = \frac{k_y}{|\mathbf{k}|},  \quad k_3 = \frac{k_z}{|\mathbf{k}|},  \quad \epsilon = \frac{d}{|\mathbf{k}|},  \quad \nu = \frac{\alpha}{|\mathbf{k}|}, \quad S_x\left(\lambda, \epsilon, \nu\right)= 1 +\frac{\epsilon}{\lambda + \nu}.
$$
Correspondingly, if there are $\Re{\lambda} > 0$, the PML is unstable.

We insert the plane wave solution \eqref{eq:planewave} in the PML and obtain the dispersion relation
{
\begin{align}\label{eq:PML_dispersion_relation}
&F\left(-i \lambda,  \frac{1}{S_x\left(\lambda, \epsilon, \nu\right)}k_1, k_2, k_3\right):= \Big|-i\lambda\mathbf{I}  + \frac{1}{S_x\left(\lambda, \epsilon, \nu\right)}k_x\mathbf{P}\mathbf{A}_{x} + \sum_{\xi = y,z}{k_{\xi}}\mathbf{P} \mathbf{A}_{\xi}\Big| = 0,
\end{align}
}
with $S_x\left(\lambda, \epsilon, \nu\right) = 1 + \frac{\epsilon}{\lambda + \nu}$.
The scaled eigenvalue $\lambda$ is the root of the complicated nonlinear dispersion relation \eqref{eq:PML_dispersion_relation} for the PML.
 The roots $\lambda$ can be difficult to determine. However standard perturbation arguments yield the following well-known result, the so-called \emph{geometric stability condition}, see Definition \ref{def:geometric_condition} and \cite{Duruthesis2012,Be_etAl,AppeloKreiss2006}.
 
\begin{figure}[h!]
\tikzset{every picture/.style={line width=1.5pt}} 

\begin{tikzpicture}[x=0.7pt,y=0.7pt,yscale=-1,xscale=1]
%
\begin{scope}
\draw    (98,252) -- (97.01,38) ;
\draw [shift={(97,36)}, rotate = 449.73] [color={rgb, 255:red, 0; green, 0; blue, 0 }  ][line width=0.75]    (10.93,-3.29) .. controls (6.95,-1.4) and (3.31,-0.3) .. (0,0) .. controls (3.31,0.3) and (6.95,1.4) .. (10.93,3.29)   ;
\draw    (98,252) -- (335,252) ;
\draw [shift={(337,252)}, rotate = 180] [color={rgb, 255:red, 0; green, 0; blue, 0 }  ][line width=0.75]    (10.93,-3.29) .. controls (6.95,-1.4) and (3.31,-0.3) .. (0,0) .. controls (3.31,0.3) and (6.95,1.4) .. (10.93,3.29)   ;
\draw    (117,173) .. controls (187,124) and (223,163) .. (240,220) ;
\draw    (98,252) -- (169.82,153.62) ;
\draw [shift={(171,152)}, rotate = 486.13] [color={rgb, 255:red, 0; green, 0; blue, 0 }  ][line width=0.75]    (10.93,-3.29) .. controls (6.95,-1.4) and (3.31,-0.3) .. (0,0) .. controls (3.31,0.3) and (6.95,1.4) .. (10.93,3.29)   ;
\draw    (171,152) -- (149.51,70.93) ;
\draw [shift={(149,69)}, rotate = 435.15] [color={rgb, 255:red, 0; green, 0; blue, 0 }  ][line width=0.75]    (10.93,-3.29) .. controls (6.95,-1.4) and (3.31,-0.3) .. (0,0) .. controls (3.31,0.3) and (6.95,1.4) .. (10.93,3.29)   ;
%
\draw (240,210) node [anchor=north west][inner sep=0.75pt]   [align=left] {$S$};
\draw (142,197) node [anchor=north west][inner sep=0.75pt]   [align=left] {$V_p$};
\draw (165,97.4) node [anchor=north west][inner sep=0.75pt]    {$V_g$};
\draw (55,57.4) node [anchor=north west][inner sep=0.75pt]    {$k_y/\omega$};
\draw (283,260) node [anchor=north west][inner sep=0.75pt]    {$k_x/\omega$};
\draw (205,57.4) node [anchor=north west][inner sep=2.75pt]    {$\mathbf{Unstable}$};
\end{scope}
%
%
\begin{scope}[shift = {(320,0)}]
\draw    (98,252) -- (97.01,38) ;
\draw [shift={(97,36)}, rotate = 449.73] [color={rgb, 255:red, 0; green, 0; blue, 0 }  ][line width=0.75]    (10.93,-3.29) .. controls (6.95,-1.4) and (3.31,-0.3) .. (0,0) .. controls (3.31,0.3) and (6.95,1.4) .. (10.93,3.29)   ;
\draw    (98,252) -- (335,252) ;
\draw [shift={(337,252)}, rotate = 180] [color={rgb, 255:red, 0; green, 0; blue, 0 }  ][line width=0.75]    (10.93,-3.29) .. controls (6.95,-1.4) and (3.31,-0.3) .. (0,0) .. controls (3.31,0.3) and (6.95,1.4) .. (10.93,3.29)   ;
\draw    (98,252) -- (169.82,153.62) ;
\draw [shift={(171,152)}, rotate = 486.13] [color={rgb, 255:red, 0; green, 0; blue, 0 }  ][line width=0.75]    (10.93,-3.29) .. controls (6.95,-1.4) and (3.31,-0.3) .. (0,0) .. controls (3.31,0.3) and (6.95,1.4) .. (10.93,3.29)   ;
\draw    (171,152) -- (244.33,104.09) ;
\draw [shift={(246,103)}, rotate = 506.84] [color={rgb, 255:red, 0; green, 0; blue, 0 }  ][line width=0.75]    (10.93,-3.29) .. controls (6.95,-1.4) and (3.31,-0.3) .. (0,0) .. controls (3.31,0.3) and (6.95,1.4) .. (10.93,3.29)   ;
\draw    (149,121) .. controls (149,146) and (196,177) .. (245,161) ;
\draw (240,152) node [anchor=north west][inner sep=0.75pt]   [align=left] {$S$};
\draw (142,197) node [anchor=north west][inner sep=0.75pt]   [align=left] {$V_p$};
\draw (201,102.5) node [anchor=north west][inner sep=0.75pt]    {$V_g$};
\draw (55,57.4) node [anchor=north west][inner sep=0.75pt]    {$k_y/\omega$};
\draw (283,260) node [anchor=north west][inner sep=0.75pt]    {$k_x/\omega$};
\draw (205,57.4) node [anchor=north west][inner sep=2.75pt]    {$\mathbf{Stable}$};
\end{scope}
\end{tikzpicture}
\caption{The geometric stability condition of the PML}
\label{fig:stabilitycondition}
\end{figure}
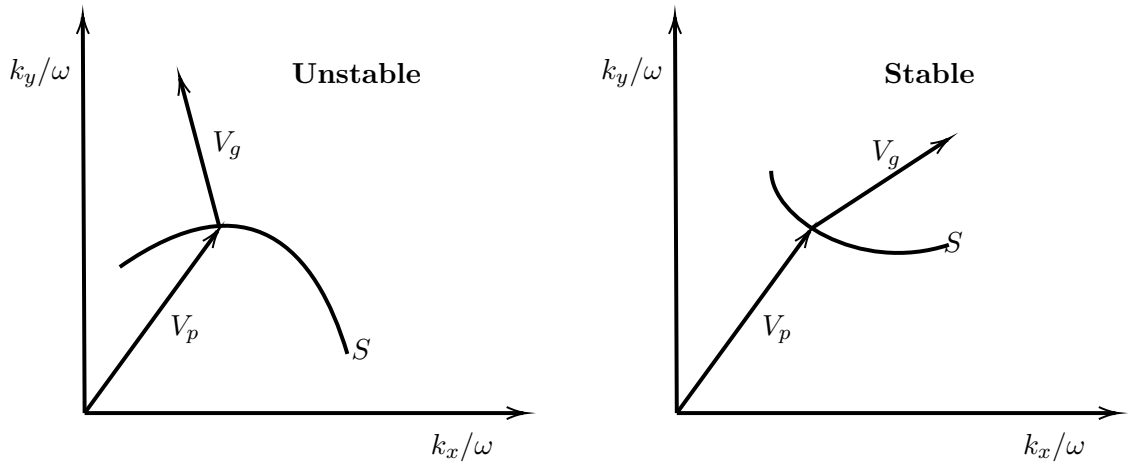
\begin{figure}[h!]
    \begin{subfigure}{0.5\textwidth}
    \centering
    \includegraphics[width=\linewidth]{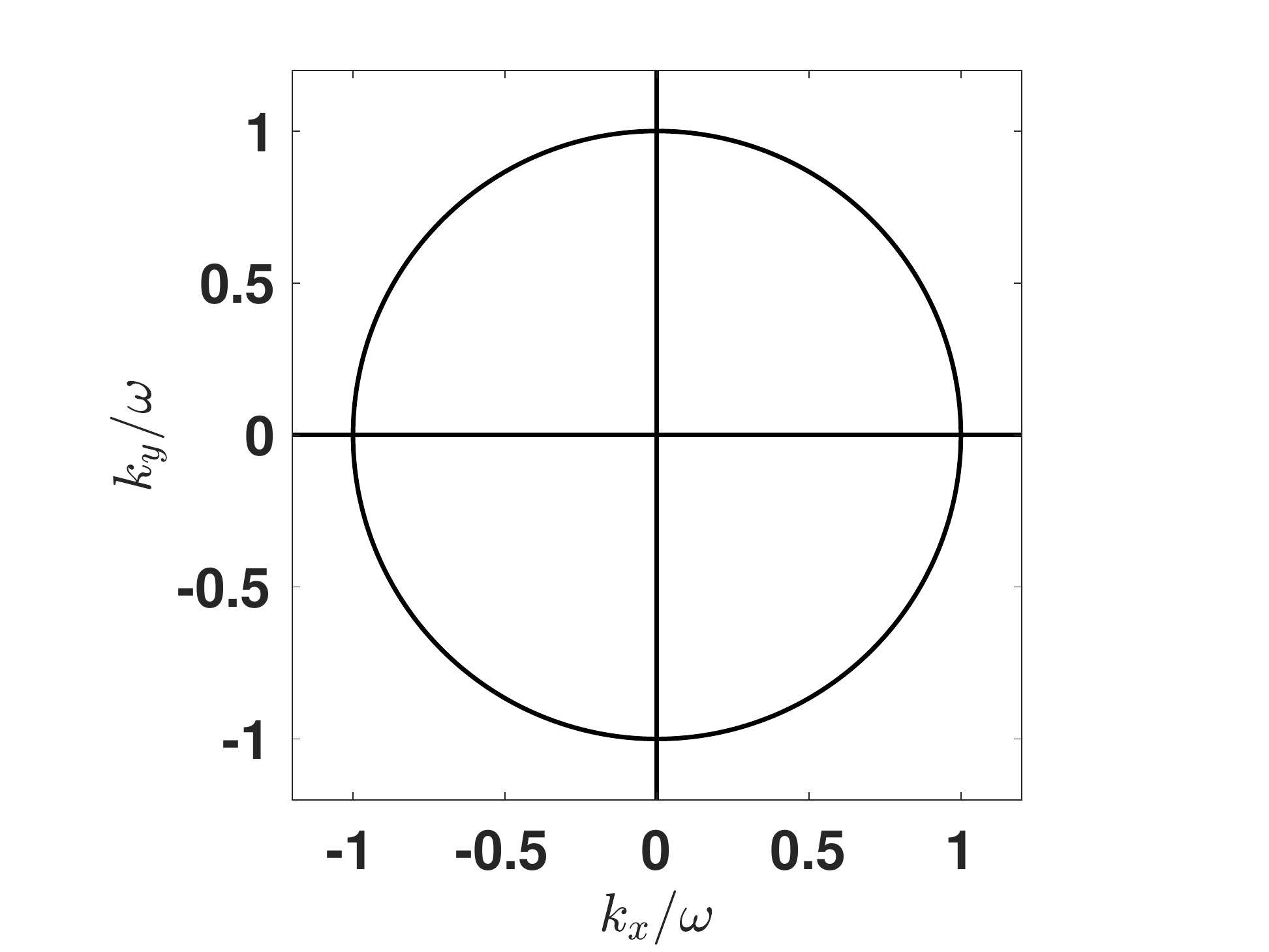} 
    \caption{Acoustic medium}
    \label{fig:Acousticmedium}
    \end{subfigure}
    \begin{subfigure}{0.5\textwidth}
    \centering
    \includegraphics[width=\linewidth]{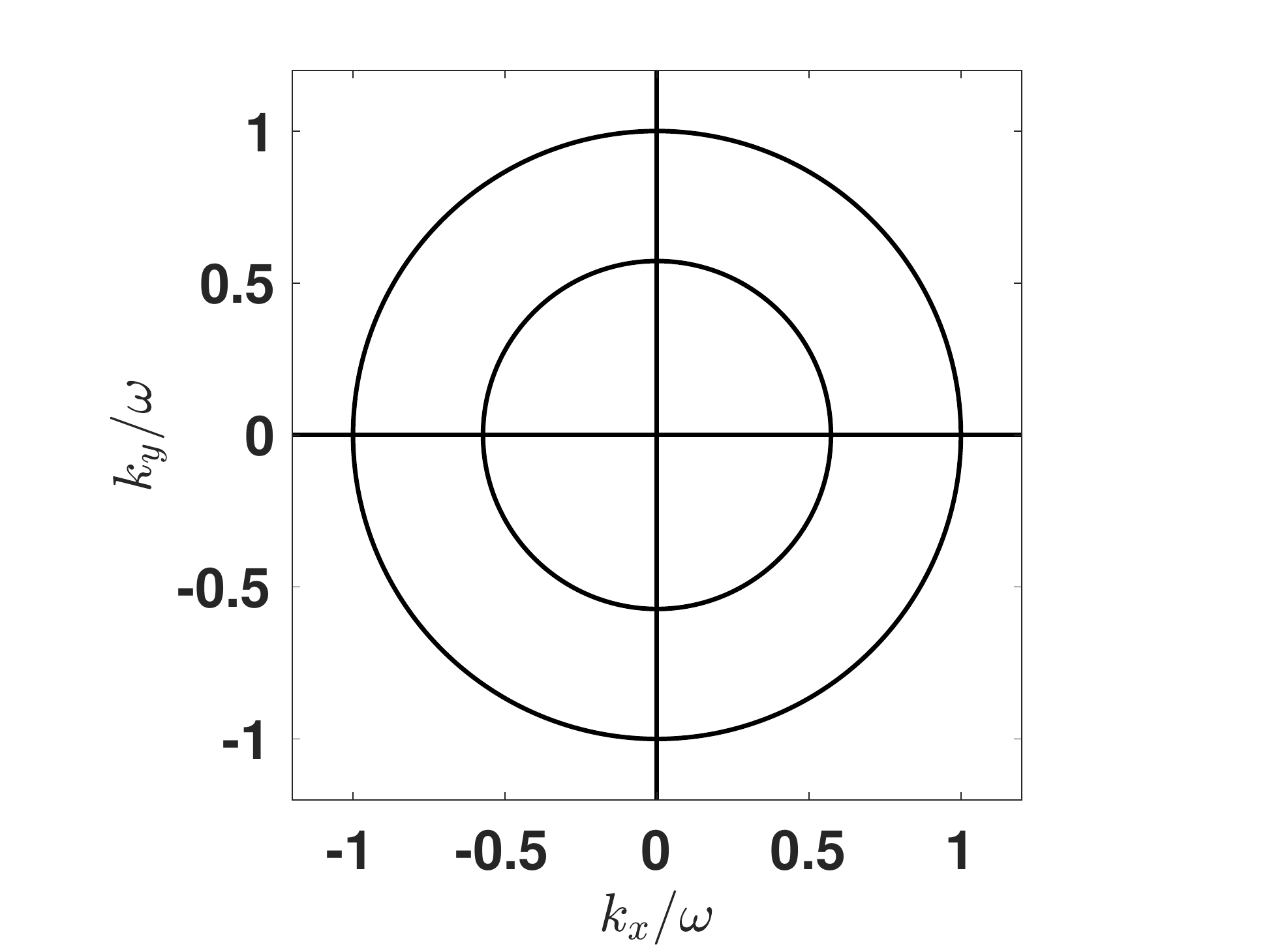} 
    \caption{Isotropic elastic wave solid}
    \label{fig:Isotropicelastic}
    \end{subfigure}
    \begin{subfigure}{0.5\textwidth}
    \centering
    \includegraphics[width=\linewidth]{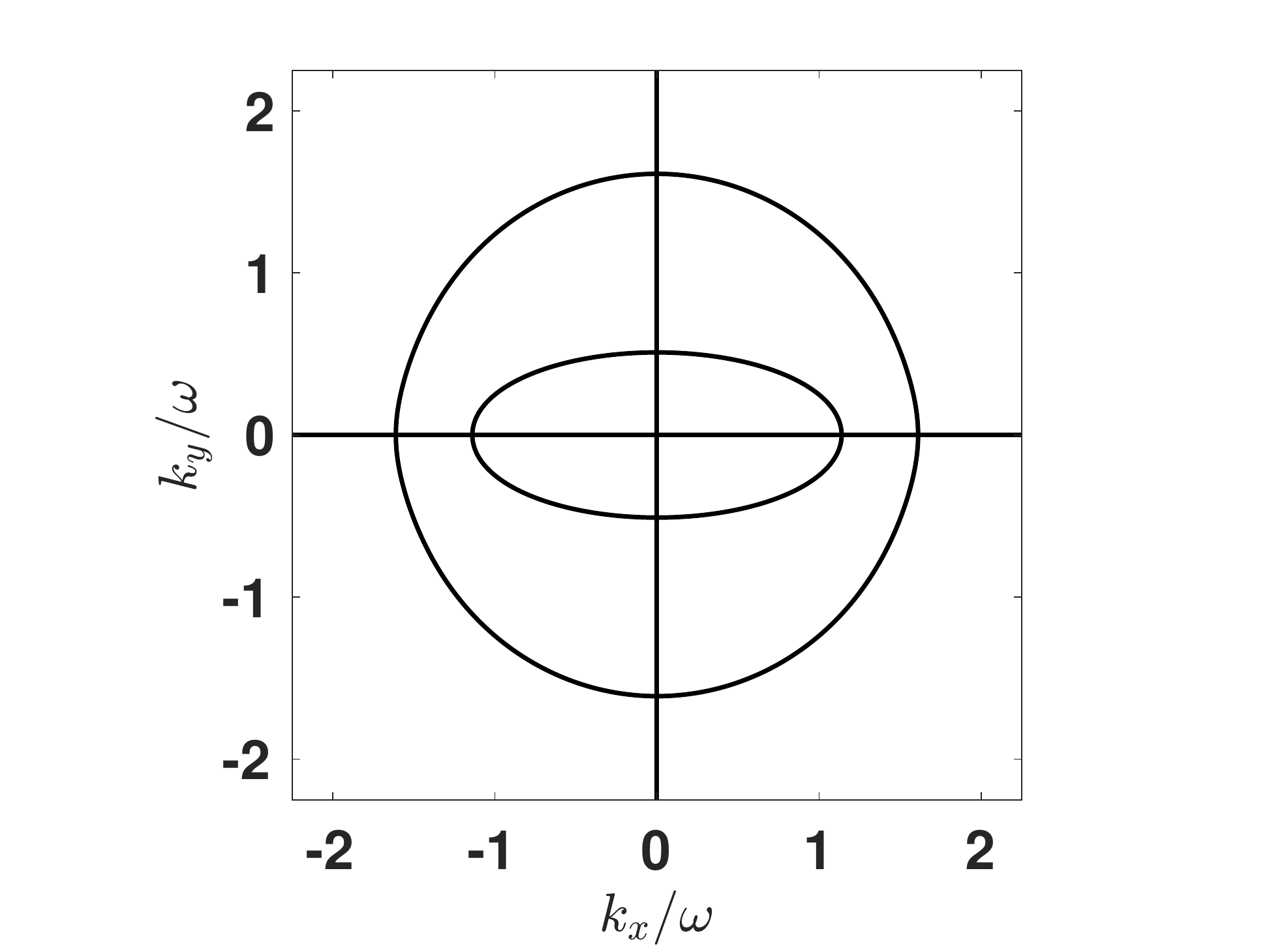} 
    \caption{AM1: Stable anisotropic elastic solid}
    \label{fig:AM1}
    \end{subfigure}
    \begin{subfigure}{0.5\textwidth}
    \centering
    \includegraphics[width=\linewidth]{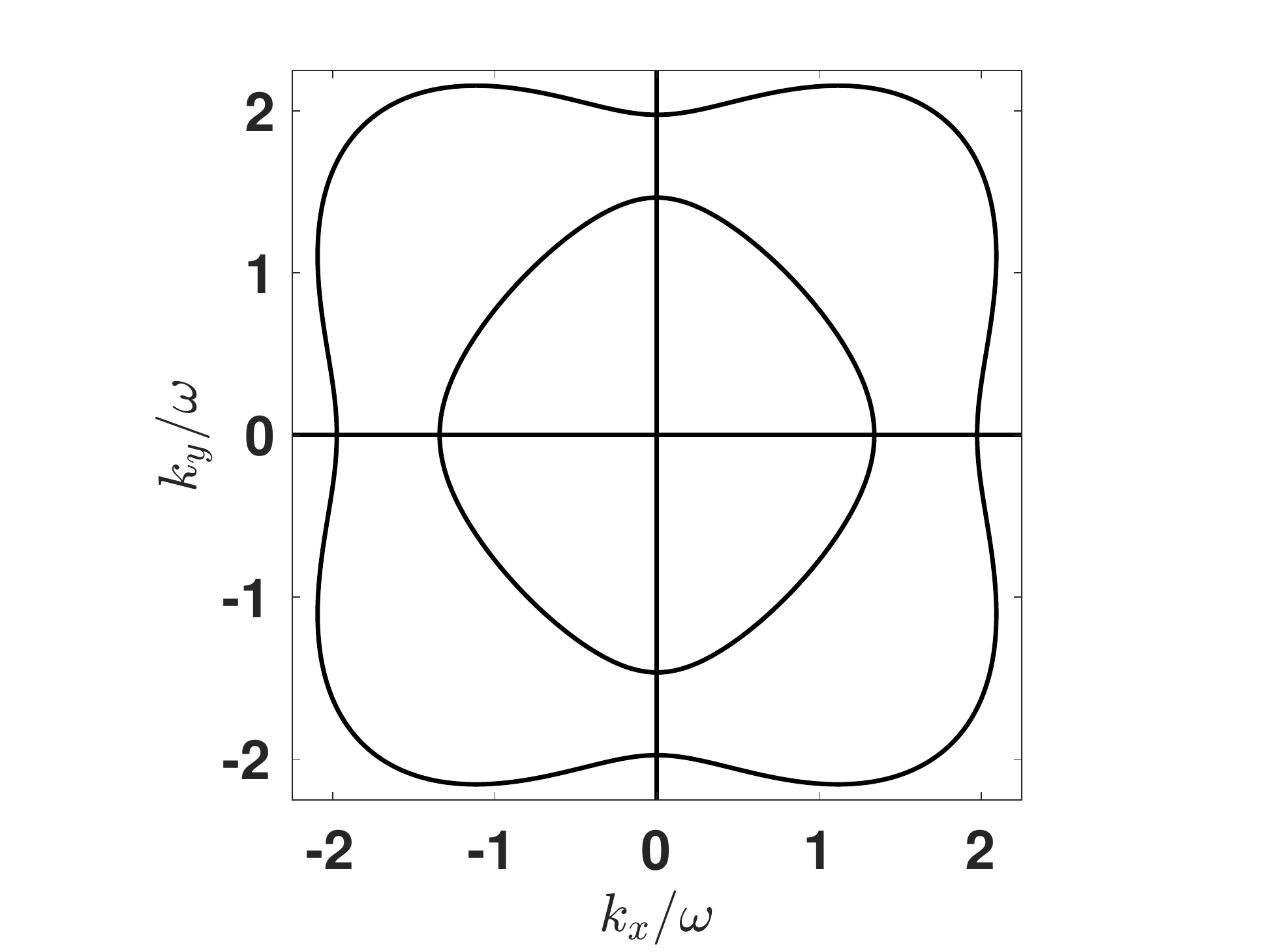} 
    \caption{AM2: Unstable anisotropic elastic solid}
    \label{fig:AM2}
    \end{subfigure}
    \caption{Slowness diagrams showing the solutions of the dispersion relation for different wave media in two space dimensions.}
     \label{fig:dispersion_relation}
\end{figure}

\begin{theorem}[Necessary  condition for stability]\label{theo:high_frequency_stability}
Consider the constant coefficient PML, with $d_x = d >0$, $\alpha_x = \alpha \ge 0$ and $\gamma_x = \gamma> 0$. 
If the geometric stability condition in the $x$-direction is violated, then there are growing modes with $\Re{\lambda} > 0$ at all sufficiently high frequencies.
\end{theorem}
The theorem states that if the underlying medium violates the {\it geometric stability condition}, then there are high frequency wave modes with exponentially growing amplitude. See also Figure \ref{fig:stabilitycondition} for a schematic explanation of Theorem \ref{theo:high_frequency_stability}.
\begin{remark}
At intermediate frequencies the parameter $\gamma_x$ can be used to stabilise unstable modes, by compressing the grid using $\gamma_x < 1$.
 However at sufficiently high frequencies the instability persists. See  \cite{Duruthesis2012,DuKr} for details.
\end{remark}
In Figure \ref{fig:dispersion_relation} we display the slowness diagrams of different wave media. Note that the geometric stability condition is satisfied for the first three diagrams in Figure \ref{fig:dispersion_relation}(a)--(c).
The last diagram Figure \ref{fig:dispersion_relation}(d) which is an anisotropic elastic medium violates the geometric stability condition in certain parts of the slowness curve.

A typical  example where the geometric stability condition is violated is in  advective acoustics or linearised  Euler equations with non-vanishing mean flow \cite{HU1996201,TAM1998213,Lions2002,ABARBANEL1999266,HU2001455,Hagstrom2003,DIAZ20063820}. Other examples often include certain crystals and anisotropic elastic media \cite{Komatitsch_etal2011,Duruthesis2012,DuKr,Be_etAl}. See for example the anisotropic material AM2 displayed in Figure \ref{fig:dispersion_relation} d). 
For advective acoustics or linearised  Euler equations with non-vanishing mean flow, a space-time linear transformation can allow for the derivation of stable PML, avoiding the catastrophic exponential growth for the Cauchy PML \cite{ABARBANEL1999266,HU2001455, Hagstrom2003,DIAZ20063820}.
There is however no cure yet for general systems such as linear anisotropic  elastodynamic equation.

In isotropic media and many anisotropic media the geometric stability condition is satisfied. In particular, it can be shown for isotropic elastic materials that  there are no growing modes for any frequency.
However, in numerical simulations instability may persist even when the mode analysis indicates that the Cauchy problem is stable. One possibility is that the introduction of boundary conditions can ruin the stability or well-posedness of the PML  \cite{SkAdCr}. A second possible source of numerical growth can result from unstable numerical approximation of the PML IBVP. Next we will review the theory of well-posedness and stability for PML IBVP.




\subsubsection{Stability analysis of PML IBVP}
The mode analysis for PMLs has been extended to include boundary conditions, see \cite{Duruthesis2012,DuKrSIAM,KDuru2016,DuruKozdonKreiss2016}. The central result is that the PML IBVP will not support exponentially growing solutions as long as neither the underlying undamped IBVP nor the PML Cauchy problem support exponentially growing solutions.  
We will summarise the analysis below. For more elaborate discussion we refer the reader to the references above.


As before we assume constant coefficients and consider a PML with damping only in the horizontal $(x)$ direction. We analyze a left half-space problem with a boundary at $x=0$, and a lower half-space problem with a boundary at $y=0$, seperatly. In both cases the boundary condition is of the same form as in the original equation. We will look for simple wave like solutions 
\begin{align}\label{eq:wave-like_IBVP2a}
\mathbf{U} = \boldsymbol{\phi}_0(x)e^{ik_y y + ik_z z + s t} ,  \quad \int_{-\infty}^{0}|\boldsymbol{\phi}_0(x)|^2 dx < \infty,
\end{align}
and
\begin{align}\label{eq:wave-like_IBVP2b}
\mathbf{U} = \boldsymbol{\phi}_0(y)e^{ik_x x + ik_z z + s t} ,\quad \int_{-\infty}^{0}|\boldsymbol{\phi}_0(y)|^2 dy < \infty,
\end{align}
respectively, where  $s = \Re{s} + i\Im{s}$.


We consider only physical media where the \emph{geometric stability condition} is satisfied.

\paragraph{PML IBVPs.}
 We consider the PML in a domain with boundaries and boundary conditions, that is  2a) left half-space problem, and 2b) lower half-space problem.
 
As above the PML is not symmetric, we consider first the boundary terminating the PML, to be precise \\
2a) The left half-spce PML problem:  $(x, y, z) \in  (-\infty, 0) \times (-\infty, \infty)\times (-\infty, \infty)$, with $\mathbf{L}\mathbf{U} = 0$,  $x =  0$.
We look for simple wave-like solutions of the form

Next, we consider \\
2b) The lower half-space PML problem:  $(x, y, z) \in  (-\infty, \infty) \times (-\infty, 0)\times (-\infty, \infty)$, with $\mathbf{L}\mathbf{U} =0$,  $y =  0$.
We also  look for simple wave solutions of the form

\begin{definition}
The IBVP PML  2a) or 2b) is not stable if there are nontrivial solutions $\mathbf{U}$ of the form  \eqref{eq:wave-like_IBVP2a} or \eqref{eq:wave-like_IBVP2b} with $\Re{s} > 0$.
\end{definition}

We make the following technical assumption
\begin{assumption}\label{assumpt:technical_assumption}
In the absence of  the PML, $d_x  = 0$, the IBVPs are stable, that is  \emph{all nontrivial eigen-pairs are in the stable half of the complex plane, $ \Re{s} \le 0$}. 
\end{assumption}
Note that the underlying IBVPs without the PML, $d_x =0$, satisfy strictly non-growing energy estimates, \eqref{eq:cauchy_estimate_first_0} and \eqref{eq:cauchy_estimate_second_0}. Therefore, Assumption \ref{assumpt:technical_assumption} is satisfied, otherwise the energy in the medium will grow, thereby contradicting the energy estimate. The results below would determine how the eigenvalues will be perturbed {when the PML is present, $d_x >0$}.
%
%
%
%
\begin{theorem}\label{theo:IBVP_problem_2a}
  Consider the constant coefficient PML IBVP  2a) with $d_x =d\ge 0$ subject to the boundary condition  $\mathbf{L}\mathbf{U} = 0 $ at $x =  0$. If Assumption \ref{assumpt:technical_assumption} holds and the geometric stability condition is satisfied,  then  the IBVP has no nontrivial  solution  of the form \eqref{eq:wave-like_IBVP2a} with $\Re{s}>0$.
\end{theorem}
\begin{theorem}\label{theo:IBVP_problem_2b}
Consider the constant constant coefficient lPML IBVP  2b)   with $d_x = d\ge0$,  subject to the boundary condition  $\mathbf{L}\mathbf{U} = 0 $ at $y =  0$.  If Assumption \ref{assumpt:technical_assumption} holds and the geometric stability condition is satisfied then { the PML damping $d_x > 0$ will move all eigenvalues further into the stable complex plane}, $\Re{s}<0$.
\end{theorem}
The proofs of Theorems \ref{theo:IBVP_problem_2a} and  \ref{theo:IBVP_problem_2b} can be adapted from  \cite{KDuru2016,DuruKozdonKreiss2016,DuKrSIAM}.

Theorems \ref{theo:high_frequency_stability}--\ref{theo:IBVP_problem_2b}  together prove that if the geometric stability condition is satisfied the constant coefficient PML  IBVP does not support exponentially growing solutions. However these results are  { too technical to be extended to the analysis of numerical approximations of the PML IBVP}. 

Below, we will review  extensions of the results using the energy method. The energy estimates enable the development of stable and accurate numerical methods for the PML using SBP finite difference methods and DG methods. More detailed discussions can be found in \cite{KDuru2016,DuruKozdonKreiss2016,ElasticDG_PML2019,DuruGabKreiss2019}.
\begin{remark}
We remark that the analysis above and Theorems \ref{theo:high_frequency_stability}--\ref{theo:IBVP_problem_2b} hold for both first order systems and second order hyperbolic systems.
In particular some of these results, for example Theorem \ref{theo:IBVP_problem_2b}, were initially formulated for the 2D elastic wave equation in second order form \cite{DuKrSIAM}, however the results hold for first order systems as well \cite{DuruKozdonKreiss2016} and extend to 3D \cite{ElasticDG_PML2019,DuruGabKreiss2019}.
\end{remark}

\begin{remark}
We also note that Theorems \ref{theo:high_frequency_stability}--\ref{theo:IBVP_problem_2b} proves that the PML is stable in the presence of boundary waves such as Rayleigh waves propagating on the surface of elastic solids.
The proof of stability of the PML in the presence  of interface waves such as Stoneley waves and Scholte waves is still an open problem. However, we believe that a similar technique for boundary waves can be used to establish stability for interface wave modes.
\end{remark}
\subsection{Energy methods}
The PML is generally asymmetric. It becomes extremely difficult to derive energy estimates that can be used to design stable numerical methods for the PML in truncated domains. 
In many settings, a straight forward approach would yield an exponentially growing estimate in physical space. This would indicate well-posedness, but exponential growth is not optimal for an absorbing model. A well-functioning PML should not support exponential growth.

We note that at constant coefficient an energy estimate for the acoustic wave equation in second order form became possible \cite{Baffet2019} for a specific time-domain formulation of the PML.  However, we will  briefly review the energy method in the Laplace space in general media, which is applicable to most PML models and useful for developing  stable and accurate numerical methods.

To begin with we reformulate the IBVP by introducing new variables by subtracting the product of the initial function and $e^{-t}$ from each of the unknown functions. The new unknown functions satisfy the same system as the original unknowns but with homogeneous initial data and inhomogeneous source terms in the equations. Denote the source terms by $\mathbf{F}_Q(x,y,z, t) =(\mathbf{F}_U(x,y,z, t), \mathbf{F}_{w_\xi}(x,y,z, t))^T.$
Laplace transformation in time  of the PML equations \eqref{eq:elastic_pml_1A}-\eqref{eq:elastic_pml_2A},  yields

 {
\begin{equation}\label{eq:elastic_pml_1A_Laplace}
\begin{split}
\mathbf{P}^{-1} s\widetilde{\mathbf{U}} = \sum_{\xi = x, y, z}\left[\frac{1}{\gamma_\xi}\mathbf{A}_{\xi}\frac{\partial{\widetilde{\mathbf{U}}}}{\partial \xi}  - d_{\xi}\left(\xi\right)\widetilde{\mathbf{w}}_{\xi}\right] + \widetilde{\mathbf{F}}_{U},
  \end{split}
  \end{equation}
\begin{equation}\label{eq:elastic_pml_2A_Laplace}
\begin{split}
s\widetilde{\mathbf{w}}_{\xi}
 = \frac{1}{\gamma_\xi}\mathbf{A}_{\xi}\frac{\partial{\widetilde{\mathbf{U}}}}{\partial \xi} -   \left(\alpha_{\xi}\left(\xi\right) + d_{\xi}\left(\xi\right)\right)\widetilde{\mathbf{w}}_{\xi} +  \widetilde{\mathbf{F}}_{w_\xi},
  \end{split}
  \end{equation}
}
Further, we use \eqref{eq:elastic_pml_2A_Laplace} and eliminate the auxiliary variables $\widetilde{\mathbf{w}}_{\xi}$ in \eqref{eq:elastic_pml_1A_Laplace}.  We obtain
\begin{equation}\label{eq:elastic_pml_3D_Laplace}
\begin{split}
\mathbf{P}^{-1} s\widetilde{\mathbf{U}}=   \sum_{\xi = x, y, z}\frac{1}{S_{\xi}} \mathbf{A}_{\xi}\frac{\partial{\widetilde{\mathbf{U}} }}{\partial \xi}  + \mathbf{P}^{-1}\widetilde{\mathbf{F}}, \quad \widetilde{\mathbf{F}} =     \mathbf{P}\left( \widetilde{\mathbf{F}}_{U}  - \sum_{\xi = x,y, z} \frac{d_{\xi} }{s + \alpha_{\xi} + d_{\xi}} \widetilde{\mathbf{F}}_{w_\xi}\right).
  \end{split}
  \end{equation}
 Transforming the boundary conditions \eqref{eq:bc_first_order} and the interface conditions  \eqref{eq:interface_condition}  gives
\begin{equation}\label{eq:BC_General2_Laplace}
\begin{split}
L\widetilde{\mathbf{U}} = 0, \quad \xi = \pm 1, 
 \end{split}
\end{equation}
such that 
$$
\frac{1}{2}\widetilde{\mathbf{U}}^\dagger \mathbf{A}_{\xi}\widetilde{\mathbf{U}} \le 0,
$$
and 
\begin{align}\label{eq:interface_condition_laplace}
\mathbb{I}\mathbb{T}\left(\mathbf{A}_{\xi}, \widetilde{\mathbf{U}} \right) =0, \quad \xi = 0,
\end{align}
with   $\mathrm{IT}\left( \widetilde{\mathbf{U}}^\dagger\mathbf{A}_{\xi} \widetilde{\mathbf{U}} \right) =0$, and $\widetilde{\mathbf{U}}^{\dagger}$ denotes the conjugate transpose of $\widetilde{\mathbf{U}}$.
 
 For Laplace transformed systems, it is necessary to extend the weighted scalar product \eqref{eq:weighted_scalar_product} and the corresponding norm to complex functions.
 For complex functions, we define the weighted scalar product
 \begin{align}\label{eq:elastic_scalarproduct_laplacePML}
\left(\widetilde{\mathbf{U}}, \widetilde{\mathbf{F}}\right)_{{P}} =  \int_{\Omega} \frac{1}{2}\left[\widetilde{\mathbf{U}}^{\dagger}\mathbf{P}^{-1} \widetilde{\mathbf{F}}\right]   dxdydz, \quad \|\widetilde{\mathbf{U}} (\cdot ,\cdot ,\cdot , s)\|_{P}^2 = \left(\widetilde{\mathbf{U}}, \widetilde{\mathbf{U}}\right)_{{P}}.
 \end{align}
Again, $\widetilde{\mathbf{U}}^{\dagger}$ denotes the complex conjugate transpose of $\widetilde{\mathbf{U}}$. 

 As before we will consider the acoustic wave equation and proceed later to the general linear hyperbolic system system.
  We will now review the results presented in \cite{ElasticDG_PML2019,DuruGabKreiss2019}
\subsubsection{Energy estimate of the PML for the acoustic wave equation}
In \cite{DuruGabKreiss2019}, the PML for the acoustics wave equation was considered with the standard metric parameters, $\gamma_\xi = 1$, $\alpha = 0$ and $d_\xi(\xi) \ge 0$. 
The estimate can be extended to the general metric parameters \eqref{eq:PML_metric} used in this review. Consider specifically the acoustic wave equation with  $\mathbf{U} = \left(p,\mathbf{v}\right)^T$, having
  \begin{equation}\label{eq:acoustic_pml_3D_Laplace_00}
\begin{split}
\frac{1}{\kappa} {s \widetilde{p}}  + {\grad}_{d}  \cdot \widetilde{\mathbf{v}}= \frac{1}{\kappa }\widehat{F}_p, \quad
\rho {s \widetilde{\mathbf{v}} } + {\grad}_{d} \widetilde{p} = \rho  \mathbf{  \widehat{f} }, 
  \end{split}
  \end{equation}
   where  $ {\grad}_{d} = \left(1/S_{x} \partial/\partial{{x}}, 1/S_{y} \partial/\partial{{y}}, 1/S_{z} \partial/\partial{{z}}\right)^T$ where $\widehat{F}_p$ and $\mathbf{  \widehat{f} }$ are source terms that depend on the initial condition.
Further, we eliminate the velocity fields, in \eqref{eq:acoustic_pml_3D_Laplace_00}, having
 \begin{equation}\label{eq:acoustic_wave_pml_3D_Laplace}
\begin{split}
\frac{1}{\kappa} {s^2 \widetilde{p}}  - {\grad}_{d}  \cdot  \left(\frac{1}{\rho}{\grad}_{d} \widetilde{p}\right)&= \frac{s}{\kappa } \widehat{F}_p  - {\grad}_{d} \cdot \mathbf{\widehat{f}}, \\ 
  \end{split}
  \end{equation}

with the boundary conditions at ${\xi} = \pm 1$
\begin{equation}\label{eq:boundary_condition_acoustic_laplace}
\frac{1+r_{\xi}}{2Z}s\widetilde{p}  \pm \frac{1-r_{\xi}}{2}\frac{1}{S_{\xi}} \frac{\partial \widetilde{p}}{\partial \xi} = 0, \quad \text{with} \quad |r_{\xi}| \le 1.
 \end{equation}
The jump condition \eqref{eq:elastic_acoustic_interface_first} and \eqref{eq:interface_condition_Laplace} translate to
 \begin{align}\label{eq:jump_x_laplace}
\widetilde{p}^{-} = \widetilde{p}^{+} = \widetilde{p},  \quad \lJump \frac{1}{S_{\xi} \rho} \frac{\partial \widetilde{p}}{\partial \xi}  \rJump = 0.
 \end{align}

For   $s \ne 0$ and $d_{\xi} \ge 0$, 
denote $s = a + i b$ and
{
\begin{align}\label{eq:scaled_variables_1}
   \Re\left(\frac{\left(sS_{\xi}\right)^*}{S_{\xi}}\right) =  a + \epsilon_{\xi}(s, d_{\xi}), \quad \epsilon_{\xi}(s, d_{\xi}) :=  d_{\xi} b^2\frac{(a+\alpha)^2 + b^2 + (a+\alpha)a + \alpha d_\xi  }{\left((a+\alpha)(\alpha + a + d_\xi )\right)^2 + (d_\xi b)^2} \ge 0,
\end{align} 
}
with
  \begin{align}
  \Re\left(\frac{1}{S_{\xi}}\right) =  \frac{ (a+\alpha)(a+\alpha+d_{\xi})+ b^2}{\gamma_\xi\left( \left(a+\alpha+d_{\xi}\right)^2+b^2\right)} > 0.
  \end{align}
Define the energies 
\begin{equation}\label{eq:acoustic_pml_3D_energy}
\begin{split}
&\widetilde{E}^2_p\left(s, d_{\xi}\right) 
 = \Big\|{s}  \widetilde{p}\Big\|^2_{1/\kappa}+\sum_{\xi = x,y,z} \Big\|\frac{1}{S_{\xi}} \frac{\partial \widetilde{p} }{\partial {\xi}}\Big\|^2_{1/\rho} > 0, \\
&\widetilde{E}_f^2\left(s, d_{\xi}\right)  =  \Big\|{s}   \widehat{F}_p \Big\|^2_{1/\kappa}+\sum_{\xi = x,y,z} \Big\| \frac{1}{S_{\xi}}\frac{\partial \widehat{f}_{\xi}}{\partial {\xi}}\Big\|_{\kappa}^2 > 0.
  \end{split}
  \end{equation}
  The following theorem and corollary were proven in \cite{DuruGabKreiss2019}
\begin{theorem}\label{Theo:Stability_PML_Laplace_Acoustics}
Consider  the PML equation  in the Laplace space \eqref{eq:acoustic_wave_pml_3D_Laplace} with  constant damping  $d_{\xi}  \ge 0$,  and source terms  $\widetilde{\mathbf{F}}(x,y,z, s)$, $\widetilde{\mathbf{f}} $,    subject to  homogeneous initial data,  the boundary conditions \eqref{eq:boundary_condition_acoustic_laplace},  and the jump condition \eqref{eq:jump_x_laplace} at discontinuities in material parameters.  
Let  the energy norms  $\widetilde{E}_p\left(s, d_{\xi}\right) > 0$, $\widetilde{E}_f\left(s, d_{\xi}\right) > 0$ in the Laplace space be defined in \eqref{eq:acoustic_pml_3D_energy}. 
%
For any $\Re{s}=a>0$ we have
\begin{equation}\label{eq:energy_estimate_pml_laplace_acoustic_cont}
\begin{split}
&a\widetilde{E}_p^2\left(s, d_{\xi}\right) + \sum_{\xi = x,y,z} \Big\|\frac{1}{S_{\xi}} \frac{\partial \widetilde{p} }{\partial {\xi}}\Big\|^2_{\epsilon_{\xi}/\rho}  + \widetilde{\mathrm{BT}}\left(s, d_{\xi}\right)    \le   2 \widetilde{E}_p\left(s, d_{\xi}\right)  \widetilde{E}_f\left(s, d_{\xi}\right), \\
& \widetilde{\mathrm{BT}}  =   \int_{-1}^{1}\int_{-1}^{1}\sum_{\xi = x,y,z}\Re\left(\frac{1}{S_{\xi}}\right)\left(\frac{s^*}{S_{\xi} \rho} \widetilde{p}^*\frac{\partial \widetilde{p}}{\partial \xi}\Big|_{\xi =-1}- \frac{s^*}{S_{\xi} \rho} \widetilde{p}^*\frac{\partial \widetilde{p}}{\partial \xi}\Big|_{\xi =1}   \right)\frac{dxdydz}{d\xi}  \ge 0.
  \end{split}
\end{equation}
\end{theorem}

\begin{corollary}\label{eq:remark_homogenous_continuous}
Consider the same PML  IBVP,  as in Theorem \ref{Theo:Stability_PML_Laplace_Acoustics}. The IBVP, with $d_{\xi} \ge 0$, is asymptotically stable in the sense that no exponentially growing solutions are supported.
\end{corollary}
\subsubsection{Energy estimate of the PML for general first order systems}
The discussions here are the contributions from \cite{ElasticDG_PML2019}. For  general systems we consider the scaled equation
\begin{equation}\label{eq:elastic_pml_3D_Laplace_source}
\begin{split}
&\mathbf{P}^{-1} \left(sS_x\right)\widetilde{\mathbf{U}} =   \sum_{\xi = x, y, z}\frac{S_x}{S_{\xi}} \mathbf{A}_{\xi}\frac{\partial{\widetilde{\mathbf{U}} }}{\partial \xi}  + \mathbf{P}^{-1}\widetilde{\mathbf{F}}, 
\\
&\widetilde{\mathbf{F}} =     \mathbf{P}\left(S_x \widetilde{\mathbf{F}}_{U}  - \sum_{\xi = x,y, z} \frac{d_{\xi}S_x }{s + \alpha_{\xi} + d_{\xi}} \widetilde{\mathbf{F}}_{w_\xi}\right),
  \end{split}
  \end{equation}
with  the boundary conditions \eqref{eq:BC_General2_Laplace}.
{  The form \eqref{eq:elastic_pml_3D_Laplace_source}  is chosen so that with suitable assumptions the spatial operator remains  in its original form, and allows  using the standard integration by parts methodology.}
For an energy estimate we also need the coefficient $\left(sS_x\right)$ in the left hand side to have a positive real part. Not that
\begin{align*}
&\Re{\left(sS_x\right)} =  {\gamma_x}\left(a + \left( \frac{a\left(a+\alpha\right) + b^2}{{|s + \alpha|^2}}\right)d_x\right) \ge a \gamma_x, \\
&0 < \frac{1}{\Re{\left(sS_x\right)}} = \frac{1}{\gamma_x\left(a + \left( \frac{a\left(a+\alpha\right) + b^2}{{|s + \alpha|^2}}\right)d_x \right)} \le \frac{1}{a\gamma_x}.
\end{align*}
\paragraph{The PML strip problem.}
Consider specifically the $x$-dependent  PML strip problem, that is  \eqref{eq:elastic_pml_1A}-\eqref{eq:elastic_pml_2A}  with  $d_x \ge  0$, $\gamma_x = \gamma >0$  and $\partial /\partial \xi =  0$, $d_\xi = 0$, $\gamma_\xi = 1$ for $\xi = y, z$.
This simplification results in the 1D PML problem
\begin{equation}\label{eq:weak_Laplace_1D}
\begin{split}
 \left(sS_x\right) \mathbf{P}^{-1} \widetilde{\mathbf{U}} = {A}_x\frac{\partial{\widetilde{\mathbf{U}} }}{\partial x} +  \mathbf{P}^{-1}\widetilde{\mathbf{F}}.
  \end{split}
  \end{equation}
  Equation \eqref{eq:weak_Laplace_1D} lives in a 3D domain defined in \eqref{eq:physical_domain}, but has been simplified by restricting the initial data and the forcing $\widetilde{\mathbf{F}}$  to functions that vary only in 1D, the $x$-axis.
  The following theorem was proven in \cite{ElasticDG_PML2019}.
\begin{theorem}\label{theo:PML_1D_Strip_Laplace}
Consider the 1D   PML equation \eqref{eq:weak_Laplace_1D} in the Laplace space, with  piecewise constant $d_x(x) = d\ge 0$, constant grid stretching $\gamma_x = \gamma >0$ and  $\alpha_x  = \alpha\ge 0$,  subject to the boundary conditions \eqref{eq:BC_General2_Laplace}, and the interface condition \eqref{eq:interface_condition_laplace} at discontinuities of $d_x(x)$. For all $s$ such that  $\Re(s) \ge a > 0$ we have
{
\begin{equation}\label{eq:energy_estimate_pml_laplace_strip}
\begin{split}
& \|\sqrt{\Re(s S_x)}\widetilde{\mathbf{U}} (\cdot ,\cdot ,\cdot , s)\|_{P}^2 \le \|\widetilde{\mathbf{U}} (\cdot ,\cdot ,\cdot , s)\|_{P} \|\widetilde{\mathbf{F}} (\cdot ,\cdot ,\cdot , s)\|_{P} + \widetilde{\mathrm{BT}}, \quad  \Re(s S_x) \ge a\gamma > 0, \\
& \widetilde{\mathbf{F}} =     \mathbf{P}\left(S_x \widetilde{\mathbf{F}}_{U}  -  \frac{d }{s + \alpha} \widetilde{\mathbf{F}}_{w_x}\right), \quad \widetilde{\mathrm{BT}} = \int_{\widetilde{\Gamma}}\frac{1}{2}\left[\widetilde{\mathbf{U}}^\dagger{A}_x\widetilde{\mathbf{U}}\right]\Big|_{x = -1}^{x =1}dydz  \le 0.
\end{split}
\end{equation}
}
\end{theorem}

\begin{remark}\label{remark:pml_strip}
It is important to note that   similar energy  estimates \eqref{eq:energy_estimate_pml_laplace} are also valid for the  PML strip problems in the $y$-axis or $z$-axis, that is  when $ d_y = d>0, d_x = d_z=0$ or  $ d_z = d>0, d_x = d_y=0$.
\end{remark}

\paragraph{The PML  edge problem.}
Next, we consider specifically the  $xy$--edge PML  problem, that is \eqref{eq:elastic_pml_1A}-\eqref{eq:elastic_pml_2A} with $d_x  = d_y = d>0$, $\gamma_x =\gamma_y = \gamma >0$,  $\alpha_x = \alpha_y = \alpha\ge 0$ and $d_z = 0$, $\partial /\partial z =0$.  
This simplification results in  the 2D $xy$-edge PML problem, 
\begin{equation}\label{eq:weak_Laplace_2D}
\begin{split}
 \left(sS_x\right) \mathbf{P}^{-1} \widetilde{\mathbf{U}}  =  \sum_{\xi = x, y}\mathbf{A}_{\xi}\frac{\partial{\widetilde{\mathbf{U}} }}{\partial \xi}  + \mathbf{P}^{-1}\widetilde{\mathbf{F}}.
  \end{split}
  \end{equation}
Note also that equation \eqref{eq:weak_Laplace_2D} lives in a 3D domain defined in \eqref{eq:physical_domain}, but has been simplified by restricting the initial data and the forcing $\widetilde{\mathbf{F}}$  to functions that vary only in 2D, the $xy$-plane. We have the result which was proven in \cite{ElasticDG_PML2019}.
\begin{theorem}\label{theo:PML_2D_Edge_Laplace}
Consider the 2D   PML equation \eqref{eq:weak_Laplace_2D} in the Laplace space, with  piecewise constant $d_x(x) = d_y(y) = d\ge 0$, constant grid stretching $\gamma_x =\gamma_y = \gamma >0$  and  $\alpha_x = \alpha_y = \alpha\ge 0$,  subject to the boundary conditions \eqref{eq:BC_General2_Laplace},   and the interface condition \eqref{eq:interface_condition_laplace} at discontinuities of $d_x(x)$, $d_y(y)$. For all $s$ such that  $\Re(s) \ge a> 0$ we have
{
\begin{equation}\label{eq:energy_estimate_pml_laplace_edge}
\begin{split}
&\|\sqrt{\Re(s S_x)}\widetilde{\mathbf{U}} (\cdot ,\cdot ,\cdot , s)\|_{P} \le \|\widetilde{\mathbf{U}} (\cdot ,\cdot ,\cdot , s)\|_{P} \|\widetilde{\mathbf{F}} (\cdot ,\cdot ,\cdot , s)\|_{P} + \widetilde{\mathrm{BT}},  \quad \Re(s S_x) \ge a\gamma > 0,\\
&  \widetilde{\mathbf{F}} =     \mathbf{P}\left(S_x \widetilde{\mathbf{F}}_{U}  -  \sum_{\xi = x, y}\frac{d }{s + \alpha} \widetilde{\mathbf{F}}_{w_\xi}\right) , 
\quad
\widetilde{\mathrm{BT}} = \sum_{\xi = x, y}\int_{\widetilde{\Gamma}} \frac{1}{2 }\left[\widetilde{\mathbf{U}}^\dagger\mathbf{A}_{\xi} \widetilde{\mathbf{U}}\right]\Big|_{-1}^{1}\frac{dxdydz}{d\xi}  \le 0.
\end{split}
\end{equation}
}
\end{theorem}

   \begin{remark}\label{remark:pml_edge}
It is also noteworthy   that  similar energy  estimates \eqref{eq:energy_estimate_pml_laplace_edge} are valid for the  PML edge problems in the $xz$-edge or $yz$-edge, that is  when $ d_x = d_z = d>0, d_y=0$ or  $ d_y = d_z = d>0, d_x =0$.
\end{remark}

\paragraph{The PML  corner problem.}

Consider the corresponding PML corner problem,  \eqref{eq:elastic_pml_1A}-\eqref{eq:elastic_pml_2A}   with $d_x  = d_y = d_z = d>0$, $\gamma_x =\gamma_y = \gamma_z = \gamma >0$, and $\alpha_x  = \alpha_y = \alpha_z = \alpha >0$.
Then, all PML metrics are identical $S_y = S_x = S_z$. We have $S_x/S_{\xi} = 1$ and 
\begin{equation}\label{eq:weak_Laplace_Corner}
\begin{split}
\left(sS_x\right) \mathbf{P}^{-1} \widetilde{\mathbf{U}}  =  \sum_{\xi = x, y,z}\mathbf{A}_{\xi}\frac{\partial{\widetilde{\mathbf{U}} }}{\partial \xi}  +  \mathbf{P}^{-1}\widetilde{\mathbf{F}} .
  \end{split}
  \end{equation}
\begin{theorem}\label{theo:PML_3D_Corner_Laplace}
Consider the 3D   PML equation \eqref{eq:weak_Laplace_Corner} in the Laplace space, with  piecewise constant $d_x(x) = d_y(y) = d_z(z) =d\ge 0$, constant grid stretching $\gamma_x =\gamma_y=\gamma_z = \gamma >0$  and  $\alpha_x = \alpha_y =  \alpha_z =  \alpha\ge 0$,  subject to the boundary conditions \eqref{eq:BC_General2_Laplace}, and the interface condition \eqref{eq:interface_condition_laplace} at discontinuities of $d_x(x)$, $d_y(y)$, $d_z(z)$. For all $s$ such that  $\Re(s) \ge a > 0$ we have
%
\begin{equation}\label{eq:energy_estimate_pml_laplace_corner}
\begin{split}
&\|\sqrt{\Re(s S_x)}\widetilde{\mathbf{U}} (\cdot ,\cdot ,\cdot , s)\|_{P}^2 \le \|\widetilde{\mathbf{U}} (\cdot ,\cdot ,\cdot , s)\|_{P} \|\widetilde{\mathbf{F}} (\cdot ,\cdot ,\cdot , s)\|_{P} + \widetilde{\mathrm{BT}},  \quad \Re(s S_x)\ge a\gamma > 0,\\
&  \widetilde{\mathbf{F}} =     \mathbf{P}\left(S_x \widetilde{\mathbf{F}}_{U}  -   \sum_{\xi = x, y,z}\frac{d }{s + \alpha}\widetilde{\mathbf{F}}_{w_\xi}\right),  
\quad
\widetilde{\mathrm{BT}} = \sum_{\xi = x, y, z}\int_{\widetilde{\Gamma}} \frac{1}{2}\left[\widetilde{\mathbf{U}}^\dagger\mathbf{A}_{\xi} \widetilde{\mathbf{U}}\right]\Big|_{-1}^{1}\frac{dxdydz}{d\xi}  \le 0. 
\end{split}
\end{equation}
\end{theorem}

We introduce the energy norms  in the physical space 
\begin{equation}\label{eq:elastic_pml_3D_energy_physical}
\begin{split}
 \|{\mathbf{U}} (\cdot ,\cdot ,\cdot , t)\|_{P}^2  = \|\mathcal{L}^{-1}\widetilde{\mathbf{U}} (\cdot ,\cdot ,\cdot , s)\|_{P}^2, \quad
  \|{\mathbf{F}} (\cdot ,\cdot ,\cdot , t)\|_{P}^2  = \|\mathcal{L}^{-1}\widetilde{\mathbf{F}} (\cdot ,\cdot ,\cdot , s)\|_{P}^2 ,
  \end{split}
  \end{equation}
\begin{theorem}\label{Theo:estimate_PML}
Consider the energy estimate in the Laplace space
\begin{equation}\label{eq:energy_estimate_pml_laplace}
\begin{split}
& \|\widetilde{\mathbf{U}} (\cdot ,\cdot ,\cdot , s)\|_{P}^2 \le \frac{1}{a\gamma }\|\widetilde{\mathbf{U}} (\cdot ,\cdot ,\cdot , s)\|_{P} \|\widetilde{\mathbf{F}} (\cdot ,\cdot ,\cdot , s)\|_{P}, \quad \Re(s) \ge a  > 0.
\end{split}
\end{equation}
For any $a>0$, $\gamma >0$ and $T>0$ we have
  \begin{align}\label{eq:estimate_0_0}
  \int_0^{T} e^{-2at}  \|{\mathbf{U}} (\cdot ,\cdot ,\cdot , t)\|_{P}^2 dt \le  \frac{1}{\left(a\gamma\right)^2}\int_0^{T} e^{-2at} \left( \| {\mathbf{F}} (\cdot ,\cdot ,\cdot , t)\|_{P}^2\right) dt.
  \end{align}
\end{theorem}
Since the boundary terms $\widetilde{\mathrm{BT}}$ in
\eqref{eq:energy_estimate_pml_laplace_acoustic_cont},  \eqref{eq:energy_estimate_pml_laplace_corner}, 
\eqref{eq:energy_estimate_pml_laplace_edge} 
and 
\eqref{eq:energy_estimate_pml_laplace_strip} are never positive, we can use Theorem \ref{Theo:estimate_PML} to invert the the estimates in Theorem \ref{Theo:Stability_PML_Laplace_Acoustics}, 
\ref{theo:PML_1D_Strip_Laplace}, 
\ref{theo:PML_2D_Edge_Laplace} 
and 
\ref{theo:PML_3D_Corner_Laplace}, 
and 
get an energy estimate in the physical space.
The estimate seemingly allows the energy of the solution to grow  exponentially with time for general data. However, by choosing $a>0$ in relation to the length 
of the time interval of interest, a bound  involving only algebraic growth in time follows.

\begin{remark}
For general systems the energy estimates for the first order formulation, Theorems \ref{theo:PML_1D_Strip_Laplace}--\ref{theo:PML_3D_Corner_Laplace}, are not directly applicable to general second order formulation, in particular when mixed spatial derivatives are present. For general second order systems, such as displacement formulation of linear elastodynamics, we may need other simplifying assumptions and model reductions in order to formulate results that will enable the development of accurate and robust numerical methods. 
\end{remark}

\section{Numerical  analysis of the discrete PML}\label{sec:s5}
In this section, we review the derivation of stable numerical methods and numerical analysis of the PML. We will formulate the ideas in a general DG framework.
However, the results extend to finite difference approximation based on the SBP-SAT method.
As shown in \cite{KDuru2016,DuruKozdonKreiss2016,ElasticDG_PML2019,DuruGabKreiss2019}, the goal is to design the numerical method for the PML IBVP in a manner that allows the derivation of discrete energy estimates analogous to the continuous energy estimates presented in Theorems \ref{Theo:Stability_PML_Laplace_Acoustics}--\ref{Theo:estimate_PML}.

As presented in \cite{ElasticDG_PML2019,DuruGabKreiss2019} we will  use the physically motivated numerical fluxes developed in \cite{DuruGabrielIgel2017,Duru_exhype_2_2019} to patch DG elements together into the global domain.  For the undamped problem \eqref{eq:first_order_linear_hyp_pde}, the physically motivated numerical flux is upwind by construction and gives an energy estimate analogous to \eqref{eq:cauchy_estimate_first_0}.  We believe that the general technique can be extended other DG methods with a different but stable numerical flux. 

\subsection{Integral formulation of the boundary and inter-element procedures}
Consider  the spatial domain $(x,y,z)  \in \Omega = [-1, 1]\times[-1, 1]\times[-1, 1]$ and discretise it into $K\times L\times M$ elements, where the $klm$-th element is denoted by $\Omega_{klm} = [x_k, x_{k+1}]\times [y_l, y_{l+1}]\times [z_m, z_{m+1}]$, for $k = 1, 2, \dots, K$, $l = 1, 2, \dots, L$, $m = 1, 2, \dots, M$ with $x_1 = -1$,  $y_1 = -1$, $z_1 = -1$ and $x_{K+1} = 1$, $y_{L+1} = 1$, $z_{M+1} = 1$.
%
The physical element $\Omega_{lmn} = [x_k, x_{k+1}]\times [y_l, y_{l+1}]\times [z_m, z_{m+1}]$ is mapped to a reference element $(q, r, s) \in \widetilde{\Omega} = [-1, 1]^{3}$ using the linear transformation 
{
\begin{align}\label{eq:transf}
x = x_k + \frac{\Delta{x}_k}{2}\left(1 + q \right),  \quad
y = y_l + \frac{\Delta{y}_l}{2}\left(1 + r \right), \quad 
z = z_m + \frac{\Delta{z}_m}{2}\left(1 + s \right),  
\end{align}
}
with
{
$ \Delta{x}_k = x_{k+1} - x_k, \quad  \Delta{y}_l = y_{l+1} - y_l, \quad  \Delta{z}_m = y_{m+1} - y_m.$
}
The nonzero metric derivatives and the Jacobian of the transformation are
$$
 q_x = \frac{2}{\Delta{x}_k}, \quad r_y = \frac{2}{\Delta{y}_l}, \quad s_z = \frac{2}{\Delta{z}_m}, \qquad J = \frac{\Delta{x}_k}{2}\frac{\Delta{y}_l }{2}\frac{\Delta{z}_m}{2} > 0.
$$
The elemental volume integral yields
\begin{align}
\int_{\Omega_{klm}} f(x,y,z) dxdydz = \int_{\widetilde{\Omega}} f(q,r,s) Jdqdrds,
\end{align}
We introduce  the scalar products defined by the volume and surface integrals
\begin{align}\label{eq:scalar_products}
\left(\boldsymbol{\phi}, \boldsymbol{\psi}\right) = \int_{\widetilde{\Omega}}\boldsymbol{\phi}^T\boldsymbol{\psi} Jdqdrds, \quad 
\left(\boldsymbol{\phi}, \boldsymbol{\psi}\right)_{\Gamma_{\xi}} = \int_{\widetilde{\Gamma}}\boldsymbol{\phi}^T\boldsymbol{\psi} J\frac{dqdrds}{d\xi},
\end{align}
and  adopt the notation 
\begin{align*}
& \mathbf{w}_{q} = \mathbf{w}_{x}, \quad \mathbf{w}_{r} = \mathbf{w}_{y}, \quad \mathbf{w}_{s} = \mathbf{w}_{z}, \quad \mathbf{A}_{q} = \frac{q_x}{\gamma_x} \mathbf{A}_{x} , \quad \mathbf{A}_{r} = \frac{r_y}{\gamma_y} \mathbf{A}_{y} , \quad \mathbf{A}_{s} = \frac{s_z}{\gamma_z} \mathbf{A}_{z}, \nonumber \\
& {d}_{q} =  {d}_{x} , \quad {d}_{r} =  {d}_{y} , \quad {d}_{s} =  {d}_{z}, \quad {\alpha}_{q} =  {\alpha}_{x} , \quad {\alpha}_{r} =  {\alpha}_{y} , \quad {\alpha}_{s} =  {\alpha}_{z}.
\end{align*}

    The elemental integral formulation reads
\begin{align}\label{eq:weak_01}
\left(\boldsymbol{\phi}, \mathbf{P}^{-1}\frac{\partial \mathbf{U}}{\partial t}\right) = \sum_{\xi = q, r, s}\left[\left(\boldsymbol{\phi}, \mathbf{A}_{\xi}\frac{\partial{\mathbf{U}}}{\partial \xi}  - d_{\xi}\left(\xi\right)\mathbf{w}_{\xi}\right) - \left(\left(\boldsymbol{\phi}\Big|_{\xi = -1}, \mathbf{FL}_{\xi} \right)_{\Gamma_{\xi}} + \left(\boldsymbol{\phi}\Big|_{\xi = 1}, \mathbf{FR}_{\xi} \right)_{\Gamma_{\xi}}\right)\right],
\end{align}
\begin{align}\label{eq:weak_11}
\left(\boldsymbol{\phi}, \frac{\partial{\mathbf{w}_{\xi}}}{\partial t}\right) = \left(\boldsymbol{\phi}, \mathbf{A}_{\xi}\frac{\partial{\mathbf{U}}}{\partial \xi} -   \left(\alpha_{\xi}(\xi) + d_{\xi}\left(\xi\right)\right)\mathbf{w}_{\xi}\right) - \underbrace{\theta_{\xi}\left(\left(\boldsymbol{\phi}\Big|_{\xi = -1}, \mathbf{FL}_{\xi} \right)_{\Gamma_{\xi}} + \left(\boldsymbol{\phi}\Big|_{\xi = 1}, \mathbf{FR}_{\xi} \right)_{\Gamma_{\xi}}\right)}_{\text{PML stabilizing flux fluctuations}},
\end{align}
where $\mathbf{FL}_{\xi}$ and $\mathbf{FR}_{\xi}$ are physics based flux fluctuation vectors designed to enforce the boundary and inter-element conditions. Please see \cite{DuruGabrielIgel2017,Duru_exhype_2_2019, ElasticDG_PML2019,DuruGabKreiss2019} for more details.
%
We note however for exact solutions that satisfy the IBVP the flux fluctuations vanish identially $\mathbf{FL}_{\xi} =\mathbf{FR}_{\xi} = 0$. 
The stabilising PML parameters $\theta_\xi$ to be determined by requiring stability of the discrete PML, in the sense corresponding to Theorems \ref{theo:PML_1D_Strip_Laplace}, \ref{theo:PML_2D_Edge_Laplace}, \ref{theo:PML_3D_Corner_Laplace} and \ref{Theo:estimate_PML}.
\begin{remark}
In \cite{DuKrSIAM,DuruKozdonKreiss2016,ElasticDG_PML2019,DuruGabKreiss2019} we have demonstrated  for DG/SBP approximation that appropriate choice of the stabilising PML parameters $\theta_\xi$ is necessary for the stability of the discrete PML for several systems. These important results will be summarised in the next sections.
\end{remark}
 \subsection{The DG/SBP approximation of the PML}
We follow a standard finite element procedure and approximate the elemental solution by a polynomial interpolant 
\begin{equation}\label{eq:variables_elemental}
\mathbf{U}\left(q,r,s, t\right) =  \sum_{i = 1}^{N+1} \sum_{j = 1}^{N+1}\sum_{k = 1}^{N+1}\mathbf{U}_{ijk}(t) \boldsymbol{\phi}_{ijk}\left(q,r,s\right),
\end{equation}
where $\mathbf{U}_{ijk}(t)$, are the elemental degrees of freedom to be determined, 
 and $ \boldsymbol{\phi}_{ijk}(q,r,s)$ are the $ijk$-th interpolating polynomials. 
 We consider tensor products of  nodal basis with $ \boldsymbol{\phi}_{ijk}(q,r,s) = \mathcal{L}_i(q)\mathcal{L}_j(r)\mathcal{L}_k(s)$,  where $\mathcal{L}_i(q)$, $\mathcal{L}_j(r)$, $\mathcal{L}_k(s)$, are one dimensional nodal interpolating Lagrange polynomials of degree $N$.
We will only use Gauss-type quadrature rules such  that for all polynomial integrand $f(\xi)$ of degree $\le 2N-1$, the corresponding one dimensional rule is exact, 
 $\sum_{m = 1}^{N+1} f(\xi_m)h_m = \int_{-1}^{1}f(\xi) d\xi.$

 We now make a classical  Galerkin approximation by choosing test functions  in the same space as the basis functions, so that the residual is orthogonal to the space of test functions. We rearrange the elemental degrees of freedom $[\mathbf{U}_{ijk}(t) ]$ row-wise as a vector, $\mathbf{U}(t)$,  of length $n_f(N+1)^d$ where $d = 3$ is the number of space dimensions, $n_f$ ($n_f = 4$ for acoustics and $n_f = 9$ for elasticity) is the number of fields in the physical variables. We have the semi-discrete approximation 
{
\begin{equation}\label{eq:disc_elastic_pml_1A}
\begin{split}
\mathbf{P}^{-1}\frac{ d {\mathbf{U}}}{d t} = \sum_{\xi = x, y, z}\left[\mathbf{A}_{\xi}\mathbf{D}_{\xi} \mathbf{U}  - \mathbf{d}_{\xi}\mathbf{w}_{\xi} - \mathbf{H}_{\xi}^{-1}{\left({\mathbf{e}_{\xi}(-1)} \mathbf{FL}_\xi +  {\mathbf{e}_{\xi}(1)} \mathbf{FR}_\xi  \right)}\right],
  \end{split}
  \end{equation}
\begin{equation}\label{eq:disc_elastic_pml_3A}
\begin{split}
\frac{d {\mathbf{w}_{\xi}}}{ d t} 
 = \mathbf{A}_{\xi}\mathbf{D}_{\xi} \mathbf{U}  - \left(\mathbf{d}_{\xi} + \boldsymbol{\alpha}_{\xi}\right)\mathbf{w}_{\xi} - \underbrace{\theta_{\xi}\mathbf{H}_{\xi}^{-1}{\left({\mathbf{e}_{\xi}(-1)} \mathbf{FL}_\xi +  {\mathbf{e}_{\xi}(1)} \mathbf{FR}_\xi  \right)}}_{\text{PML stabilizing flux fluctuations}}.
  \end{split}
  \end{equation}
}
Here the discrete operators are defined by
\begin{equation}\label{eq:derivative_operator_3d}
\mathbf{D}_{x} = \frac{2}{\Delta{x}_k}\left(I_{f} \otimes D\otimes I\otimes I\right), \quad \mathbf{D}_y = \frac{2}{\Delta{y}_l}\left(I_{f} \otimes  I\otimes D\otimes I\right), \quad \mathbf{D}_z = \frac{2}{\Delta{z}_m}\left(I_{f} \otimes  I\otimes I\otimes D\right),
\end{equation}
and the discrete norms are given by
{\small
\begin{align*}
\mathbf{H}_x = \frac{\Delta{x}_k}{2}\left(I_{f} \otimes  H\otimes I\otimes I\right), \quad \mathbf{H}_y = \frac{\Delta{y}_l}{2}\left(I_{f} \otimes  I\otimes H\otimes I\right), \quad \mathbf{H}_z = \frac{\Delta{z}_m}{2}\left(I_{f} \otimes  I\otimes I\otimes H\right),
\end{align*}
}
where
{
\begin{equation}\label{eq:derivative_operator_1d}
 D = H^{-1} Q, \quad H = \mathrm{diag}[h_1, h_2, \cdots, h_{P+1}], \quad Q_{ij} = \sum_{m = 1}^{P+1} h_m \mathcal{L}_i(q_m)  {\mathcal{L}_j^{\prime}(q_m)},
\end{equation}
}
is a  spectral difference approximation of the first derivative, in one space dimension.
We also have
\begin{equation*}
\mathbf{P} = \left({P} \otimes  I\otimes I\otimes I\right), \quad \mathbf{A}_{\xi} = \left({A}_{\xi} \otimes  I\otimes I\otimes I\right),  \quad \mathbf{H} = \mathbf{H}_x\mathbf{H}_y\mathbf{H}_z,
\end{equation*}
where $I$ is the $(N+1)\times(N+1)$ identity matrix, $I_{f} $ is the $n_f \times n_f$ identity matrix, and $\otimes$ denotes the Kronecker product. 

We also introduce the projection matrices
\begin{align*}
&\mathbf{e}_x(\eta) = \left(I_{n_f} \otimes  \boldsymbol{e}(\eta)\otimes I\otimes I\right), \quad \mathbf{e}_y(\eta) = \left(I_{n_f} \otimes  I\otimes \boldsymbol{e}(\eta)\otimes I\right), 
\\
&\mathbf{e}_z(\eta) = \left(I_{n_f} \otimes  I\otimes I\otimes \boldsymbol{e}(\eta)\right), \quad  \mathbf{B}_{\eta}(\psi, \xi) = \mathbf{e}_{\eta}(\psi) \mathbf{e}_{\eta}^T(\xi),
\end{align*}
where
\[
  \boldsymbol{e}(\eta) = [\mathcal{L}_i(\eta), \mathcal{L}_i(\eta), \cdots, \mathcal{L}_{P+1}(\eta)]^T.
\]
The derivative oopertaors \eqref{eq:derivative_operator_3d} satisfy the discrete SBP property
\begin{align}\label{eq:disc_sbp}
 \mathbf{D}_\xi = -\mathbf{H}_\xi^{-1}\mathbf{D}_\xi^T\mathbf{H}_\xi + \mathbf{H}_\xi^{-1}\left(\mathbf{B}_\xi\left(1,1\right)-\mathbf{B}_\xi\left(-1,-1\right)\right), \quad \xi = x, y, z.
\end{align}
\begin{remark}
Note that the spatial discrete operators \eqref{eq:derivative_operator_1d}, \eqref{eq:derivative_operator_3d} satisfying the SBP property \eqref{eq:disc_sbp} have been derived using a spectral approach and have full accuracy, (that is $(P+1)$th accuracy) within the element. Standard SBP finite difference operators on equidistant grids have also been used with reduced accuracy close to the boundaries. However, the stability results are identical.
\end{remark}
Introduce the discrete inner products and norm
\begin{align*}
\Big\langle {\mathbf{U}},  {\mathbf{V}}\Big\rangle_{\mathbf{H}} = {\mathbf{V}}^{T} \mathbf{H} {\mathbf{U}}, \quad \|{\mathbf{U}}\|_{\mathbf{H}}^2 = \Big\langle {\mathbf{U}},  {\mathbf{U}}\Big\rangle_{\mathbf{H}} >0,
\end{align*}
and the discrete energy norm
\begin{align*}
\|\mathbf{U}\|^2_{hP} :=  \Big\langle {\mathbf{U}},  P^{-1}{\mathbf{U}}\Big\rangle_{\mathbf{H}} > 0.
\end{align*}
We have 
 \begin{theorem}\label{theo:discrete_case_no_damping}
 Consider the  semi-discrete  approximation   \eqref{eq:disc_elastic_pml_1A}-\eqref{eq:disc_elastic_pml_3A}. When all PML absorption functions vanish,  $d_{\xi} = 0$, the solution of the semi-discrete approximation satisfies the energy identity
 {\small
 \begin{equation*}
\frac{d}{dt} \|\mathbf{U}\|^2_{hP} \le 0.
\end{equation*}
}
 \end{theorem}
 The proof of Theorem \ref{theo:discrete_case_no_damping} can be adapted from \cite{DuruGabrielIgel2017,Duru_exhype_2_2019}. We will not repeat it here.
 
 Note that Theorem \ref{theo:discrete_case_no_damping} is not applicable when the PML is present, $d_\xi \ne 0$ for any $\xi = x, y, z$.
\subsection{Discrete energy estimate for the PML in Laplace space}
Taking the Laplace transform, in time, of the semi-discrete  problem  \eqref{eq:disc_elastic_pml_1A}-\eqref{eq:disc_elastic_pml_3A}  gives
{
\begin{equation}\label{eq:disc_elastic_pml_1A_Laplace}
\begin{split}
\mathbf{P}^{-1} s{  \widetilde{\mathbf{U}}} = \sum_{\xi = x, y, z}\left[\mathbf{A}_{\xi}\mathbf{D}_{\xi} \widetilde{\mathbf{U}}  - \mathbf{d}_{\xi}\mathbf{w}_{\xi} - \mathbf{H}_{\xi}^{-1}{\left({\mathbf{e}_{\xi}(-1)} \widetilde{\mathbf{FL}}_{\xi}+  {\mathbf{e}_{\xi}(1)} \widetilde{\mathbf{FR}}_{\xi} \right)}\right] + \widetilde{\mathbf{F}}_{U} ,
  \end{split}
  \end{equation}
\begin{equation}\label{eq:disc_elastic_pml_3A_Laplace}
\begin{split}
 s\widetilde{\mathbf{w}}_{\xi}
 = \mathbf{A}_{\xi}\mathbf{D}_{\xi} \widetilde{\mathbf{U}}  - \left(\mathbf{d}_{\xi} + \boldsymbol{\alpha}_{\xi}\right)\widetilde{\mathbf{w}}_{\xi}- \underbrace{\theta_{\xi}\mathbf{H}_{\xi}^{-1}{\left({\mathbf{e}_{\xi}(-1)} \widetilde{\mathbf{FL}}_{\xi}+  {\mathbf{e}_{\xi}(1)} \widetilde{\mathbf{FR}}_{\xi}  \right)}}_{\text{PML stabilizing flux fluctuations}} +   \widetilde{\mathbf{F}}_{w_\xi}.
  \end{split}
  \end{equation}
}

Next, we use \eqref{eq:disc_elastic_pml_3A_Laplace} and  eliminate the auxiliary variable $\widetilde{\mathbf{w}}_{\xi}$ from  \eqref{eq:disc_elastic_pml_1A_Laplace}. We have
{
\begin{equation}\label{eq:Laplace_disc_elastic_pml_1A}
\begin{split}
\mathbf{P}^{-1} \left(sS_x\right) \widetilde{\mathbf{U}} = & \sum_{\xi = x, y, z} \frac{S_x}{S_\xi}\left[\frac{1}{2}\left(\mathbf{A}_{\xi}\mathbf{D}_{\xi}  - \mathbf{H}_{\xi}^{-1}\mathbf{A}_{\xi}\mathbf{D}_{\xi}^{T}\mathbf{H}_{\xi}    + \mathbf{H}_{\xi}^{-1}\mathbf{A}_{\xi}\left(\mathbf{B}_{\xi}\left(1,1\right)\right) -  \mathbf{B}_{\xi}\left(-1,-1\right)\right)\widetilde{\mathbf{U}} \right]\\
- &\sum_{\xi = x, y, z} \frac{S_x}{S_\xi}\left[\mathbf{H}_{\xi}^{-1}{\left({\mathbf{e}_{\xi}(-1)} \widetilde{\mathbf{FL}}_{\xi} +  {\mathbf{e}_{\xi}(1)} \widetilde{\mathbf{FR}}_{\xi}  \right)} \right] + \mathbf{P}^{-1} \widetilde{\mathbf{F}} \\
+& \underbrace{\sum_{\xi = x, y, z} \frac{\mathbf{d}_\xi S_x \left(1-\theta_{\xi}\right)}{S_\xi \left(s + \alpha_{\xi}\right)}\mathbf{H}_{\xi}^{-1}{\left({\mathbf{e}_{\xi}(-1)} \widetilde{\mathbf{FL}}_{\xi} +  {\mathbf{e}_{\xi}(1)} \widetilde{\mathbf{FR}}_{\xi}  \right)}}_{\text{destabilizing PML  flux term}}. 
  \end{split}
  \end{equation}
}
Note that in \eqref{eq:Laplace_disc_elastic_pml_1A}, we  have made use of  the discrete integration-by-parts property \eqref{eq:disc_sbp} of the spatial derivative operator $\mathbf{D}_{\xi}$.
The last term in the right hand side of \eqref{eq:Laplace_disc_elastic_pml_1A} is a destabilising PML  flux term at element faces. For a DG method in one element or a standard SBP finite difference approximation in one block the destabilising terms will appear at external boundaries. And for a multi-element/block DG/SBP discretisation the  destabilising PML  flux terms will also appear every internal element boundaries in the PML. In particular, the DG  is an element based method, these destabilizing flux terms appear almost everywhere in the PML, including at internal and external element faces. The destabilizing flux terms can be eliminated using the parameter $\theta_\xi $. We set the penalty parameter $\theta_\xi =1$, extending the numerical implementation of the boundary conditions and inter-element conditions to the auxiliary differential equations. 
 The last  term in  the right hand side of \eqref{eq:Laplace_disc_elastic_pml_1A} vanishes, having
 {
 \begin{equation}\label{eq:Laplace_disc_elastic_pml_1A_0}
\begin{split}
\mathbf{P}^{-1} \left(sS_x\right) \widetilde{\mathbf{U}} = & \sum_{\xi = x, y, z} \frac{S_x}{S_\xi}\left[\frac{1}{2}\left(\mathbf{A}_{\xi}\mathbf{D}_{\xi}  - \mathbf{H}_{\xi}^{-1}\mathbf{A}_{\xi}\mathbf{D}_{\xi}^{T}\mathbf{H}_{\xi}    + \mathbf{H}_{\xi}^{-1}\mathbf{A}_{\xi}\left(\mathbf{B}_{\xi}\left(1,1\right)\right) -  \mathbf{B}_{\xi}\left(-1,-1\right)\right)\widetilde{\mathbf{U}} \right]\\
- &\sum_{\xi = x, y, z} \frac{S_x}{S_\xi}\left[\mathbf{H}_{\xi}^{-1}{\left({\mathbf{e}_{\xi}(-1)} \widetilde{\mathbf{FL}}_{\xi} +  {\mathbf{e}_{\xi}(1)} \widetilde{\mathbf{FR}}_{\xi}  \right)} \right] + \mathbf{P}^{-1} \widetilde{\mathbf{F}}.
  \end{split}
  \end{equation}
  }
  We will now make the review specific for the acoustic wave equation and the elastic wave equations.
  \subsubsection{Discrete energy estimate of the PML for the acoustic wave equation}
As before we will  consider first the acoustic wave equation, where $\widetilde{\mathbf{U}} = (\mathbf{\widetilde{p}}, \widetilde{\mathbf{v}})^T$, and summarise the discrete stability results.
More elaborate discussions can be found in \cite{KDuru2016,DuruGabKreiss2019}.
For simplicity we consider a single element DG approximation, by restricting \eqref{eq:Laplace_disc_elastic_pml_1A_0} to one element, and divide through by $S_x$, we have 
{
\footnotesize
  \begin{equation}\label{eq:acoustic_pml_3D_Laplace_disc_BC1}
\begin{split}
\frac{1}{\kappa} {s \mathbf{\widetilde{p}}}  + \widetilde{{\grad}}_{D}  \cdot \widetilde{\mathbf{v}}= & \frac{1}{\kappa}\mathbf{\widehat{f}}_p-\sum_{\xi }\frac{1}{S_{\xi}}\mathbf{H}_{\xi}^{-1}{\left(\frac{\mathbf{B}_{\xi}(-1, -1)}{Z} \left(\frac{1-r_{\xi}}{2} Z \mathbf{\widetilde{v}}_{\xi} + \frac{1+r_{\xi}}{2}\mathbf{\widetilde{p}} \right) -  \frac{\mathbf{B}_{\xi}(1, 1)}{Z} \left(\frac{1-r_{\xi}}{2} Z \mathbf{\widetilde{v}}_{\xi} - \frac{1+r_{\xi}}{2}\mathbf{\widetilde{p}} \right)\right)},
  \end{split}
  \end{equation}
    \begin{equation}\label{eq:acoustic_pml_3D_Laplace_disc_BC2}
\begin{split}
\rho {s \widetilde{\mathbf{v}} } + \widetilde{{\grad}}_{D}  \mathbf{\widetilde{p}} &= \rho \mathbf{ \widehat{f}} -\sum_{\xi }\frac{ \mathbf{n}}{S_{\xi}}\otimes\mathbf{H}_{\xi}^{-1}{\left({\mathbf{B}_{\xi}(-1, -1)} \left(\frac{1-r_{\xi}}{2} Z \mathbf{\widetilde{v}}_{\xi} + \frac{1+r_{\xi}}{2}\mathbf{\widetilde{p}} \right) +  {\mathbf{B}_{\xi}(1, 1)} \left(\frac{1-r_{\xi}}{2} Z \mathbf{\widetilde{v}}_{\xi} - \frac{1+r_{\xi}}{2}\mathbf{\widetilde{p}} \right)\right)}.
%
  \end{split}
  \end{equation}
  }
  In \eqref{eq:acoustic_pml_3D_Laplace_disc_BC1}--\eqref{eq:acoustic_pml_3D_Laplace_disc_BC2},  only external boundaries are present.
 Introducing the modified  discrete operators
{
\footnotesize
 \begin{equation}\label{eq:elemental_disc_operators_laplace}
\begin{split}
  &\widetilde{\mathbf{D}}_{\xi} =  \left(\mathbf{D}_{\xi}  + 
 \frac{1+r_{\xi} }{2}\mathbf{H}_{\xi} ^{-1} \left(\mathbf{B}_{\xi} \left(-1,-1\right) - \mathbf{B}_{\xi} \left(1,1\right)\right)\right), \quad
 \widetilde{\mathbf{D}}_{0\xi}  = \left(\mathbf{D}_{\xi}  + \frac{1-r_{\xi}}{2}\mathbf{H}_{\xi}^{-1}\left(\mathbf{B}_{\xi}(-1, -1) - \mathbf{B}_{\xi}(1, 1)\right)\right),
 \\
 &\widetilde{\mathbf{H}}_{\xi}
 =  \mathbf{H}  \widetilde{\mathbf{I}}_{\xi},
 \quad 
   \widetilde{\mathbf{I}}_{\xi} = \left(\mathbf{I} + \frac{c\left(1-r_{\xi}\right)}{2sS_{\xi}} \mathbf{H}_{\xi}^{-1} \left(\mathbf{B}_{\xi}(-1, -1)  +  \mathbf{B}_{\xi}(1, 1) \right)\right)^{-1},
  \end{split}
  \end{equation}
 }
with $c = Z/\rho$, and  eliminating the velocity fields  having
{ 
\begin{equation*}
\begin{split}
 & s^*\mathbf{H} s {\boldsymbol{\kappa}}^{-1}s\widetilde{\boldsymbol{p}}
+ 
\sum_{\xi}
\left(\frac{1}{S_{\xi}}\widetilde{\mathbf{D}}_{\xi}\right)^{\dagger}  \left(\frac{(s^*S_{\xi}^*)}{\boldsymbol{\rho}S_{\xi}} \widetilde{\mathbf{H}}_{\xi}\right)  \left(\frac{1}{S_{\xi}}\widetilde{\mathbf{D}}_{\xi}\right) \widetilde{\boldsymbol{p}}\\
 &=
|s|^2\sum_{\xi}\frac{1+r_{\xi}}{2ZS_{\xi}} {\mathbf{H}}{\mathbf{H}_{\xi}}^{-1}\left(\mathbf{B}_{\xi}{(-1,-1)} + \mathbf{B}_{\xi}{(1,1)}\right)\widetilde{\boldsymbol{p}}
s^*\mathbf{H} {\boldsymbol{\kappa}}^{-1}s\widehat{\mathbf{F}}_p
  -
\sum_{\xi}   s^*{\mathbf{H}} \left(\frac{1}{S_{\xi}}\widetilde{\mathbf{D}}_{0\xi}\right) \widehat{\mathbf{\widetilde{f}}}_{\xi},
\end{split}
\end{equation*}
}
where $\widehat{\mathbf{\widetilde{f}}}_{\xi} = \widetilde{\mathbf{I}}_{\xi} \widehat{\mathbf{{f}}}_{\xi}$. 

Note that the discrete operators $\widetilde{\mathbf{D}}_{0\xi}$, $\widetilde{\mathbf{D}}_{\xi}$, $\widetilde{\mathbf{H}}_{\xi}$ are modified by the boundary conditions.
On GLL collocation nodes $\widetilde{\mathbf{H}}_{\xi} $ is diagonal. In particular if we consider a hard wall $r_{\eta} = 1$, the boundary terms will vanish, we recover $\widetilde{\mathbf{H}}_{\xi} = \mathbf{H} = \mathbf{H}^T > 0$.
%
Consider $\mathbf{H} = \mathbf{H}^T > 0$, introduce the discrete scalar product and the corresponding norm
\begin{align}
\Big\langle \widetilde{\mathbf{u}},  \widetilde{\mathbf{v}}\Big\rangle_{\mathbf{H}} = \widetilde{\mathbf{v}}^{\dagger} \mathbf{H} \widetilde{\mathbf{u}}, \quad \|\widetilde{\mathbf{v}}\|_{\mathbf{H}}^2 = \Big\langle \widetilde{\mathbf{v}},  \widetilde{\mathbf{v}}\Big\rangle_{\mathbf{H}}.
\end{align}
For the modified discrete operator $ \Re{\left(\frac{(s^*S_{\xi}^*)}{S_{\xi}}\widetilde{\mathbf{H}}_{\xi}\right)}$ the scalar product can be decomposed into
{
\small
 \begin{align}\label{eq:modified_integration_operator_decomposition}
\mathbf{v}^{\dagger}\Re{\left(\frac{(s^*S_{\xi}^*)}{S_{\xi}}\widetilde{\mathbf{H}}_{\xi}\right)} \mathbf{v}  = a\left\langle\mathbf{v}, \mathbf{v} \right\rangle_{\widehat{\mathbf{H}}_{\xi}} +  \epsilon_{\xi}(s, d_{\xi})\left\langle\mathbf{v}, \mathbf{v} \right\rangle_{\widehat{\mathbf{H}}_{\xi}} + \left(1-r_{\xi}\right)\mathbf{BT}_{\mathrm{num}}^{(\xi)},
\end{align}
}
with 
$
\left\langle\mathbf{v}, \mathbf{v} \right\rangle_{\widehat{\mathbf{H}}_{\xi}} > 0,
$
and
$
\mathbf{BT}_{\mathrm{num}}^{(\xi)} > 0.
$

Introduce the discrete 
 energy norms
\begin{equation*}
\begin{split}
&\mathcal{\widetilde{E}}^2_p(s) =  \Big\langle {s}  \widetilde{\mathbf{p}}, {s}  \widetilde{\mathbf{p}}\Big\rangle_{\mathbf{H}/\kappa}+ \sum_{\xi = x,y,z} \Big\langle \frac{1}{S_{\eta}}\widetilde{\mathbf{D}}_{\xi}\widetilde{\mathbf{p}},\frac{1}{S_{\xi}}\widetilde{\mathbf{D}}_{\xi}\widetilde{\mathbf{p}}\Big\rangle_{  \widehat{\mathbf{H}}_{\xi} /\rho} > 0,  
  \end{split}
  \end{equation*}
  \begin{equation*}
\begin{split}
&\mathcal{\widetilde{E}}_f^2(s) = \Big\langle{s}  \widehat{\mathbf{F}}_{p}, {s}  \widehat{\mathbf{F}}_{p}\Big\rangle_{  \mathbf{H}/\kappa }+\sum_{\xi = x,y,z} \Big\langle \frac{1}{S_{\xi}}\widetilde{\mathbf{D}}_{0\xi}\widehat{\mathbf{\widetilde{f}}}_{\xi} , \frac{1}{S_{\xi}}\widetilde{\mathbf{D}}_{0\xi}\widehat{\mathbf{\widetilde{f}}}_{\xi}\Big\rangle_{  \kappa \mathbf{H} }  > 0.
  \end{split}
  \end{equation*}
  Note that
  \begin{align*}
  \Re\left(\frac{1}{S_{\xi}}\right) =  \frac{ (a+\alpha_\xi)(a+\alpha_\xi+d_{\xi})+ b^2}{\gamma_\xi\left( \left(a+\alpha_\xi+d_{\xi}\right)^2+b^2\right)} > 0.
  \end{align*}
  If $\Re{s} = a > 0$   then  $ \Re\left(\frac{1}{S_{\xi}}\right) > 0$ for any $d_{\xi} \ge 0$.
  The following theorem was proven in \cite{DuruGabKreiss2019}
\begin{theorem}\label{Theo:Stability_PML_Laplace_one_element}
Consider the one element DGSEM approximation of the PML  in the Laplace space, \eqref{eq:elemental_disc_operators_laplace},  with  constant damping  $d_{\xi}  \ge 0$, $\alpha_\xi =\alpha \ge 0$, $\gamma_\xi = \gamma >0$ and $\Re{s}  = a >  0$. 
If  $\theta_{\xi} = 1$,  then we have
{
\small 
\begin{equation*}
\begin{split}
&(\gamma a)\mathcal{\widetilde{E}}_p^2(s) +  \sum_{\xi = x,y,z} \Big\langle \frac{1}{S_{\xi}}\widetilde{\mathbf{D}}_{\xi}\widetilde{\mathbf{p}},\frac{1}{S_{\xi}}\widetilde{\mathbf{D}}_{\xi}\widetilde{\mathbf{p}}\Big\rangle_{  \mathcal{\epsilon}_{\xi}\widehat{\mathbf{H}}_{\xi} /\rho}  +\sum_{\xi} \frac{1-r_{\xi}}{\rho}\mathbf{BT}_{\mathrm{num}}^{(\xi)}+ \mathrm{BT}(s)    \le   2\mathcal{\widetilde{E}}_p(s) \mathcal{\widetilde{E}}_f(s), \\
& \mathrm{BT}(s)   =  |s|^2\sum_{\xi} \Re{\left(\frac{1}{S_{\xi}}\right)} \frac{1+r_{\xi}}{2Z} \widetilde{\boldsymbol{p}}^{\dagger}\left[{\mathbf{H}}{\mathbf{H}_{\xi}^{-1}}\left(\mathbf{B}_{\xi}{(-1,-1)} + \mathbf{B}_{\xi}{(1,1)}\right)\right]\widetilde{\boldsymbol{p}} \ge 0.
  \end{split}
\end{equation*}
}
\end{theorem}
Theorem \ref{Theo:Stability_PML_Laplace_one_element} is completely analogous to the continuous counterpart Theorem \ref{Theo:Stability_PML_Laplace_Acoustics}.
The result can be extended to multiple DG elements, see \cite{DuruGabKreiss2019} for details.

We will now review the discrete energy estimates for general systems, such as the linear elastic wave equation, as presented in \cite{ElasticDG_PML2019}.
\subsubsection{Discrete energy estimate of the PML for general systems}
We introduce the discrete weighted inner product and norm
\begin{align*}
\Big\langle{\widetilde{\mathbf{U}}, \widetilde{\mathbf{F}}}\Big\rangle_{h{P}} = \Big\langle \widetilde{\mathbf{U}},  {\mathbf{P}}^{-1}  \widetilde{\mathbf{F}}\Big\rangle_{\mathbf{H}}, \quad \|\widetilde{\mathbf{U}}\|^2_{hP} = \Big\langle{\widetilde{\mathbf{U}}, \widetilde{\mathbf{U}}}\Big\rangle_{h{P}}.
 \end{align*}
  Following \cite{ElasticDG_PML2019}, we  consider specifically the elastic wave equation in general media, and for the two elements model we introduce the surface terms, including the external boundary terms 
 \begin{align*}
 \mathcal{BT}_s\left(\widehat{\widetilde{v}}_{\eta}^{-}, \widehat{\widetilde{T}}_{\eta}^{-}\right)  =  \frac{\Delta{y}}{2}  \frac{\Delta{z}}{2} \sum_{i = 1}^{P+1}\sum_{k = 1}^{P+1}\sum_{\eta = x,y,z}\left(\left(  \widehat{\widetilde{T}}_\eta^{-*} \widehat{\widetilde{v}}_\eta^{-}  \right)\Big|_{-1}\right)_{i k} h_{i} h_{k} \le 0,
 \end{align*}
 \begin{align*}
\mathcal{BT}_s\left(\widehat{\widetilde{v}}_{\eta}^{+}, \widehat{\widetilde{T}}_{\eta}^{+}\right)  = \frac{\Delta{y}}{2}  \frac{\Delta{z}}{2} \sum_{i = 1}^{P+1}\sum_{k = 1}^{P+1}\sum_{\eta = x,y,z}\left(\left( \widehat{\widetilde{T}}_\eta^{+*} \widehat{\widetilde{v}}_\eta^{+}  \right)\Big|_{ 1}\right)_{i k} h_{i} h_{k}  \le 0,
 \end{align*}
  the interface term
  \begin{align*}
  \mathcal{IT}_s\left(\widehat{\widetilde{v}}^{\pm}, \widehat{\widetilde{T}}^{\pm} \right) = - \frac{\Delta{y}}{2}  \frac{\Delta{z}}{2} \sum_{i = 1}^{P+1}\sum_{k = 1}^{P+1}\sum_{\eta = x,y,z}\left( {\widehat{\widetilde{T}}}_\eta^*  \lJump{{\widehat{\widetilde{v}}}_\eta \rJump} \right)_{i k} h_{i} h_{k} \equiv 0,
\end{align*}
   and the fluctuation term
    \begin{align*}
  \mathcal{F}_{luc}\left({\widetilde{G}},  Z\right) = - \frac{\Delta{y}}{2}  \frac{\Delta{z}}{2} \sum_{\eta = x, y, z}\sum_{i = 1}^{P+1}\sum_{k = 1}^{P+1}\left(\left(\frac{1}{Z_{\eta}}|{\widetilde{G}}_\eta |^2\right)\Big|_{-1} + \left(\frac{1}{Z_{\eta}}|{\widetilde{G}}_\eta |^2\right)\Big|_{ 1}  \right)_{i, k} h_{i} h_{k}  \le 0.
  \end{align*}
Here $v^*$ denotes the complex conjugate of $v$, and the surface integrals have been approximated by quadrature rules.

As in the continuous setting,  we will consider first the 1D  PML strip problem, and proceed later to the edge and the corner regions.
\paragraph{The discrete PML strip problem.}
 In \eqref{eq:Laplace_disc_elastic_pml_1A_0}, let  $d_x = d\ge 0$, $d_\xi \equiv 0$, $\mathbf{D}_\xi = \widetilde{\mathbf{FL}}_{\xi} = \widetilde{\mathbf{FR}}_{\xi}= 0$  for $\xi = y, z$. 
 We arrive at
 {
 \begin{equation}\label{eq:Laplace_disc_elastic_pml_1A_1D}
\begin{split}
\mathbf{P}^{-1} \left(sS_x\right) \widetilde{\mathbf{U}} = &  \left[\frac{1}{2}\left(\mathbf{A}_{x}\mathbf{D}_{x}  - \mathbf{H}_{x}^{-1}\mathbf{A}_{x}\mathbf{D}_{x}^{T}\mathbf{H}_{x}    + \mathbf{H}_{x}^{-1}\mathbf{A}_{x}\left(\mathbf{B}_{\xi}\left(1,1\right)\right) -  \mathbf{B}_{x}\left(-1,-1\right)\right)\widetilde{\mathbf{U}} \right]\\
- &\left[\mathbf{H}_{x}^{-1}{\left({\mathbf{e}_{x}(-1)} \widetilde{\mathbf{FL}}_{x} +  {\mathbf{e}_{x}(1)} \widetilde{\mathbf{FR}}_{x}  \right)} \right] + \mathbf{P}^{-1} \widetilde{\mathbf{F}}.
  \end{split}
  \end{equation}
  }
We can prove the discrete equivalence of Theorem \ref{theo:PML_1D_Strip_Laplace}.
 \begin{theorem}\label{theo:PML_1D_Strip_Laplace_Disc}
 Consider the  semi-discrete PML equation  in the Laplace space \eqref{eq:Laplace_disc_elastic_pml_1A_1D}, with $d_x = d\ge 0$, $\alpha_x = \alpha \ge 0$, $\gamma_x = \gamma >0$ and $\Re(s) \ge a> 0$.  We have 
 $
 \Re(s S_x)\ge \gamma a > 0
 $
 and
 {
 \footnotesize
 \begin{equation*}
 \begin{split}
& \|\sqrt{\Re(s S_x)}\widetilde{\mathbf{U}}^{-}\left(s\right)\|^2_{hP} + \|\sqrt{\Re(s S_x)}\widetilde{\mathbf{U}}^{+}\left(s\right)\|^2_{hP}  \le  \|\widetilde{\mathbf{U}}^{-}\left(s\right)\|_{hP} \|\widetilde{\mathbf{F}}^{-}\left( s\right)\|_{hP}  + \|\widetilde{\mathbf{U}}^{+}\left(s\right)\|_{hP} \|\widetilde{\mathbf{F}}^{+}\left(s\right)\|_{hP}+  \widetilde{\mathrm{BT}}_h, \\
& \widetilde{\mathrm{BT}}_h =  \mathcal{BT}_s\left(\widehat{\widetilde{v}}_{\eta}^{-}, \widehat{\widetilde{T}}_{\eta}^{-}\right) + \mathcal{BT}_s\left(\widehat{\widetilde{v}}_{\eta}^{+}, \widehat{\widetilde{T}}_{\eta}^{+}\right)  +  \mathcal{IT}_s\left(\widehat{\widetilde{v}}^{\pm}, \widehat{\widetilde{T}}^{\pm} \right) +  \mathcal{F}_{luc}\left({\widetilde{G}}^{-},  Z^{-}\right)  + \mathcal{F}_{luc}\left({\widetilde{G}}^{+},  Z^{+}\right) \le 0,
\end{split}
\end{equation*}
}
where  the surface terms $\widetilde{\mathrm{BT}}_h$ are negative semi-definite.
 \end{theorem}
\paragraph{The discrete PML  edge problem.}
Consider now   the discrete  PML  \eqref{eq:Laplace_disc_elastic_pml_1A_0} in the $xy$-edge region.
 In \eqref{eq:Laplace_disc_elastic_pml_1A_0}, let  $d_x = d_y = d\ge 0$, $d_z \equiv 0$, $\mathbf{D}_z = \widetilde{\mathbf{FL}}_{z} = 0$. 
 {
 \begin{equation}\label{eq:Laplace_disc_elastic_pml_1A_2D}
\begin{split}
\mathbf{P}^{-1} \left(sS_x\right) \widetilde{\mathbf{U}} = &\sum_{\xi = x, y} \left[\frac{1}{2}\left(\mathbf{A}_{\xi}\mathbf{D}_{\xi}  - \mathbf{H}_{\xi}^{-1}\mathbf{A}_{\xi}\mathbf{D}_{\xi}^{T}\mathbf{H}_{\xi}    + \mathbf{H}_{\xi}^{-1}\mathbf{A}_{\xi}\left(\mathbf{B}_{\xi}\left(1,1\right)\right) -  \mathbf{B}_{\xi}\left(-1,-1\right)\right)\widetilde{\mathbf{U}} \right]\\
- &\sum_{\xi = x, y} \left[\mathbf{H}_{\xi}^{-1}{\left({\mathbf{e}_{\xi}(-1)} \widetilde{\mathbf{FL}}_{\xi} +  {\mathbf{e}_{\xi}(1)} \widetilde{\mathbf{FR}}_{\xi}  \right)} \right] + \mathbf{P}^{-1} \widetilde{\mathbf{F}}
  \end{split}
  \end{equation}
  }

We formulate the discrete equivalence of Theorem \ref{theo:PML_2D_Edge_Laplace}.
 \begin{theorem}\label{theo:PML_2D_Edge_Laplace_Disc}
 Consider the  semi-discrete PML equation  in the Laplace space \eqref{eq:Laplace_disc_elastic_pml_1A_2D}, with $d_x = d_y = d \ge 0$, $\alpha_x = \alpha_y = \alpha \ge 0$, $\gamma_x = \gamma_y = \gamma >0$ and $\Re(s) \ge a> 0$.  We have 
 $
 \Re(s S_x)\ge \gamma a > 0
 $
 and
 {
 \footnotesize
 \begin{equation*}
 \begin{split}
& \|\sqrt{\Re(s S_x)}\widetilde{\mathbf{U}}^{-}\left(s\right)\|^2_{hP} + \|\sqrt{\Re(s S_x)}\widetilde{\mathbf{U}}^{+}\left(s\right)\|^2_{hP}  \le  \|\widetilde{\mathbf{U}}^{-}\left(s\right)\|_{hP} \|\widetilde{\mathbf{F}}^{-}\left( s\right)\|_{hP}  + \|\widetilde{\mathbf{U}}^{+}\left(s\right)\|_{hP} \|\widetilde{\mathbf{F}}^{+}\left(s\right)\|_{hP}+  \widetilde{\mathrm{BT}}_h,\\
& \widetilde{\mathrm{BT}}_h =  \mathcal{BT}_s\left(\widehat{\widetilde{v}}_{\eta}^{-}, \widehat{\widetilde{T}}_{\eta}^{-}\right) + \mathcal{BT}_s\left(\widehat{\widetilde{v}}_{\eta}^{+}, \widehat{\widetilde{T}}_{\eta}^{+}\right)  +  \mathcal{IT}_s\left(\widehat{\widetilde{v}}^{\pm}, \widehat{\widetilde{T}}^{\pm} \right) +  \mathcal{F}_{luc}\left({\widetilde{G}}^{-},  Z^{-}\right)  + \mathcal{F}_{luc}\left({\widetilde{G}}^{+},  Z^{+}\right) \le 0,
\end{split}
\end{equation*}
}
where the surface terms $\widetilde{\mathrm{BT}}_h$ are negative semi-definite.
%
 \end{theorem}
 \paragraph{The discrete PML  corner problem.}
Consider now  a the discrete  PML   \eqref{eq:Laplace_disc_elastic_pml_1A_0} in the corner region, where all damping functions are nonzero, $d_\xi  = d \ge 0$ for all $\xi = x, y,z$.
For this case the  PML complex metrics are identical $S_y = S_x = S_z$.  Thus $S_x/S_\xi = 1$, and from \eqref{eq:Laplace_disc_elastic_pml_1A_0}  we have
 {
 \begin{equation}\label{eq:Laplace_disc_elastic_pml_1A_3D}
\begin{split}
\mathbf{P}^{-1} \left(sS_x\right) \widetilde{\mathbf{U}} = &\sum_{\xi = x, y, z} \left[\frac{1}{2}\left(\mathbf{A}_{\xi}\mathbf{D}_{\xi}  - \mathbf{H}_{\xi}^{-1}\mathbf{A}_{\xi}\mathbf{D}_{\xi}^{T}\mathbf{H}_{\xi}    + \mathbf{H}_{\xi}^{-1}\mathbf{A}_{\xi}\left(\mathbf{B}_{\xi}\left(1,1\right)\right) -  \mathbf{B}_{\xi}\left(-1,-1\right)\right)\widetilde{\mathbf{U}} \right]\\
- &\sum_{\xi = x, y, z} \left[\mathbf{H}_{\xi}^{-1}{\left({\mathbf{e}_{\xi}(-1)} \widetilde{\mathbf{FL}}_{\xi} +  {\mathbf{e}_{\xi}(1)} \widetilde{\mathbf{FR}}_{\xi}  \right)} \right] + \mathbf{P}^{-1} \widetilde{\mathbf{F}}.
  \end{split}
  \end{equation}
  }

We will now state the discrete equivalence of  Theorem \ref{theo:PML_3D_Corner_Laplace}.
 \begin{theorem}\label{theo:PML_3D_Corner_Laplace_Disc}
 Consider the  semi-discrete PML equation  in the Laplace space \eqref{eq:Laplace_disc_elastic_pml_1A_3D}, with $d_\xi  = d\ge 0$,  $\alpha_\xi  = \alpha  \ge 0$, $\gamma_\xi =\gamma >0$ and $\Re(s) \ge  a > 0$.  We have 
 $
 \Re(s S_x)\ge \gamma a > 0
 $
 and
 {
 \footnotesize
 \begin{equation*}
 \begin{split}
& \|\sqrt{\Re(s S_x)}\widetilde{\mathbf{U}}^{-}\left(s\right)\|^2_{hP} + \|\sqrt{\Re(s S_x)}\widetilde{\mathbf{U}}^{+}\left(s\right)\|^2_{hP}  \le  \|\widetilde{\mathbf{U}}^{-}\left(s\right)\|_{hP} \|\widetilde{\mathbf{F}}^{-}\left( s\right)\|_{hP}  + \|\widetilde{\mathbf{U}}^{+}\left(s\right)\|_{hP} \|\widetilde{\mathbf{F}}^{+}\left(s\right)\|_{hP}+  \widetilde{\mathrm{BT}}_h,  \\
& \widetilde{\mathrm{BT}}_h =  \mathcal{BT}_s\left(\widehat{\widetilde{v}}_{\eta}^{-}, \widehat{\widetilde{T}}_{\eta}^{-}\right) + \mathcal{BT}_s\left(\widehat{\widetilde{v}}_{\eta}^{+}, \widehat{\widetilde{T}}_{\eta}^{+}\right)  +  \mathcal{IT}_s\left(\widehat{\widetilde{v}}^{\pm}, \widehat{\widetilde{T}}^{\pm} \right) +  \mathcal{F}_{luc}\left({\widetilde{G}}^{-},  Z^{-}\right)  + \mathcal{F}_{luc}\left({\widetilde{G}}^{+},  Z^{+}\right) \le 0,
\end{split}
\end{equation*}
}
where the surface terms are negative semi-definite.
 \end{theorem}
A discrete equivalence of Theorem \ref{Theo:estimate_PML} also holds: 
\begin{corollary}\label{corollary:estimate_PML_Disc}
Consider the semi-discrete energy estimate in the Laplace space
\begin{equation*}
\begin{split}
& \|\widetilde{\mathbf{U}} (\cdot ,\cdot ,\cdot , s)\|_{hP}^2 \le \frac{1}{\gamma a }\|\widetilde{\mathbf{U}} (\cdot ,\cdot ,\cdot , s)\|_{hP} \|\widetilde{\mathbf{F}} (\cdot ,\cdot ,\cdot , s)\|_{hP}, \quad \Re(s)\ge a>0.
\end{split}
\end{equation*}
For any $a> 0$, $\gamma >0$ and $T>0$ we have
  \begin{align*}
  \int_0^{T} e^{-2at}  \|{\mathbf{U}} (\cdot ,\cdot ,\cdot , t)\|_{hP}^2 dt \le  \frac{1}{\left(\gamma a\right)^2}\int_0^{T} e^{-2at} \left( \| {\mathbf{F}} (\cdot ,\cdot ,\cdot , t)\|_{hP}^2\right) dt.
  \end{align*}
\end{corollary}
We can use Theorem \ref{corollary:estimate_PML_Disc} to invert the the estimates in Theorems \ref{Theo:Stability_PML_Laplace_one_element}, \ref{theo:PML_1D_Strip_Laplace_Disc}, \ref{theo:PML_2D_Edge_Laplace_Disc} and \ref{theo:PML_3D_Corner_Laplace_Disc} to get a discrete energy estimate in the physical space.
As before the estimate seemingly allows the energy of the solution to grow  exponentially with time for general data. However, by choosing $a>0$ in relation to the length  of the time interval of interest, a bound  involving only algebraic growth in time follows.

The following remarks are of important note.

\begin{remark}
For general first order systems the continuous energy estimates, Theorems \ref{theo:PML_1D_Strip_Laplace}--\ref{theo:PML_3D_Corner_Laplace}, and the discrete energy estimates, Theorems \ref{theo:PML_1D_Strip_Laplace_Disc}--\ref{theo:PML_3D_Corner_Laplace_Disc}, are valid for all media parameters including those that violate the so-called {\it geometric stability condition}. These results may not eliminate instabilities that exist in the continuous problems but they ensure robustness and guarantee longtime numerical stability for the discrete PML, in particular for physical models satisfying the geometric stability conditions.
\end{remark}

\begin{remark}
We also remark that the discrete energy estimates, Theorems \ref{theo:PML_1D_Strip_Laplace_Disc}--\ref{theo:PML_3D_Corner_Laplace_Disc}, are valid for all discrete derivative operators that satisfy  the SBP property \eqref{eq:disc_sbp}. These include, but not limited to, DG spectral difference operators defined on Gaussian-type quadrature nodes and standard SBP finite difference operator on equidistant grids.
\end{remark}
\begin{remark}
For the general second order systems, \eqref{eq:linear_2nd_wave_PML_1A}--\eqref{eq:linear_2nd_wave_PML_3A} with the boundary condition \eqref{eq:bc_second_order_operator} and the interface conditions 
\eqref{eq:elastic_acoustic_interface_second_PML}, such as displacement formulation of linear elastodynamics, there is no  general result yet to guide the derivation of provably stable numerical methods for the PML.
The development of such a general procedure will be useful in designing reliable and efficient numerical methods for second order systems, such as multi-block SBP finite difference methods \cite{DURU201437,doi:10.1137/130947210,doi:10.1007/s10915-014-9817-1,MATTSSON20088753}, and energy-based and symmetric interior penalty DG methods \cite{doi:10.1137/05063194X,doi:10.1137/140973517,APPELO2018362,Antonietti2016,10.1111/j.1365-246X.2008.03915.x}.
\end{remark}
\section{Numerical problems, What can go wrong? How it can be fixed?}\label{sec:s6}
We will perform some numerical experiments to review and demonstrate what can go wrong with a numerical PML and show how it can be fixed.
To do this we consider the computational setup, Figure \ref{fig:Waveguide}, a semi-infinite waveguide 
 {
\begin{figure} [h!]
 \centering \includegraphics[width=0.6\textwidth]{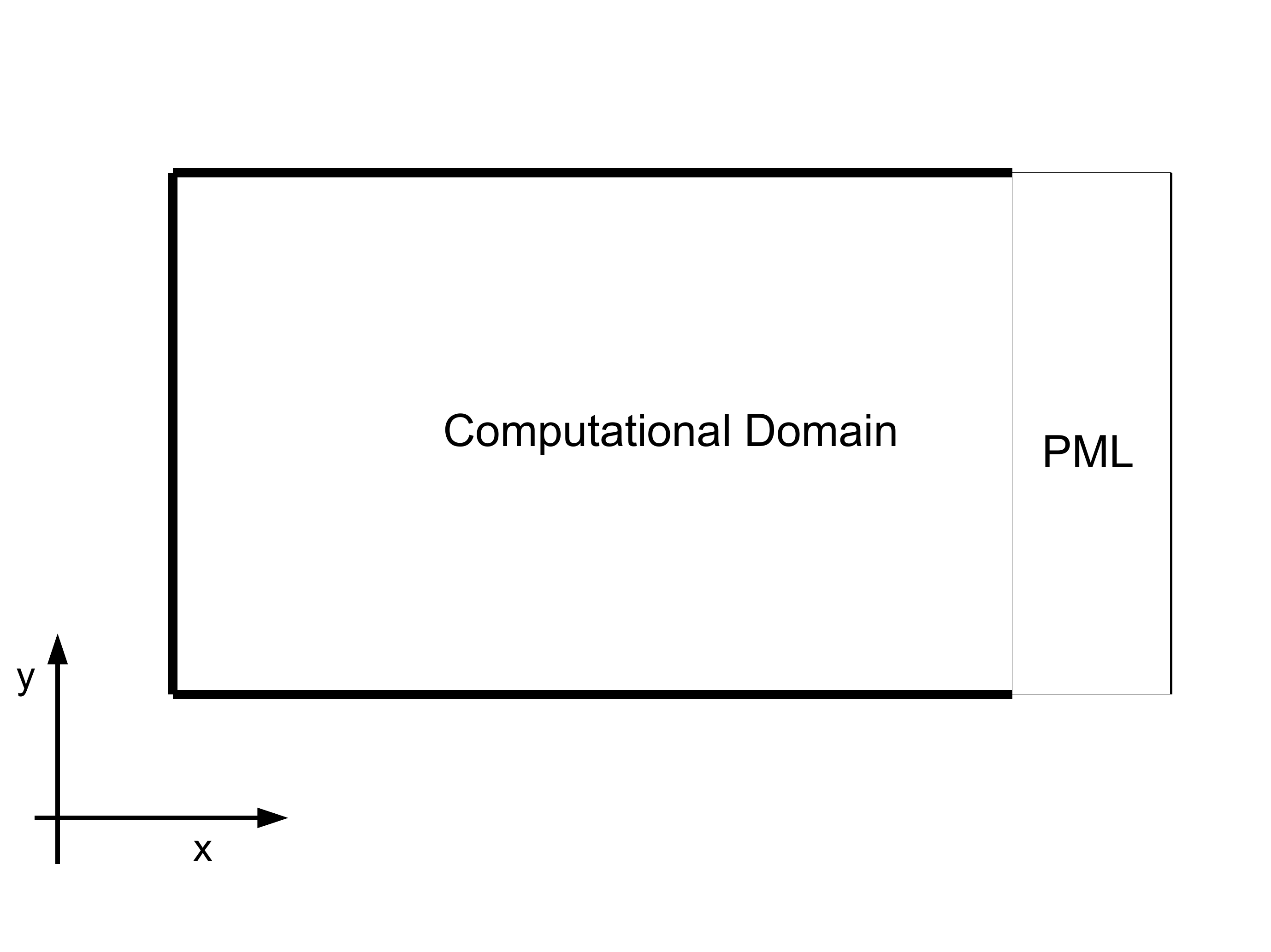}
\caption{\textit{A semi-infinite waveguide truncated by the PML of finite width $\delta_x > 0$.}}
 \label{fig:Waveguide}
\end{figure}
}
We truncate the waveguide by the PML of finite width $\delta_x > 0$, having $\left(x,y\right) \in \Omega = [-L_x, L_x+\delta_x]\times[0, L_y]$.
 The domain $\Omega$ can be an electromagnetic medium, an acoustic medium or an elastic medium, including isotropic and anisotropic elastic properties, inside the waveguide. 
 Throughout we will consider  medium parameters for which the {\it geometric stability condition} is satisfied and the unbounded constant coefficient PML is stable.
 
 We will consider DG discretisations but the results extend to SBP finite difference methods.
 The  discrete PMLs for different wave media  exhibit very similar numerical  instabilities.
 As will be shown  these numerical difficulties can be cured by the theory developed in \cite{DuruKozdonKreiss2016,ElasticDG_PML2019,KDuru2016,DuruGabKreiss2019}, and summarised in previous sections of this paper.

 We consider a 2D problem with $\partial /\partial z = 0$, $d_z(z) = 0$, the velocity component $v_z$ and the auxiliary variables $\mathbf{w}_z$ drop out. For linear elasticity the stress field components $\sigma_{xz}, \sigma_{yz}, \sigma_{zz}$ also drop out.
 We will also consider the  vertical PML  strip problem (with $d_x(x) \ge 0$, $d_y(y) = 0$), that is a PML in the $x$-direction truncating the right boundary. 
 In all experimental setups, we will review the  effectiveness and robustness  of the PML stabilising parameter  $\theta_x=1$, as presented \cite{DuruKozdonKreiss2016,ElasticDG_PML2019,KDuru2016,DuruGabKreiss2019}.
 
 The damping profile is a cubic monomial
\begin{equation}\label{eq:damping_func}
\begin{split}
&d_x\left(x\right) = \left \{
\begin{array}{rl}
0 \quad {}  \quad {}& \text{if} \quad x \le L_x,\\
d_0\Big(\frac{x-L_x}{\delta_x}\Big)^3  & \text{if}  \quad x \ge L_x ,
\end{array} \right.
\end{split}
\end{equation}
where $d_0> 0$ is the damping strength.
We set the CFS $\alpha = 0.15$, the damping strength
\begin{equation}\label{eq:damping_coef}
d_0 = \frac{4c_p}{2\delta_x}\ln{\frac{1}{{tol}}},
\end{equation}
where $c_p$ denotes the P-wave speed, and  ${tol}=10^{-3}$ is the magnitude of the relative PML error \cite{SjPe}. 
We note that evanescent waves can be problematic for standard PML models \cite{805891,465042,535836}, unless the width of  the layer is expanded. Since the derivation of the formula \eqref{eq:damping_coef} is based on damping of propagating waves, and it must be considered as a heuristic formula when evanescent waves are present. However, it has been demonstrated in \cite{Roden_and_Gedney_2000}  that the inclusion of the complex frequency shift, CFS, in the PML metric will enhance the absorption of evanescent waves. Also numerical simulations presented in \cite{DuruKozdonKreiss2016,ElasticDG_PML2019,DuruGabKreiss2019,DuruFungWilliams2020,ExaHyPE2019} suggest that the formula \eqref{eq:damping_coef} is a good way to choose the PML damping coefficient for practical simulations.

We will evolve the Gaussian 
initial wave profile
\begin{align}\label{eq:init_condition_pressure}
f(x, y) = e^{-\log\left({2}\right)\frac{x^2 + (y-Ly/2)^2}{9}},
\end{align}
At the top and bottom boundaries, $y = 0, L_y$,  we set a free-surface boundary condition with $r_y = 1$ in \eqref{eq:boundary_condition_acoustics}--\eqref{eq:boundary_condition_elastic}, the left boundary at $x = -L_x$ is a soft/clamped  wall with $r_x = -1$ in \eqref{eq:boundary_condition_acoustics}--\eqref{eq:boundary_condition_elastic} and the PML boundary at $x = L_x+\delta_x$ is terminated with the classical first order ABC, by setting  $r_x = 0$ in \eqref{eq:boundary_condition_acoustics}--\eqref{eq:boundary_condition_elastic}. 
When the PML damping vanishes, the set-up will correspond to a semi-infinite waveguide terminated by the classical first order ABC, and will reflect all waves with non-normal incident angle. The first order ABC is exact for waves impinging at normal incidence on the boundary. Specifically, for linear elasticity  this is the first order Lysmer-Kuhlemeyer ABC which is exact for plane shear and pressure waves with normal incidence angles at the boundary \cite{LysmerKuhlemeyer1969,Baffe2012}. We will compare the accuracy of the PML and the ABC.

We set the domain width $L_\xi = 50$ km, the PML width $\delta_x = 10$ km.
We discretise the domain with a uniform element size $\Delta{x} = \Delta{y} = 5$ km, spanning the PML with two DG elements, and approximate the solution by a polynomial of degree $N = 4$, with the sub-cell resolution $h= \Delta{x}/(N+1) = 1$~km. 
Time-integration is performed using the high order  ADER  scheme \cite{Toro1999,DumbserKaeser2006}  of the same order of accuracy with the spatial discretization.
We use the time step
\begin{align}
\Delta{t} = \frac{CFL}{ \sqrt{d} \left(2N+1\right) c_p} \min{\left(\Delta{x}, \Delta{y}\right)},
\end{align}
with the $CFL = 0.9$ number and the P-wave speed $c_p$ [km/s] and $d=2$ is the spatial dimension. The final time is $t = 200 $ s.

In order to access numerical PML errors we compute a reference solution in a larger domain $\left(x,y\right) \in \Omega_{ref} = [-L_x, 3 L_x]\times[0, L_y]$ and for a limited time such that reflection from the reference domain do not re-enter the interior computational domain. By comparing the reference solution and the PML solution in the interior, $\left(x,y\right) \in  [-L_x,  L_x]\times[0, L_y]$, we get an accurate measure of the PML error.


 \subsection{Acoustic waveguide}
\begin{figure}[h!]
\begin{subfigure}{\textwidth}
    \centering
 \stackunder[5pt]{\includegraphics[width=0.305\textwidth]{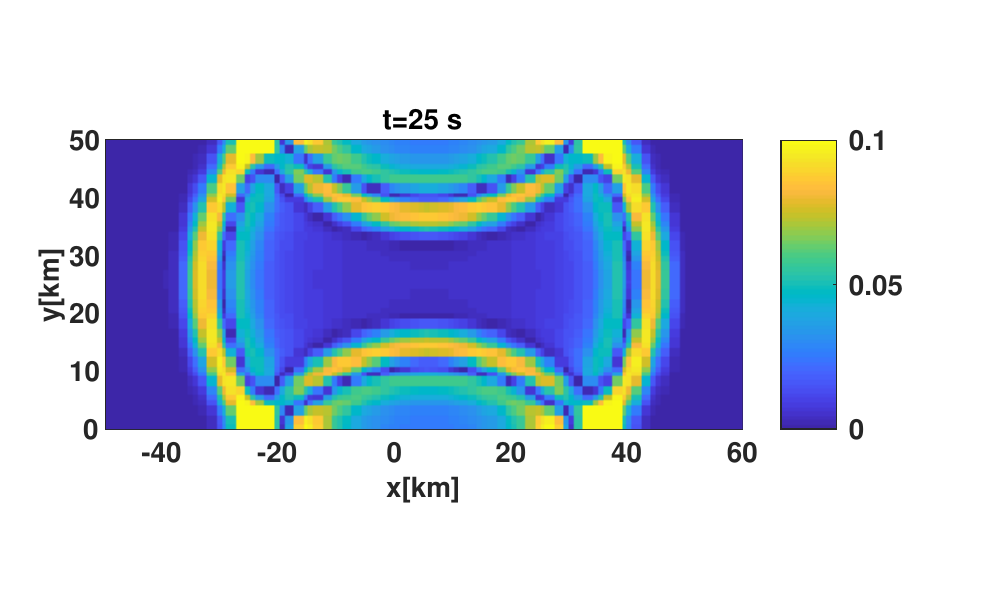}}{}%
\hspace{-1.25cm}%
\stackunder[5pt]{\includegraphics[width=0.305\textwidth]{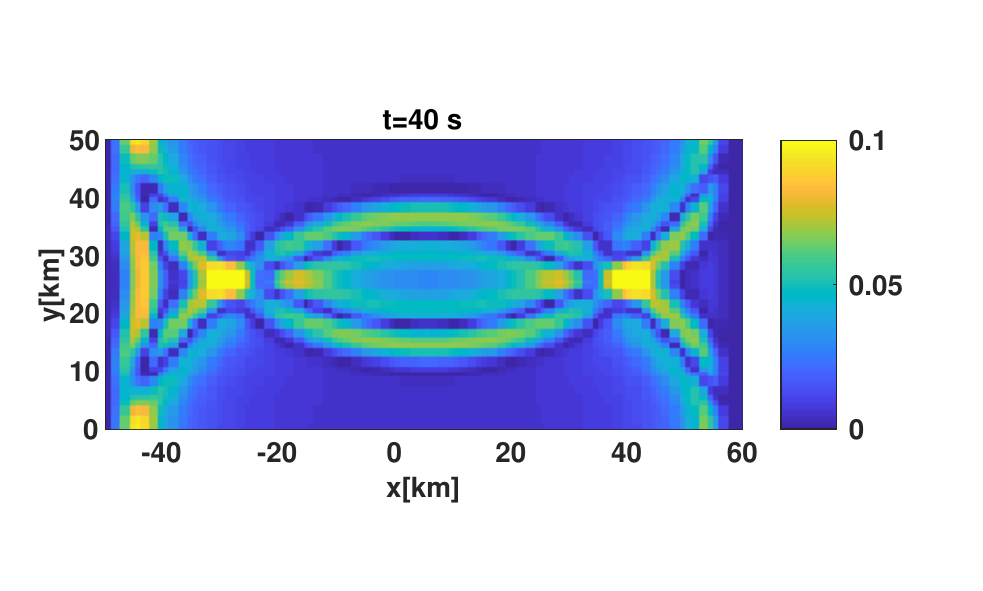}}{}%
\hspace{-1.25cm}%
\stackunder[5pt]{\includegraphics[width=0.305\textwidth]{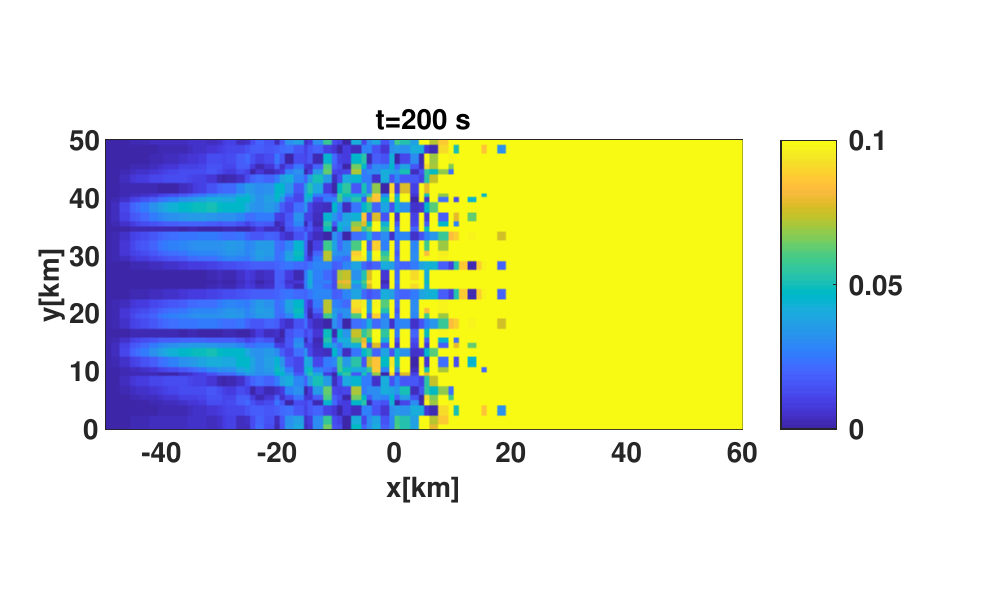}}{$\theta_x = 0$.}%
\hspace{-1.25cm}%
\stackunder[5pt]{\includegraphics[width=0.305\textwidth]{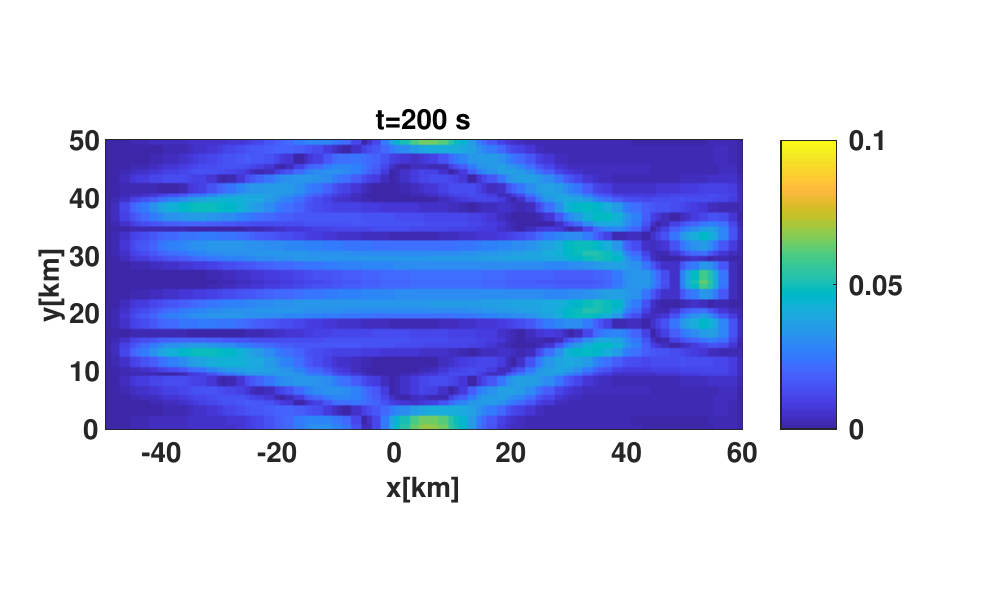}}{$\theta_x = 1$.}%
     \end{subfigure}
    \caption{Absolute pressure field $|p(x,y,t)|$ at $t = 25, 40, 200$ s for $\theta_x = 0, 1$. The penalty parameter $\theta_x = 1$ stabilises the PML solution at long times, $t = 200$ s.}
    \label{fig:pressure_t200s}
\end{figure}

\begin{figure} [h!]
 \centering
\begin{subfigure}{0.45\textwidth}
    \centering
 \includegraphics[width=\linewidth]{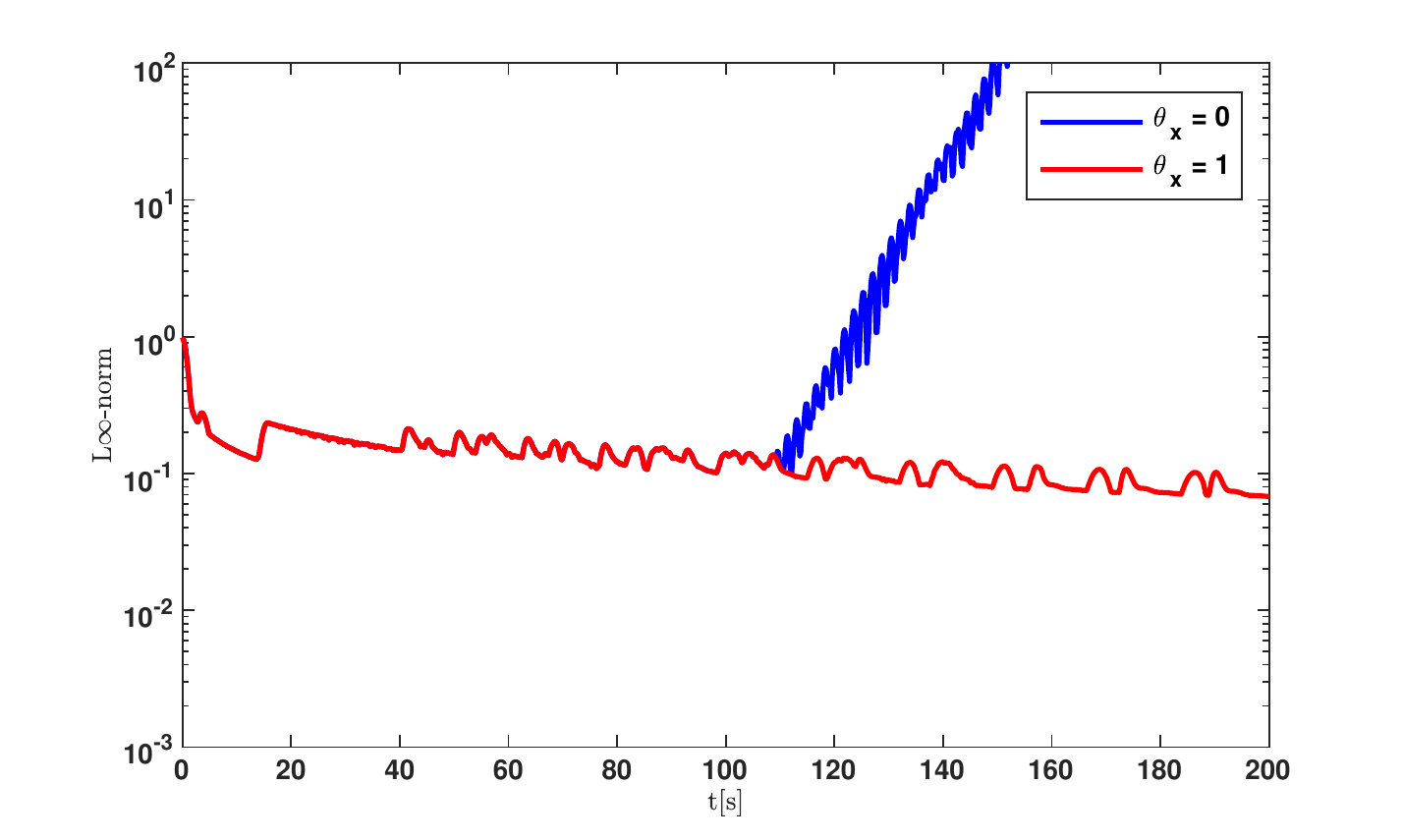}
 \caption{L$_{\infty}$ norm of $p$ for the PML with $\theta_x = 0, 1$}
  \label{fig:Time_series_pressure_norm}
  \end{subfigure}
  \begin{subfigure}{0.45\textwidth}
    \centering
 \includegraphics[width=\linewidth]{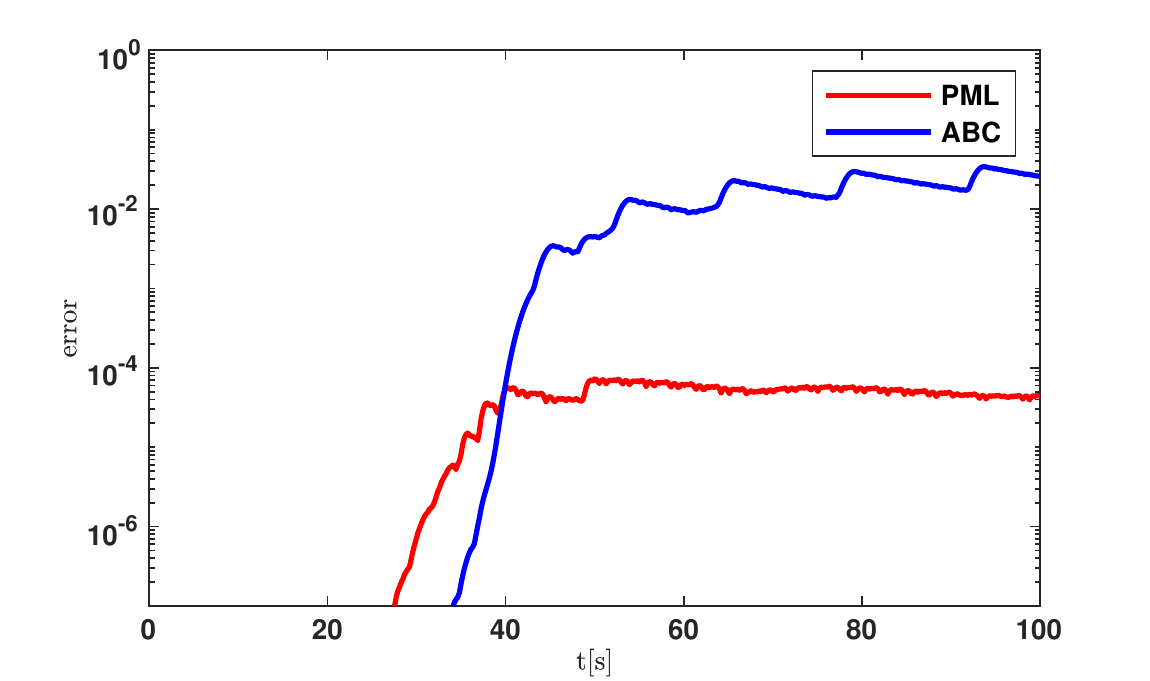}
 \caption{PML error versus ABC error for $p$}
  \label{fig:Time_series_error_pml_vs_abc}
  \end{subfigure}
 \caption{Time series of the L$_{\infty}$ norm of the pressure field for the PML with $\theta_x = 0, 1$.}
 \label{fig:Time_series_pressure}
\end{figure}

 We consider first a 2D acoustic medium with  constant acoustic wave speed $c_p = 1.484$ km/s   and constant medium density $\rho = 1$ g/cm$^3$.
We set the initial condition \eqref{eq:init_condition_pressure}, $p(x, y, t=0) = f(x,y)$, for the pressure field,
 and zero initial condition for the velocity fields and the auxiliary variables.
  The snapshots of the absolute pressure fields are plotted in Figure \ref{fig:pressure_t200s}, for $\theta_x = 0, 1$, and the time history  L$_{\infty}$-norm of the pressure field plotted in Figure  \ref{fig:Time_series_pressure_norm}.

We can see from Figure \ref{fig:pressure_t200s} that  without the stabilisation parameter, that is with $\theta_x =0$, the solution in the PML explodes after some time steps. See also Figure \ref{fig:Time_series_pressure_norm}.  As shown in \cite{DuruGabKreiss2019}, the initiation time of the explosive numerical mode depends on the quadrature rule used and the mesh resolution.  On a finer mesh the growth persists, but it starts at a much later time. 
The PML stabilisation parameter $\theta_x =1$ guarantees both accuracy and numerical stability of the PML at long times. See also Figure \ref{fig:Time_series_pressure_norm} and \ref{fig:Time_series_error_pml_vs_abc}. While the ABC is numerically stable, it however generates unacceptable large spurious reflections errors. The errors can be reduced by moving the artificial boundary further away. This however will lead to a significant increase in the computational cost.

These experiments demonstrate the significant importance of proving continuous and discrete energy estimates for the PML in the Laplace space, and its impact in ensuring long-time stability of the PML. The results extend to 2D Maxwell's equations of electrodynamics and to SBP-SAT finite difference method \cite{KDuru2016}.
The initial ideas of proving continuous and discrete energy estimates for the PML in the Laplace space were first derived for the 2D Maxwell's equation using the SBP-SAT finite difference method \cite{KDuru2016}. In \cite{DuruGabKreiss2019}, these were recently extended to  acoustic waves, and to 3D and DG methods. Although the proofs are derived for constant and piece-wise damping coefficients, however, numerical experiments demonstrate that the analysis may be valid for continuously varying damping profiles.
 \subsection{Isotropic and anisotropic elastic waveguides}
 Next we  consider the PML truncating a linear elastic waveguide. Both isotropic and anisotropic media properties will be considered.
The density of the media and elastic constants are given in Table \ref{tbl:ElasticParameters}.
In both isotropic and anisotropic elastic media we have chosen the media properties such that the maximum P-wave speed is $c_p = 6$ km/s. The dispersion relation of the isotropic and anisotropic elastic media considered here are shown in Figure \ref{fig:dispersion_relation}(b)--2(c), and indicate that the Cauchy PML problem is stable.
 
\begin{table}[h!]
  \caption{Elastic media parameters.}
  \centering{
    \begin{tabular}{|c|c|c|}
      \hline
      Parameters & Isotropic medium & Anisotropic medium (AM1)
        \\ \hline
      $c_{11}$ [GPa] & 97.20     & 20    
      \\ \hline
      $c_{12}$ [GPa]& 36.85 & 3.8
      \\ \hline
       $c_{22}$ [GPa]& 97.20 & 4 
       \\ \hline
        $c_{33}$ [GPa]& 30.17 & 2
         \\ \hline
        $\rho $ [gm/cm$^3$]& $2.7$ & $20/36$
        \\ \hline
    \end{tabular}
  }
  \label{tbl:ElasticParameters}
\end{table}

In elastic media, we set the initial condition \eqref{eq:init_condition_pressure}, $v_x(x, y, t=0) =v_y(x, y, t=0)=f(x,y)$, for the velocity field,
and zero initial conditions for the stress fields and the auxiliary variables.

\begin{figure}[h!]
 %
\begin{subfigure}{\textwidth}
    \centering
 \stackunder[5pt]{\includegraphics[width=0.305\textwidth]{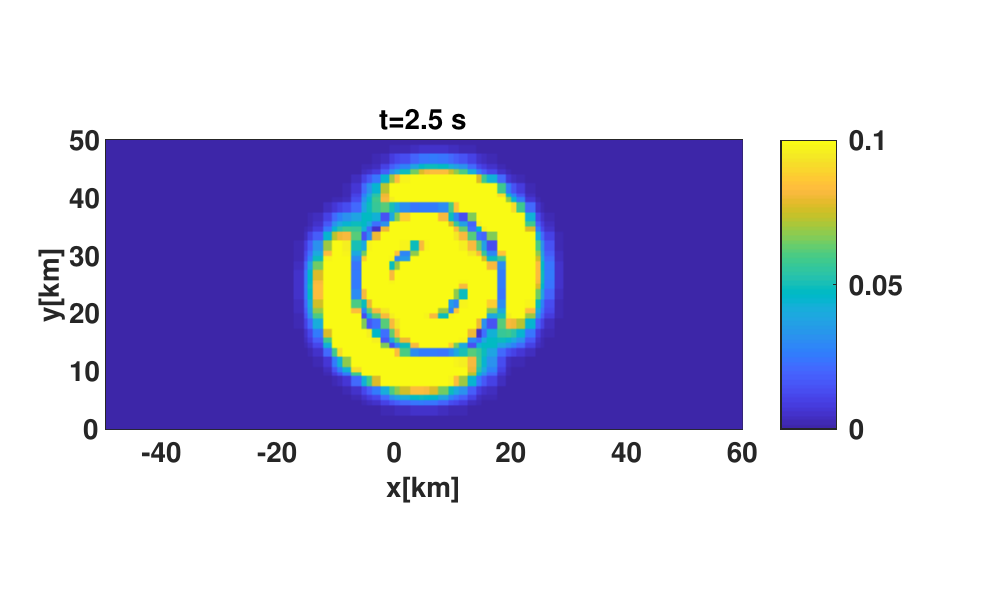}}{}%
\hspace{-1.25cm}%
\stackunder[5pt]{\includegraphics[width=0.305\textwidth]{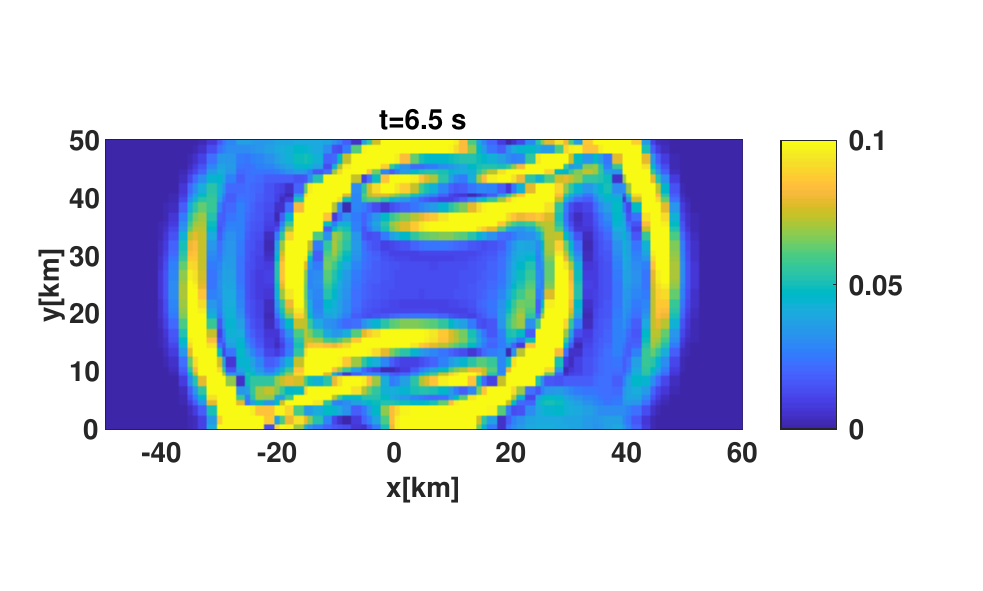}}{}%
\hspace{-1.25cm}%
\stackunder[5pt]{\includegraphics[width=0.305\textwidth]{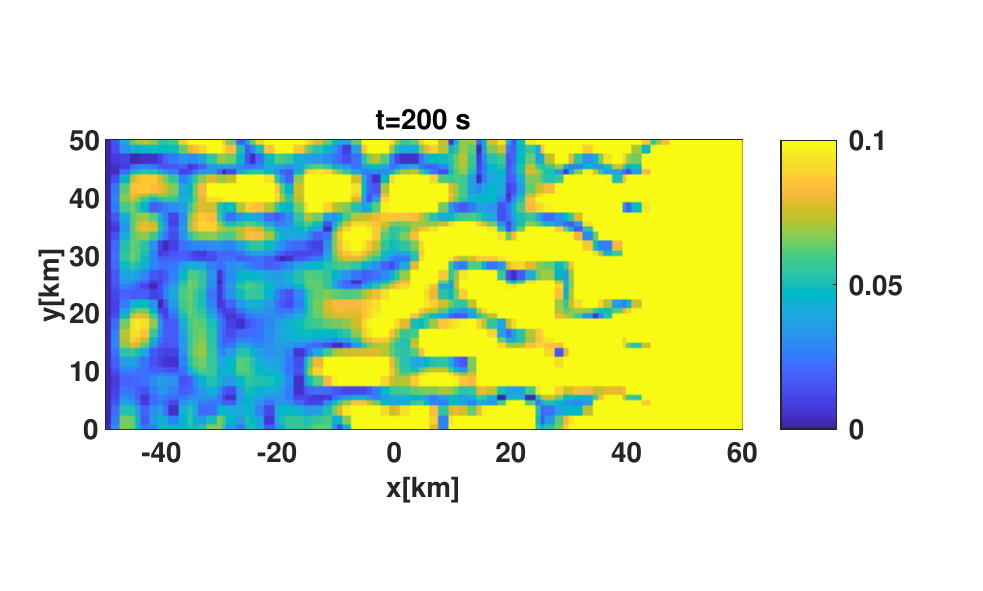}}{$\theta_x = 0$.}%
\hspace{-1.25cm}%
\stackunder[5pt]{\includegraphics[width=0.305\textwidth]{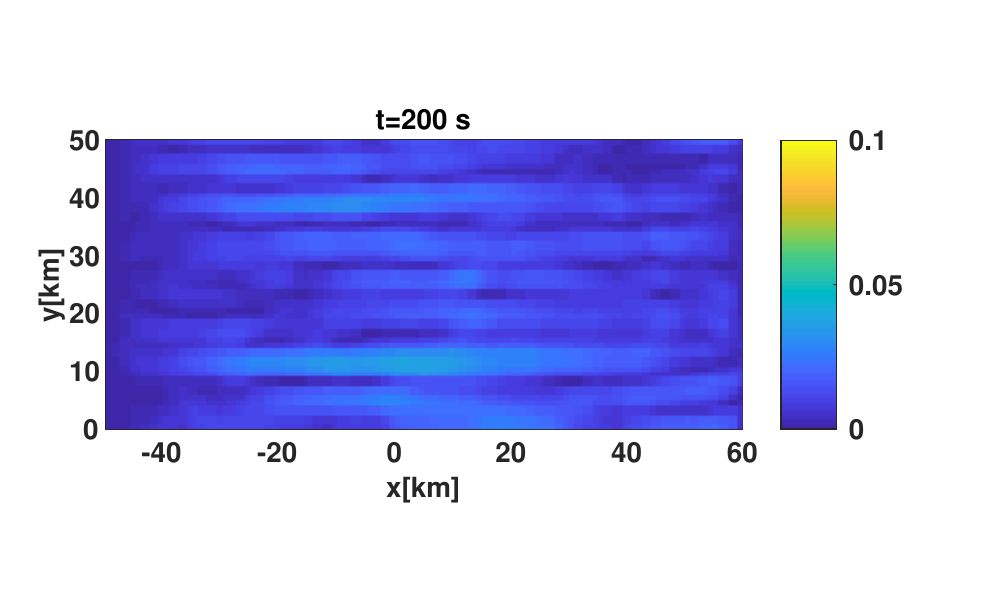}}{$\theta_x = 1$.}%
 \label{fig:isotropic_elastic_snapshot}
 \caption{Isotropic medium}
  \end{subfigure}
\begin{subfigure}{\textwidth}
    \centering
\stackunder[5pt]{\includegraphics[width=0.305\textwidth]{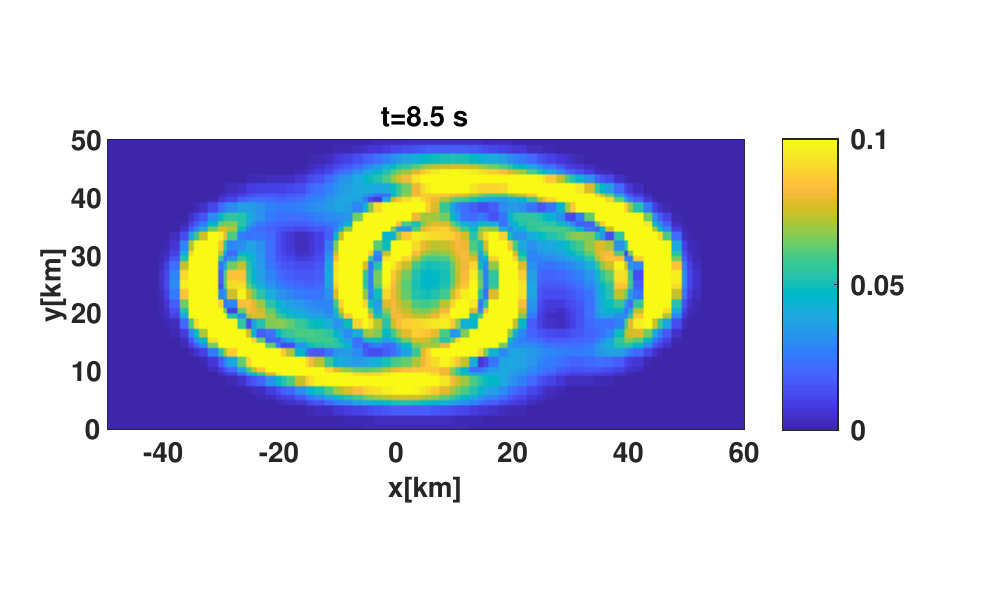}}{}%
\hspace{-1.25cm}%
 \stackunder[5pt]{\includegraphics[width=0.305\textwidth]{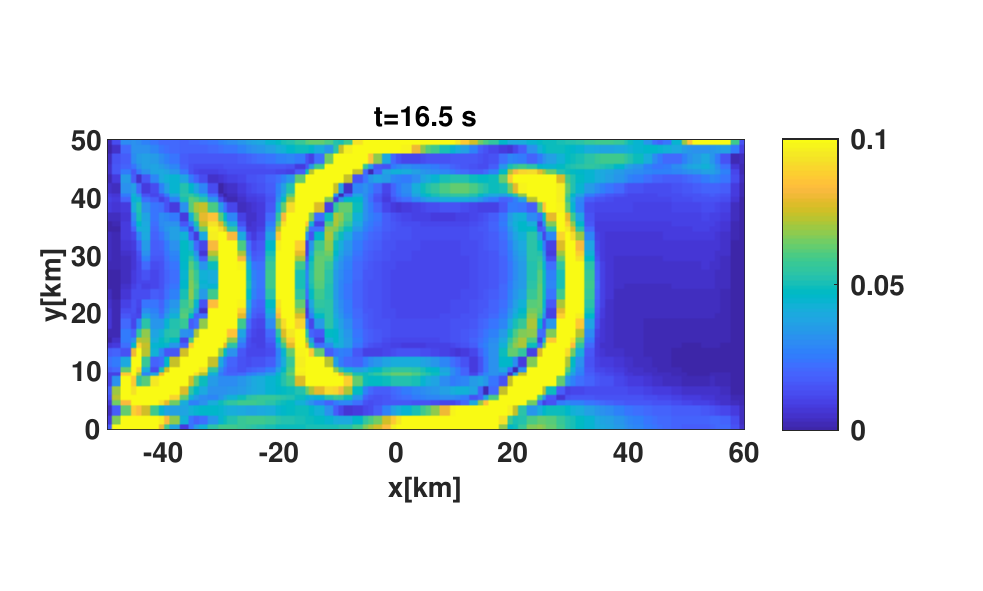}}{}%
\hspace{-1.25cm}%
\stackunder[5pt]{\includegraphics[width=0.305\textwidth]{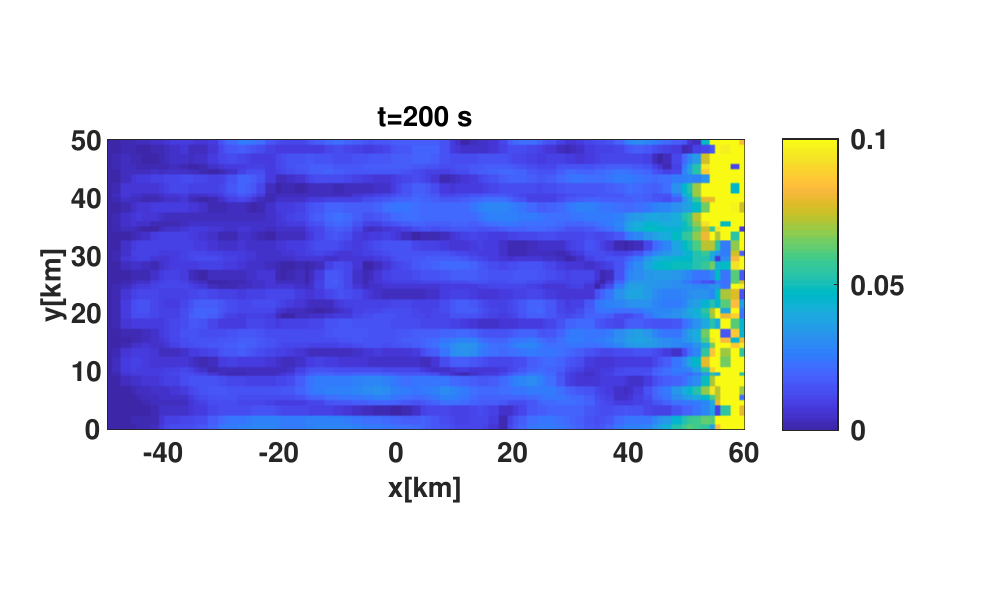}}{$\theta_x = 0$.}%
\hspace{-1.25cm}%
\stackunder[5pt]{\includegraphics[width=0.305\textwidth]{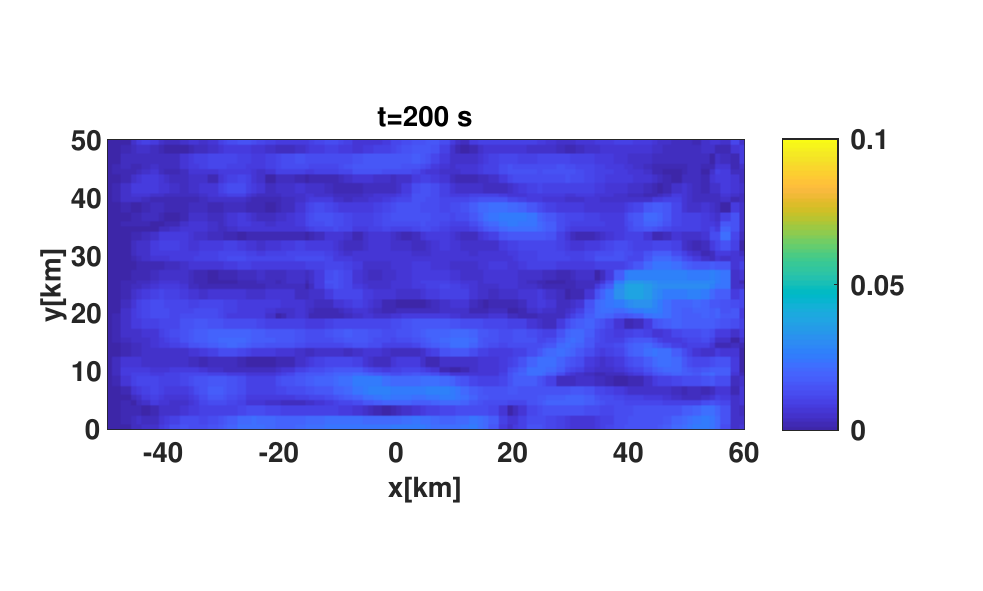}}{$\theta_x = 1$.}%
\label{fig:anisotropic_elastic_snapshot}
 \caption{Anisotropic medium}
  \end{subfigure}
%
    \caption{Snapshots of the absolute particle velocity $\sqrt{v_x^2 + v_y^2}$ in 2D elastic solids.}
     \label{fig:elastic_velocity_snapshot}
\end{figure}

\begin{figure}[h!]
	\centering
	\begin{subfigure}{.475\textwidth}
		\includegraphics[width=\textwidth]{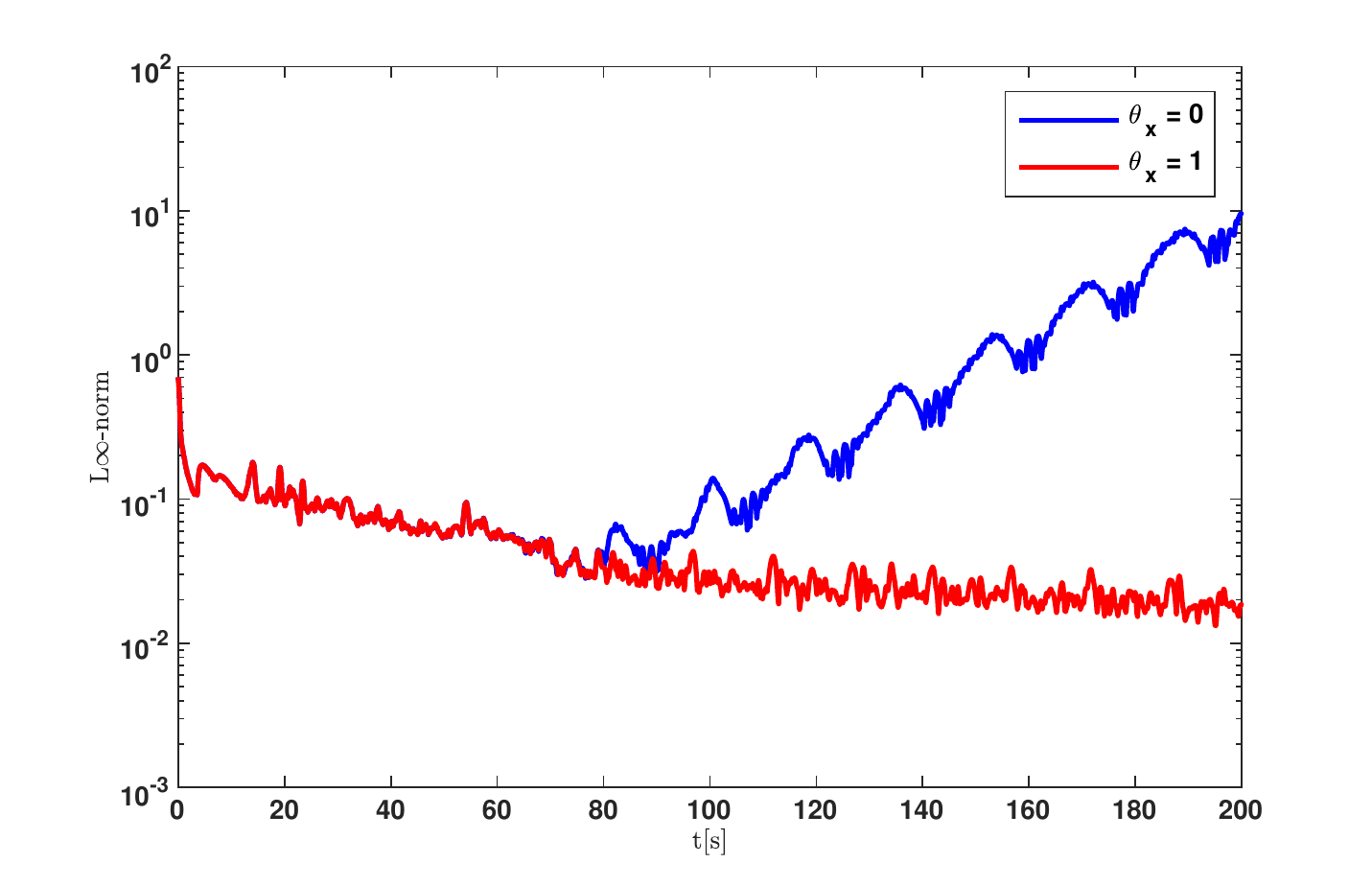}
		\caption{Isotropic medium: L$_{\infty}$ norm of $\mathbf{v}$ for the PML.}
	\end{subfigure}
	\begin{subfigure}{.515\textwidth}
		\includegraphics[width=\textwidth]{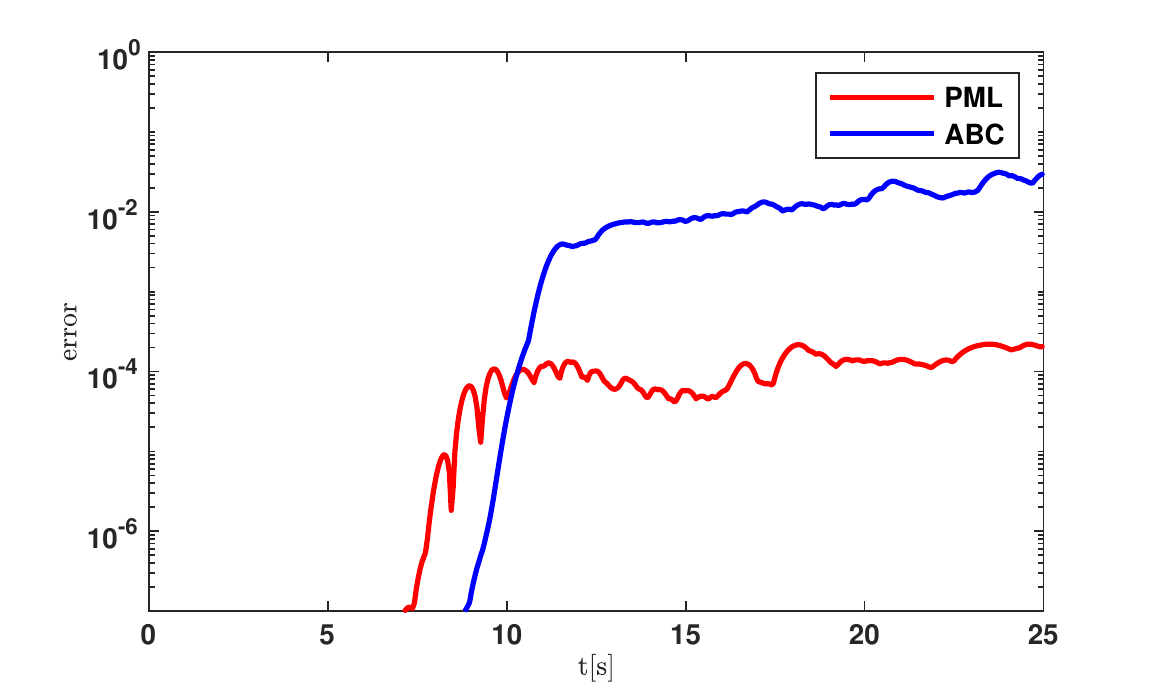}
		\caption{Isotropic medium: PML error versus ABC error.}
	\end{subfigure}
	\begin{subfigure}{.45\textwidth}
		\includegraphics[width=\textwidth]{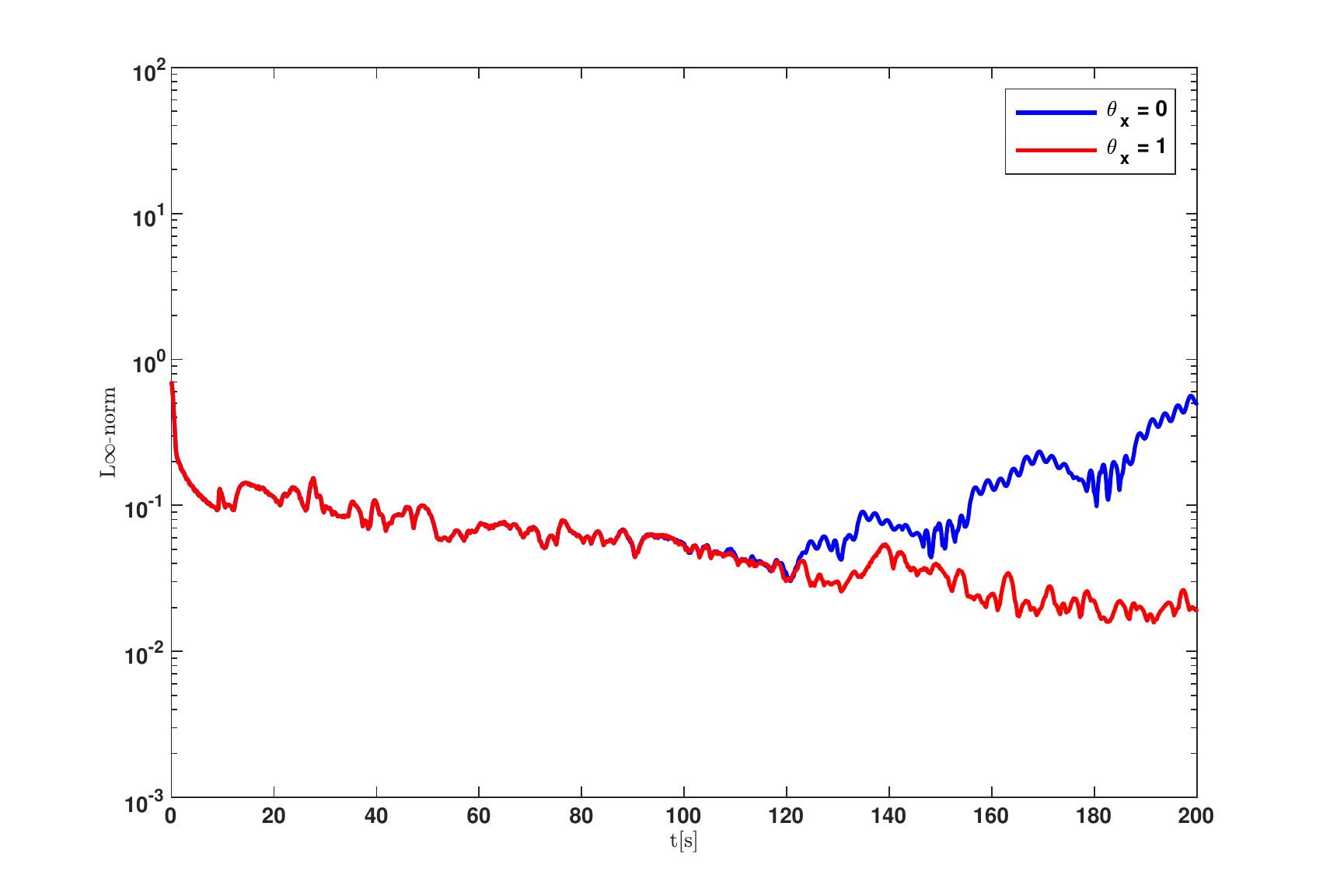}
		\caption{Anisotropic medium: L$_{\infty}$ norm of $\mathbf{v}$ for the PML.}
	\end{subfigure}
         \begin{subfigure}{.535\textwidth}
		\includegraphics[width=\textwidth]{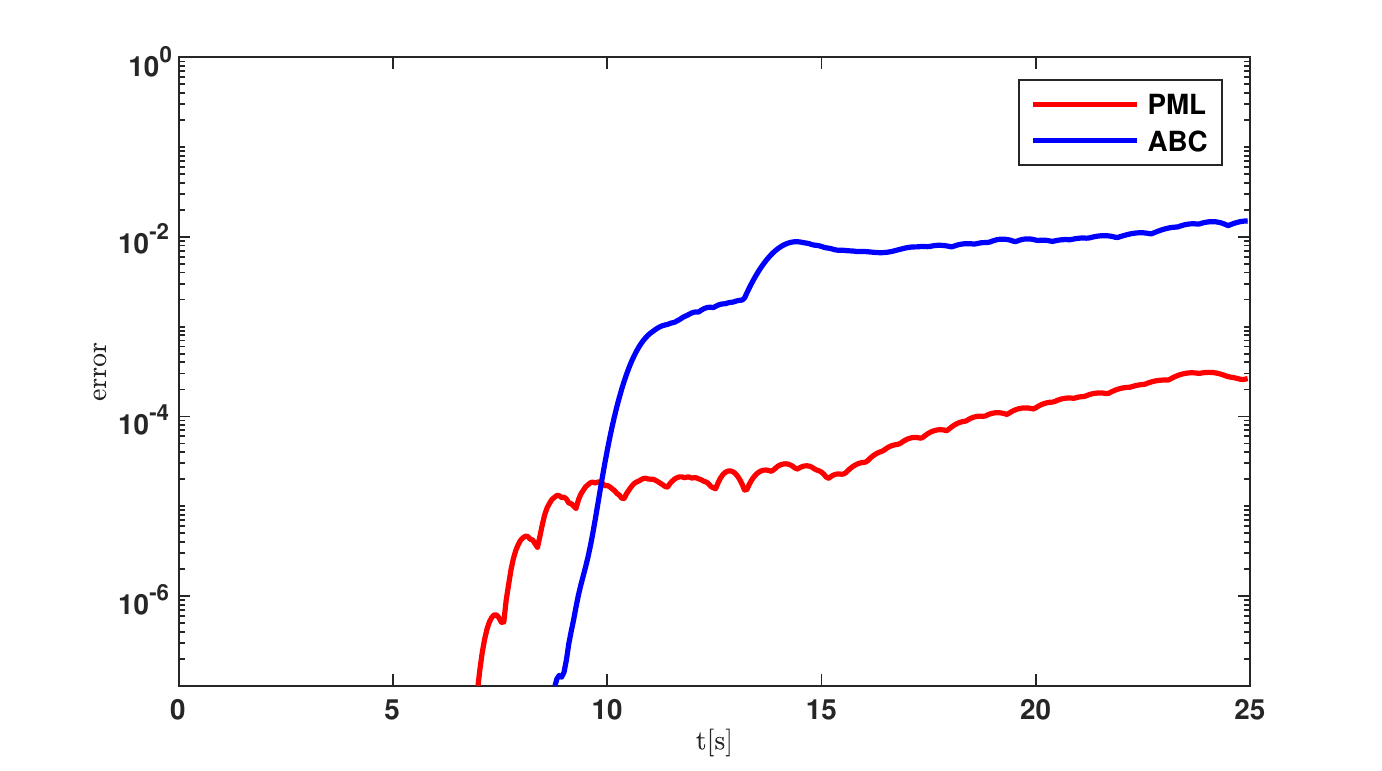}
		\caption{Anisotropic medium: PML error versus ABC error.}
	\end{subfigure}
	\caption{Time series of the error and L$_{\infty}$ norm of the velocity vector for the PML.}
	\label{fig:elastic_velocity_time_series}
\end{figure}

The snapshots of the absolute velocity are plotted in Figure \ref{fig:elastic_velocity_snapshot}, for $\theta_x = 0, 1$, and the time history  L$_{\infty}$-norm of the absolute velocity and the numerical PML errors  are plotted in Figure  \ref{fig:elastic_velocity_time_series}. 
We can also see from Figure \ref{fig:elastic_velocity_snapshot} and Figure \ref{fig:elastic_velocity_time_series} that  without the PML stabilisation parameter, that is with $\theta_x =0$, the solution in the PML explodes after some time steps.  As above, and shown in \cite{DuruKozdonKreiss2016,ElasticDG_PML2019}, the initiation time of the explosive numerical mode depends on the quadrature rule used and the mesh resolution.  On a finer mesh the growth persists, but it starts at a much later time. The consistent behaviour of the unstable mode indicate that the observed instabilities are numerical artefacts, and our stabilisation procedure provides a definite cure without any additional computational cost.

Wave propagation in an elastic media is much more complicated than in an acoustic or electromagnetic medium. For example in elastic media multiple body waves are supported including the existence of surface waves such as Rayleigh and Love waves. In elastic waveguides, with free-surface boundary conditions, the situation is much more complex as the presence of Rayleigh--Lamb wave modes can support backward propagating modes \cite{SkAdCr} which can be problematic for the PML. In \cite{DuKrSIAM,Duruthesis2012} we prove that the free-surface boundary conditions for the PML is stable and demonstrate in \cite{Duruthesis2012,DURU2014445} that the existence of backward propagating Rayleigh--Lamb wave modes in waveguides will not be problematic for a finite width PML for the 2D elastic waveguide. 

Unlike the acoustic wave equation, the numerical analysis of the PML for the elastic wave equation is technically challenging. As above the initial ideas of proving continuous and discrete energy estimates for the simplified PML model problems in the Laplace space were first derived for the 2D elastic wave equation using the SBP-SAT finite difference method \cite{DuruKozdonKreiss2016}. In \cite{ElasticDG_PML2019}, we generalise the technique to 3D and DG methods. As before the proofs are derived for model problems with  constant and piece-wise damping coefficients, however, numerical experiments demonstrate that the analysis may be valid for general 3D problem with continuously varying damping profiles.

\section{Applications}\label{sec:s7}
We will now consider some seismological application problems and demonstrate impact. We will consider both a standard 3D seismological benchmark problem, the LOH1 \cite{Seismowine, Kristekova_etal2009, Kristekova_etal2006} which has analytical exact solutions, and a real-world wave propagation propagation problem which involves the simulation of 3D seismic waves in a section of European Alpine region, with strong non-planar free-surface topography. Our algorithms and the PML have been implemented in two different freely open source  HPC software packages, WaveQLab3D \cite{DuruFungWilliams2020} and ExaHyPE \cite{ElasticDG_PML2019}, for large-scale simulation of seismic waves in geometrically complex 3D Earth models. The software package WaveQLab3D \cite{DuruFungWilliams2020} is  a  high order accurate SBP-SAT finite difference solver. ExaHyPE is a DG solver of arbitrary accuracy for large-scale numerical simulation of hyperbolic wave propagation problems on dynamically adaptive curvilinear meshes.
We will consider homogeneous initial conditions on the solution $\mathbf{U}(x,y,z,0) =0$ and  seismic sources defined by the singular double-couple moment tensor point source 
{
\begin{align}\label{eq:pointsource}
\mathbf{f}(x,y,z,t) &=   \mathbf{M} \delta_x(x-x_0)\delta_y(y-y_0)\delta_z(z-z_0)g(t), \quad 
\mathbf{M} = \left( 0, 0, 0 , 0,0,0,0,0,M_0\right)^T.
\end{align}
}
with the moment magnitude $M_0 = 10^{18} $ Nm.
Here, $ \delta_{\eta}(\eta)$ are the one dimensional Dirac delta function, $(x_0, y_0, z_0)$ is the source location, and the source time function is given by
$$
  g(t) = \frac{t}{T^2} \exp(-t/T) , \quad T = 0.1 ~\ s.
$$

\subsection{Layer over  homogeneous half-space (LOH1) 3D benchmark problem}
We consider the 3D seismological benchmark problem,  Layer Over Homogeneous Half-space (LOH1) \cite{Seismowine, Kristekova_etal2009, Kristekova_etal2006} benchmark problem, proposed by the European network: Seismic wave Propagation and Imaging in Complex media (SPICE) validation code. The LOH1 benchmark problem has an exact analytical solution which is often used to verify the accuracy of simulation codes for seismological applications.  We will review the results presented in  \cite{ElasticDG_PML2019,DuruFungWilliams2020}, for the  DG method and SBP finite difference method, where we verify the accuracy of the PML implementations for ExaHyPE and WaveQLab3D. 

The LOH1 benchmark setup has a planar free surface and an internal planar interface between a thin low velocity (soft) upper-layer and high velocity (hard) lower crust,  see  Figure \ref{fig:loh1_setup}. The material properties for the soft upper-layer and hard lower-half-space are 
\begin{align*}
&\text{soft upper crust}: \quad \rho = 2.6~\ \text{g}/\text{cm}^3, \quad c_p = 4~\ \text{km/s}, \quad c_s = 2~\ \text{km/s}, \quad 0\le  x \le 1~\ \text{km},\\
&\text{hard lower crust}: \quad \rho = 2.7~\ \text{g}/\text{cm}^3, \quad c_p = 6~\ \text{km/s}, \quad c_s = 3.464~\ \text{km/s}, \quad x  >1~\ \text{km}.
\end{align*}
Note that the medium parameters $\rho$, $c_p$ and $c_s$ are discontinuous across $x = 1$ km.

In the $y$- and $z$-direction, the domain of the problem is unbounded. 
In the positive $x$-direction (in-towards the Earth), the domain is also unbounded with the Earth's surface $x = 0$ having the free surface boundary condition, with the reflection coefficient $r_x = 1$ in \eqref{eq:boundary_condition_elastic}. 
The SPICE code validation project \cite{Seismowine} suggested to use large enough computational model, namely $\Omega_{L} =[0, 34]\times [-26,32]^2$, so as the seismograms in the receivers are not contaminated by spurious reflections from the artificial boundaries of the model.
This would correspond to the computational domain of volume  114376~km$^3$.

We will use the PML \cite{DuruKozdonKreiss2016,ElasticDG_PML2019} to absorb outgoing waves and prevent artificial reflections from the bounded computational domain. The PML allows us to sufficiently limit the domain to be the computational cube $\Omega = [0, 16.333~\text{km}]\times[-2.287~\text{km}, 14.046~\text{km}]\times[-2.287~\text{km}, 14.046~\text{km}]$  with only two DG elements around the computational boundaries where the PML is active. Please see also  Figure \ref{fig:loh1_setup}.  The computational domain $\Omega$ is only $4357.1$~km$^3$ in volume, and amounts to $ 3.8095 \%$ of the suggested large domain $\Omega_{L}$, thus saving as much as $96.19 \%  $ of  the required computational resources. Although the PML involves auxiliary variables and equations to be stored and evolved, however, the extra computational cost for evolving the auxiliary variables is very insignificant since they are only active inside the thin PML absorbing layer.

\newcommand{\ytk}{5.8}
\newcommand{\xtk}{5.8}
\newcommand{\ztk}{3.4}

\newcommand{\xrl}{1.7}
\newcommand{\yrl}{2.1}
\newcommand{\xrr}{3.7}
\newcommand{\yrr}{4.1}

\newcommand{ \lpml}{0.196}

\begin{figure}[H]
    \centering

\begin{tikzpicture}

\coordinate (O) at (0,0,0);
\coordinate (A) at (0,\ztk,0);
\coordinate (B) at (0,\ztk,\xtk);
\coordinate (C) at (0,0,\xtk);
\coordinate (D) at (\ytk,0,0);
\coordinate (E) at (\ytk,\ztk,0);
\coordinate (F) at (\ytk,\ztk,\xtk);
\coordinate (G) at (\ytk,0,\xtk);

\coordinate (A2) at (0,3.3,0);
\coordinate (B2) at (0,3.3,\xtk);
\coordinate (E2) at (\ytk,3.3,0);
\coordinate (F2) at (\ytk,3.3,\xtk);

\coordinate (Ar) at (\yrl,\ztk,\xrl);
\coordinate (Br) at (\yrl,\ztk,\xrr);
\coordinate (Er) at (\yrr,\ztk,\xrl);
\coordinate (Fr) at (\yrr,\ztk,\xrr);

\coordinate (Arb) at (\yrl,1.9954,\xrl);
\coordinate (Brb) at (\yrl,1.9954,\xrr);
\coordinate (Erb) at (\yrr,1.9954,\xrl);
\coordinate (Frb) at (\yrr,1.9954,\xrr);

\coordinate (Arl) at (\yrl,3.3,\xrl);
\coordinate (Brl) at (\yrl,3.3,\xrr);
\coordinate (Erl) at (\yrr,3.3,\xrl);
\coordinate (Frl) at (\yrr,3.3,\xrr);

 \draw[blue,fill=blue!10,opacity=0.4] (O) -- (C) -- (G) -- (D) -- cycle;
 \draw[blue,fill=blue!10,opacity=0.4] (O) -- (A) -- (E) -- (D) -- cycle;
 \draw[blue,fill=blue!10,opacity=0.4] (O) -- (A) -- (B) -- (C) -- cycle;
\draw[blue] (D) -- (E) -- (F) -- (G) -- cycle;
 \draw[blue] (C) -- (B) -- (F) -- (G) -- cycle;
 \draw[blue] (A) -- (B) -- (F) -- (E) -- cycle;

\draw[blue,fill=yellow!40,opacity=0.5] (A2) -- (B2) -- (F2) -- (E2) -- cycle;
\draw[blue,fill=yellow!40,opacity=0.7] (Arl) -- (Brl) -- (Frl) -- (Erl) -- cycle;

\draw[blue,fill=red!40,opacity=0.7] (Ar) -- (Br) -- (Fr) -- (Er) -- cycle;
\draw[blue,fill=red!40,opacity=0.7] (Arb) -- (Brb) -- (Frb) -- (Erb) -- cycle;

\draw[blue] (Ar) -- (Arb) ; 
\draw[blue] (Br) -- (Brb) ; 
\draw[blue] (Fr) -- (Frb) ; 
\draw[blue] (Er) -- (Erb) ; 

\node[circle, draw] at (2.6,3.2,3.2) () {H};

\node[draw] at (2.6+0.7348,3.4,3.2- 0.7348) () {6};
\node[draw] at (2.6+0.5764,3.4,3.2-0.8647) () {9};


\draw[black, -stealth] (2.6,3.4,3.2) -- (2.6,3.4,0) ; 
\draw[black, -stealth] (2.6,3.4,3.2) -- (5.8,3.4,3.2) ; 
\draw[black, -stealth] (2.6,3.4,3.2) -- (5.8,3.4,0) ;


\begin{scope}[shift = {(8,1,0)}, scale = 1.9954]
\begin{scope}[shift = {(-\yrl,-1.9954,-\xrl)}]
\coordinate (O) at (0,0,0);
\coordinate (A) at (0,\ztk,0);
\coordinate (B) at (0,\ztk,\xtk);
\coordinate (C) at (0,0,\xtk);
\coordinate (D) at (\ytk,0,0);
\coordinate (E) at (\ytk,\ztk,0);
\coordinate (F) at (\ytk,\ztk,\xtk);
\coordinate (G) at (\ytk,0,\xtk);

\coordinate (A2) at (0,3.3,0);
\coordinate (B2) at (0,3.3,\xtk);
\coordinate (E2) at (\ytk,3.3,0);
\coordinate (F2) at (\ytk,3.3,\xtk);

\coordinate (sAr) at (\yrl,\ztk,\xrl);
\coordinate (sBr) at (\yrl,\ztk,\xrr);
\coordinate (sEr) at (\yrr,\ztk,\xrl);
\coordinate (sFr) at (\yrr,\ztk,\xrr);

\coordinate (sArb) at (\yrl,1.9954,\xrl);
\coordinate (sBrb) at (\yrl,1.9954,\xrr);
\coordinate (sErb) at (\yrr,1.9954,\xrl);
\coordinate (sFrb) at (\yrr,1.9954,\xrr);

\coordinate (sArl) at (\yrl,3.3,\xrl);
\coordinate (sBrl) at (\yrl,3.3,\xrr);
\coordinate (sErl) at (\yrr,3.3,\xrl);
\coordinate (sFrl) at (\yrr,3.3,\xrr);


 \draw[blue,fill=yellow!40,opacity=0.7] (sArl) -- (sBrl) -- (sFrl) -- (sErl) -- cycle;

\draw[blue,fill=red!40,opacity=0.7] (sAr) -- (sBr) -- (sFr) -- (sEr) -- cycle;
\draw[blue,fill=red!40,opacity=0.7] (sArb) -- (sBrb) -- (sFrb) -- (sErb) -- cycle;

\draw[blue] (sAr) -- (sArb) ; 
\draw[blue] (sBr) -- (sBrb) ; 
\draw[blue] (sFr) -- (sFrb) ; 
\draw[blue] (sEr) -- (sErb) ;


\draw[black, -stealth,opacity=0.7] (2.6,3.4,3.2) -- (2.6,3.4,\xrl) ; 
\draw[black, -stealth,opacity=0.7] (2.6,3.4,3.2) -- (\yrr,3.4,3.2) ; 
\draw[black, -stealth,opacity=0.7] (2.6,3.4,3.2) -- (\yrr,3.4,\xrl) ;

\draw[black, dashed,opacity=0.7] (\yrl,3.4,\xrr - \lpml) -- (\yrr,3.4,\xrr - \lpml) ; 
\draw[black, dashed,opacity=0.7] (\yrl,3.4,\xrl + \lpml) -- (\yrr,3.4,\xrl + \lpml) ;

\draw[black, dashed,opacity=0.7] (\yrl + \lpml,3.4,\xrr) -- (\yrl + \lpml,3.4,\xrl) ; 
\draw[black, dashed,opacity=0.7] (\yrr - \lpml,3.4,\xrl) -- (\yrr - \lpml,3.4,\xrr) ;

\draw[black, dashed,opacity=0.7] (\yrl,1.9954 + \lpml,\xrr - \lpml) -- (\yrr,1.9954 + \lpml,\xrr - \lpml) ; 
\draw[black, dashed,opacity=0.7] (\yrl,1.9954 + \lpml,\xrl + \lpml) -- (\yrr,1.9954 + \lpml,\xrl + \lpml) ;

\draw[black, dashed,opacity=0.7] (\yrl + \lpml,1.9954 + \lpml,\xrr) -- (\yrl + \lpml,1.9954 + \lpml,\xrl) ; 
\draw[black, dashed,opacity=0.7] (\yrr - \lpml,1.9954 + \lpml,\xrl) -- (\yrr - \lpml,1.9954 + \lpml,\xrr) ;

\draw[black, dashed,opacity=0.7] (\yrl,1.9954,\xrr - \lpml) -- (\yrr,1.9954,\xrr - \lpml) ; 
\draw[black, dashed,opacity=0.7] (\yrl,1.9954,\xrl + \lpml) -- (\yrr,1.9954,\xrl + \lpml) ;

\draw[black, dashed,opacity=0.7] (\yrl + \lpml,1.9954,\xrr) -- (\yrl + \lpml,1.9954,\xrl) ; 
\draw[black, dashed] (\yrr - \lpml,1.9954,\xrl) -- (\yrr - \lpml,1.9954,\xrr) ;

\draw[black, dashed,opacity=0.7] (\yrl + \lpml,1.9954,\xrr - \lpml) -- (\yrl + \lpml,3.4,\xrr - \lpml);
\draw[black, dashed,opacity=0.7] (\yrl,1.9954,\xrr - \lpml) -- (\yrl,3.4,\xrr - \lpml);
\draw[black, dashed,opacity=0.7] (\yrl + \lpml,1.9954,\xrr) -- (\yrl + \lpml,3.4,\xrr);

\draw[black, dashed,opacity=0.7] (\yrl + \lpml,1.9954,\xrl + \lpml) -- (\yrl + \lpml,3.4,\xrl + \lpml);
\draw[black, dashed,opacity=0.7] (\yrl,1.9954,\xrl + \lpml) -- (\yrl,3.4,\xrl + \lpml);
\draw[black, dashed,opacity=0.7] (\yrl + \lpml,1.9954,\xrl) -- (\yrl + \lpml,3.4,\xrl);

\draw[black, dashed,opacity=0.7] (\yrr - \lpml,1.9954,\xrr - \lpml) -- (\yrr - \lpml,3.4,\xrr - \lpml);
\draw[black, dashed,opacity=0.7] (\yrr,1.9954,\xrr - \lpml) -- (\yrr,3.4,\xrr - \lpml);
\draw[black, dashed,opacity=0.7] (\yrr - \lpml,1.9954,\xrr) -- (\yrr - \lpml,3.4,\xrr);

\draw[black, dashed,opacity=0.7] (\yrr - \lpml,1.9954,\xrl + \lpml) -- (\yrr - \lpml,3.4,\xrl + \lpml);
\draw[black, dashed,opacity=0.7] (\yrr,1.9954,\xrl + \lpml) -- (\yrr,3.4,\xrl + \lpml);
\draw[black, dashed,opacity=0.7] (\yrr - \lpml,1.9954,\xrl) -- (\yrr - \lpml,3.4,\xrl);

 \node[circle, draw] at (2.6,3.2,3.2) () {H};

 \node[draw] at (2.6+0.7348,3.4,3.2-0.7348) () {6};
 \node[draw] at (2.6+0.5764,3.4,3.2-0.8647) () {9};

\end{scope}
\end{scope}

\draw[dotted, gray] (sAr) -- (Ar) ;
\draw[dotted, gray] (sBr) -- (Br) ;
\draw[dotted, gray] (sFr) -- (Fr) ;
\draw[dotted, gray] (sEr) -- (Er) ;

\draw[dotted, gray] (sArb) -- (Arb) ;
\draw[dotted, gray] (sBrb) -- (Brb) ;
\draw[dotted, gray] (sFrb) -- (Frb) ;
\draw[dotted, gray] (sErb) -- (Erb) ;

\end{tikzpicture}
    \caption{LOH1 problem setup (to scale) with the upper and lower block separated by the yellow interface. In blue is the suggested computational domain, in red is the computational domain we use with the PML. The PML regions are sectioned along the boundary of the enlarged red block. 
    The hypocentre (source location) is labelled H, and stations 6 and 9 are marked on the free-surface. The red region occupies approximately $ 3.8095 \%$  of the volume of the blue region. 
Due to the efficient absorption properties of the PML, the computational load for this problem is significantly reduced, thus saving as much as $96.19 \%  $ of  the required computational resources.}
    \label{fig:loh1_setup}
\end{figure}
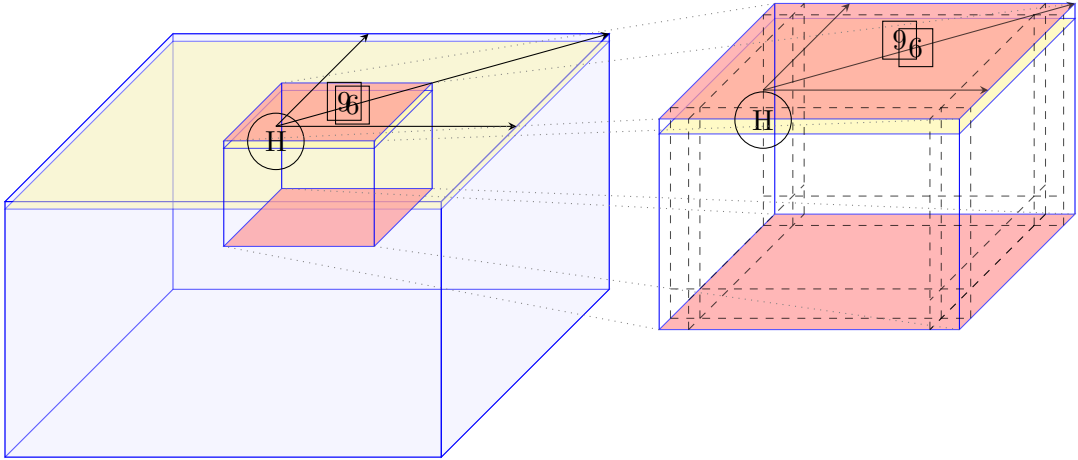

\begin{figure}[h!]
\begin{subfigure}{\textwidth}
    \centering
\stackunder[5pt]{\includegraphics[width=0.495\textwidth]{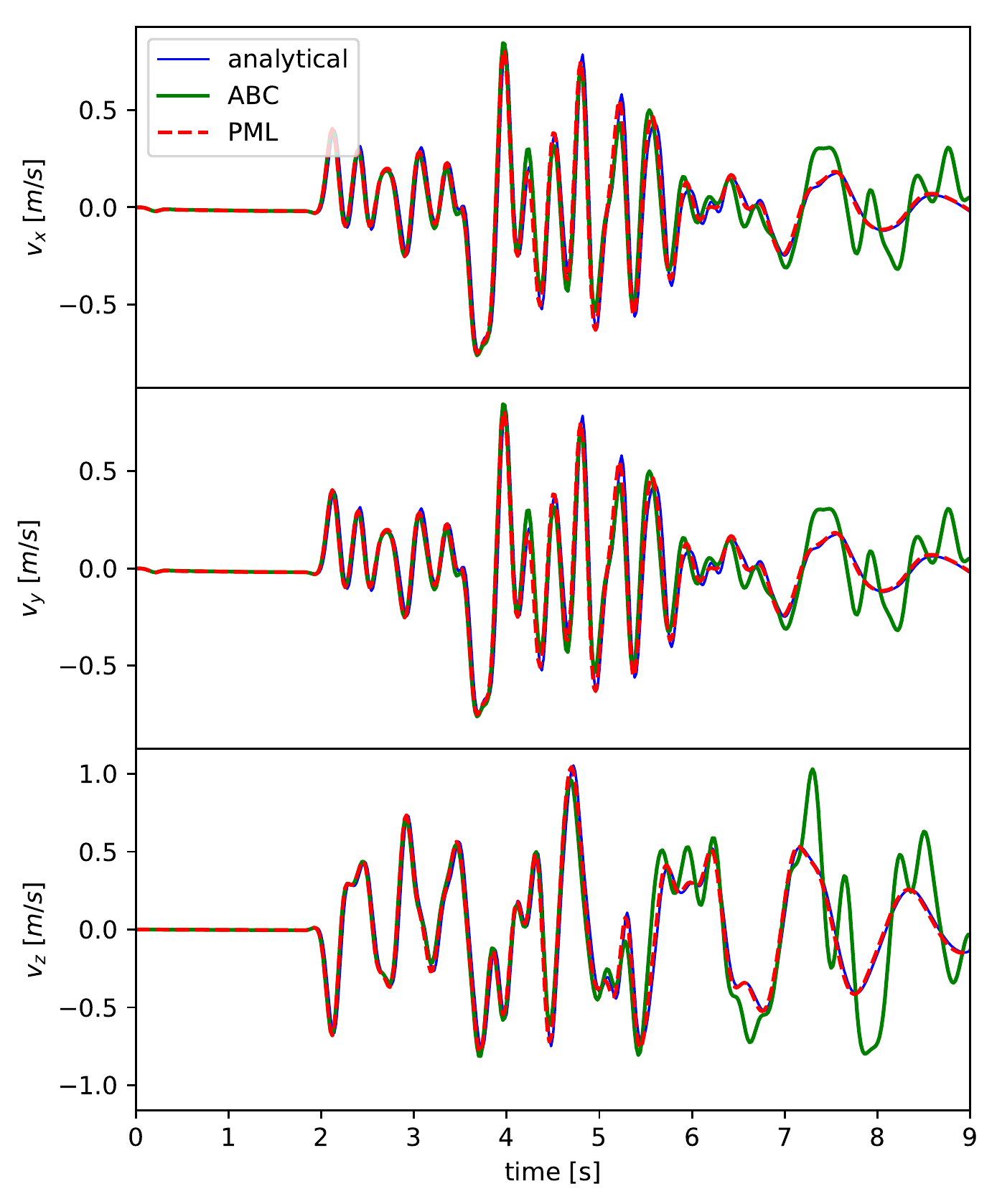}}{Receiver 6}%
\hspace{0.0cm}%
\stackunder[5pt]{\includegraphics[width=0.495\textwidth]{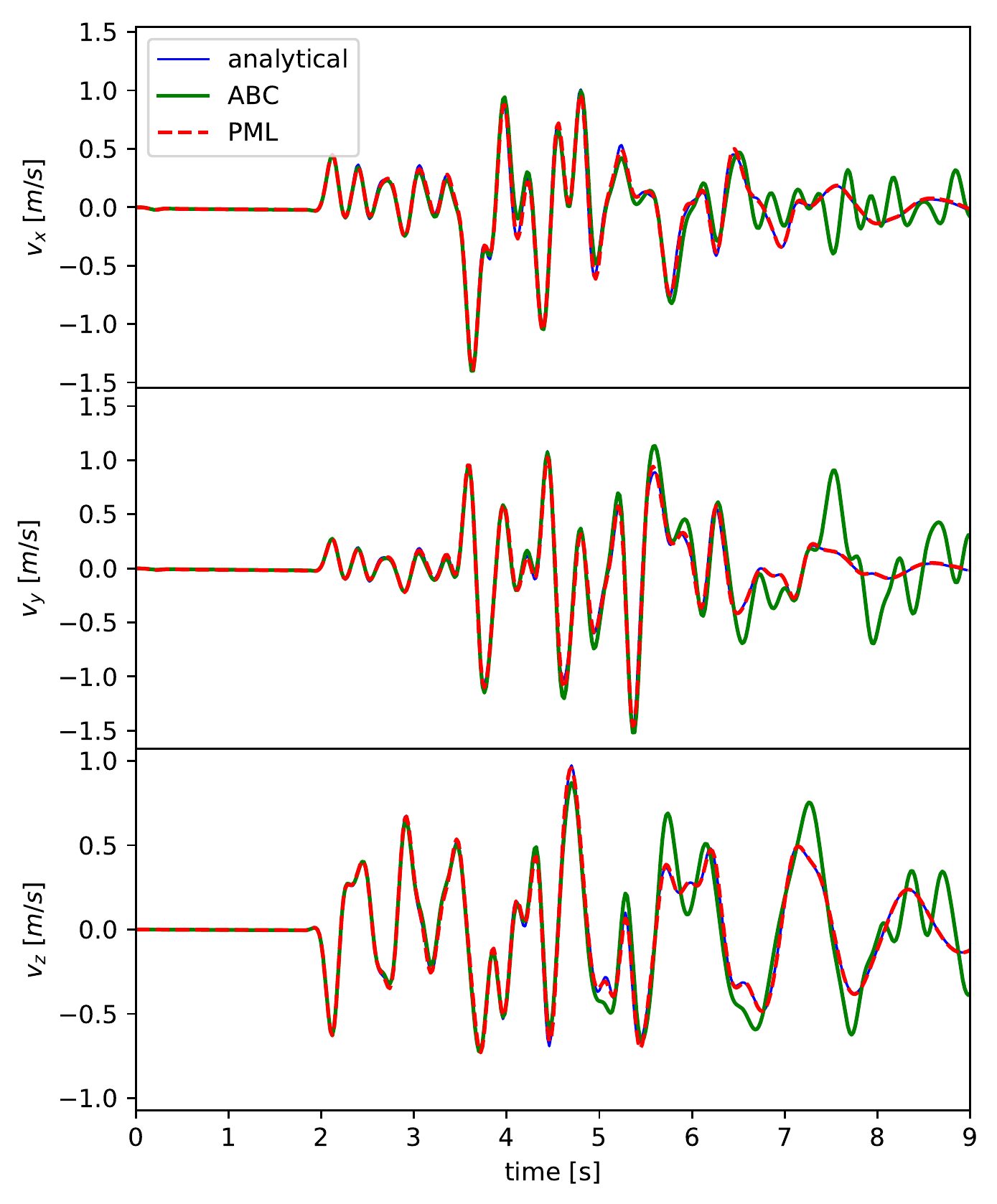}}{Receiver 9}%
     \end{subfigure}
    \caption{The LOH1 benchmark problem with degree $N = 5$ polynomial DG approximation. The PML solution coincides with the exact solution, while the ABC solution contains spurious reflections which can not dimininish with mesh refinement.}
    \label{fig:loh1_o5}
\end{figure}

We consider degree $N =5$ polynomial approximations, and run the simulation until the final time $t = 9~$s using the ExaHyPE software package \cite{ElasticDG_PML2019}. In Figure \ref{fig:loh1_o5} we display the solutions for Receivers 6 and 9.  Note that initially the ABC and PML solutions match the analytical solution very well. However, at later times the ABC solution  is polluted by numerical reflections arriving from the artificial boundaries.  As expected,  the PML solutions match the analytical solution excellently well, and remains accurate for the entire simulation duration. For the ABC, the dominant errors are the errors introduced by artificial reflections, these can never diminish with $p-$ or $h-$refinement.

These results and simulations have been replicated for the SBP-SAT finite difference method using the software package WaveQLab3D \cite{DuruFungWilliams2020}. The reader can consult the following papers \cite{ElasticDG_PML2019,Duru_exhype_2_2019,DuruFungWilliams2020} for more results and detailed discussions. 


\subsection{Large-scale numerical simulations in a 3D complex geometry}
We will now present numerical simulations in a complex geometry \cite{DuruFungWilliams2020,ExaHyPE2019}, with a geologically constrained complex non-planar free-surface topography.  Zugspitze is the tallest mountain in Germany, lying in the Wetterstein mountain range. The topography of this region is complex, with large variations in altitude across the Earths surface  \cite{Copernicus}. 

\begin{figure} [h!]
 \centering
{\includegraphics[width=.85\linewidth, angle=0]{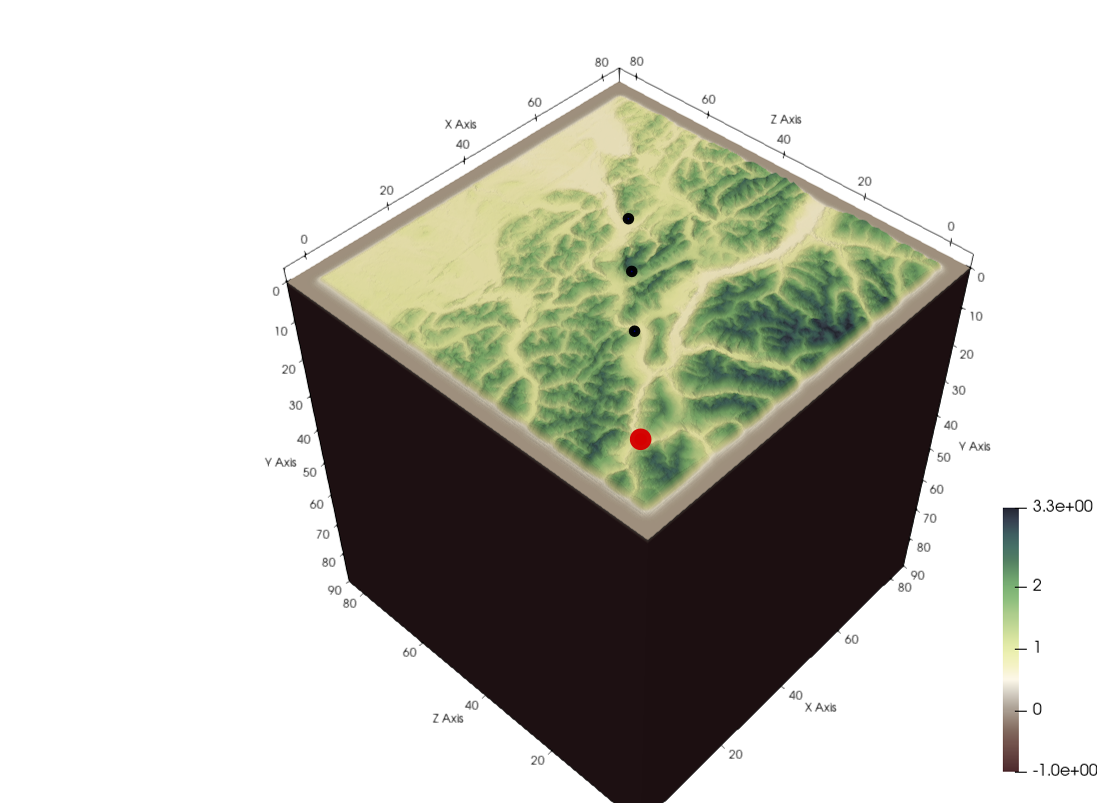}}
 \caption{The Mount Zugspitze model model set-up.  The red dot ({\red \bf $\cdot$}) at ($x=10$ km, $z=10$ km) depicts the epicenter  of a  moment tensor point source buried $10$ km at depth and the black dots ({\black \bf$\cdot$})  indicate receiver stations with on the free-surface Station 1: ($x=30$ km, $z=30$ km), Station 2: ($x=40$ km, $z=40$ km)  and Station 3: ($x=50$ km, $z=50$ km), which are the receiver stations where the solutions are sampled. Station 2 is collocated with the peak ($x=40$ km, $z=40$ km) of Mount Zugspitze. The boundaries of computational domain are surrounded by the PML to absorb outgoing waves. } 
 \label{fig:zugs_model}
\end{figure}

The Zugspitze model was set up to study the scattering effects of geometrically complex free-surface topography on the propagation of seismic wave fields in the European Alpine region.
The modelling domain is $\Omega = \bigcup_{y,z \in [-5,85] } [\widehat{X}(x,y),80] \times [-5,85]^2 $ with the $x$-co-ordinate being positive in-towards the Earth, like our previous example, and $\widehat{X}(y,z)$ parameterising the Earth's surface.
At each truncated boundary, in $y$-direction and $z$-direction, and down dip at $x = 80$~m, we have included a 5~km absorbing layer where PML boundary conditions \cite{DuruFungWilliams2020,ExaHyPE2019} are implemented  to prevent artificial reflections  from the computational boundaries from contaminating the solution. Please see Figure \ref{fig:zugs_model}. The PML relative error tolerance is $1\%$ which ensures the generation of high quality synthetic seismograms.

 
\begin{figure}[H]
    \centering
\stackunder[5pt]{\includegraphics[draft=false,width=0.95\columnwidth]{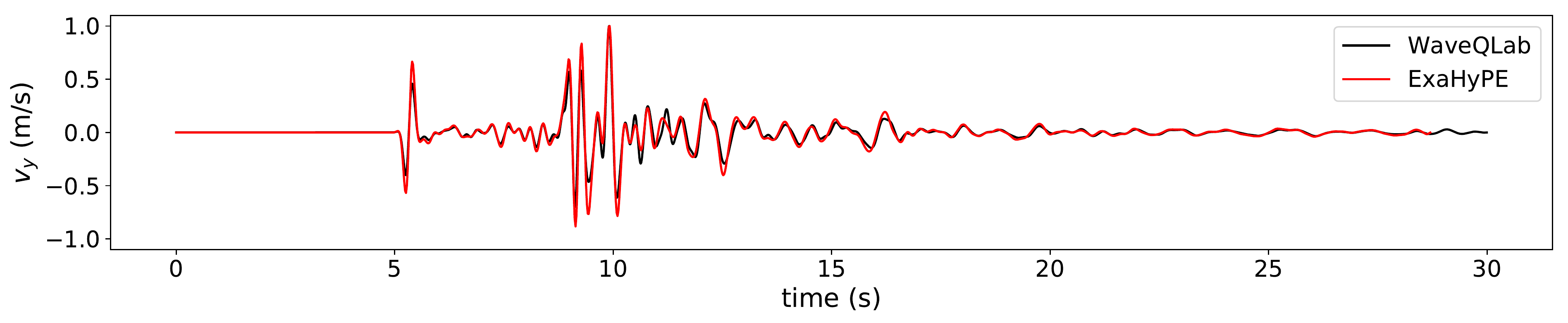}}{Station 1}%
\hspace{0.0cm}%
\stackunder[5pt]{\includegraphics[draft=false,width=0.95\columnwidth]{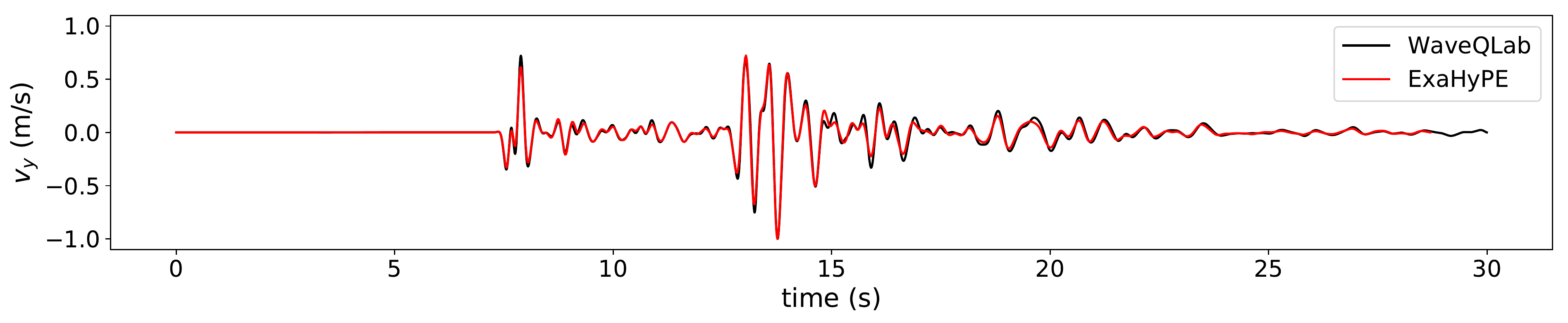}}{Station 2}%
\hspace{0.0cm}%
\stackunder[5pt]{\includegraphics[draft=false,width=0.95\columnwidth]{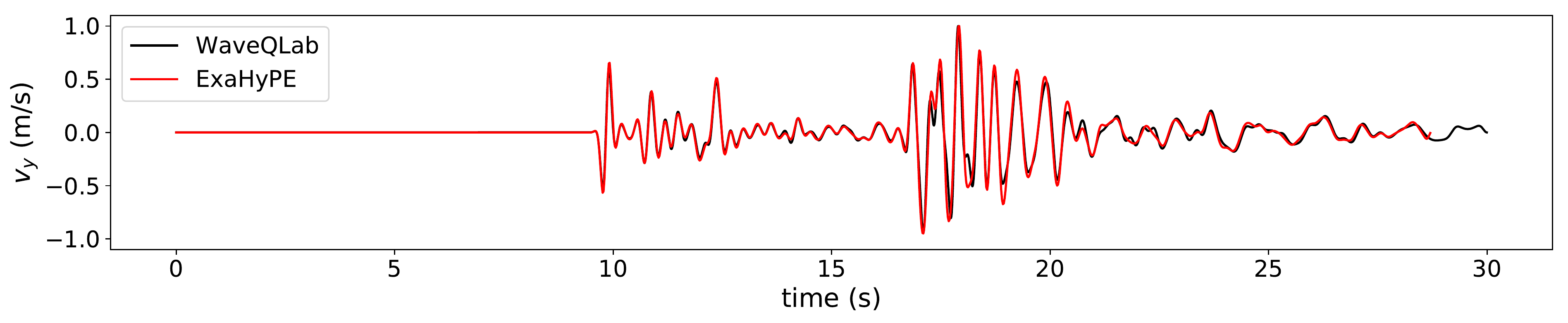}}{Station 3}%
    \caption{Seismograph from the 3 receiver stations placed on the Earths surface.}
    \label{fig:zugs_r1}
\end{figure}

We consider the homogenous material properties 

$$\rho = 2670~\ \text{kg/m}^3, \quad c_p = 6000~\ \text{m/s}, \quad c_s = 3464~\ \text{m/s}.$$
%
\begin{figure}[h!]
    \centering
    \includegraphics[width=0.485\textwidth]{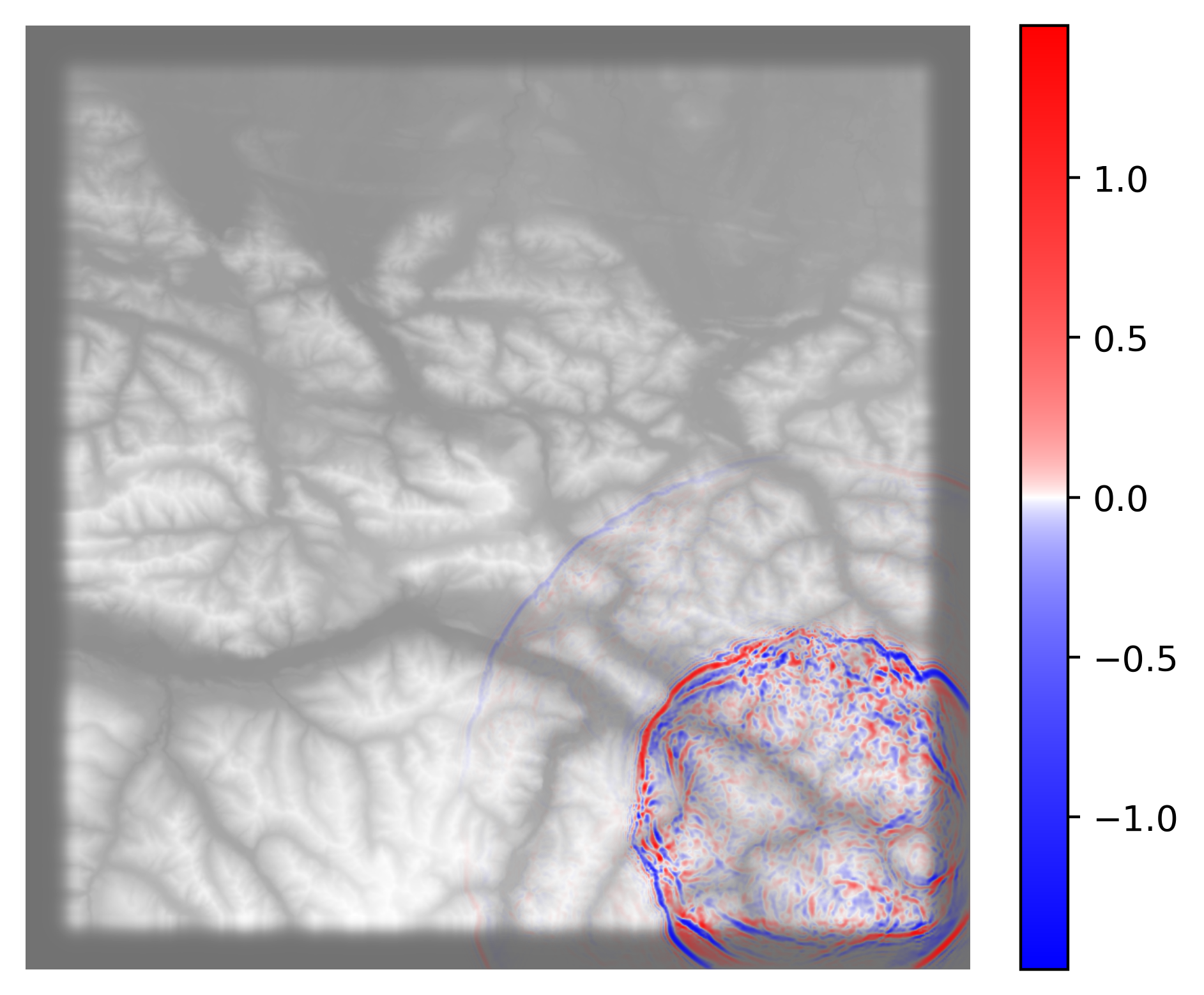}
    \includegraphics[width=0.485\textwidth]{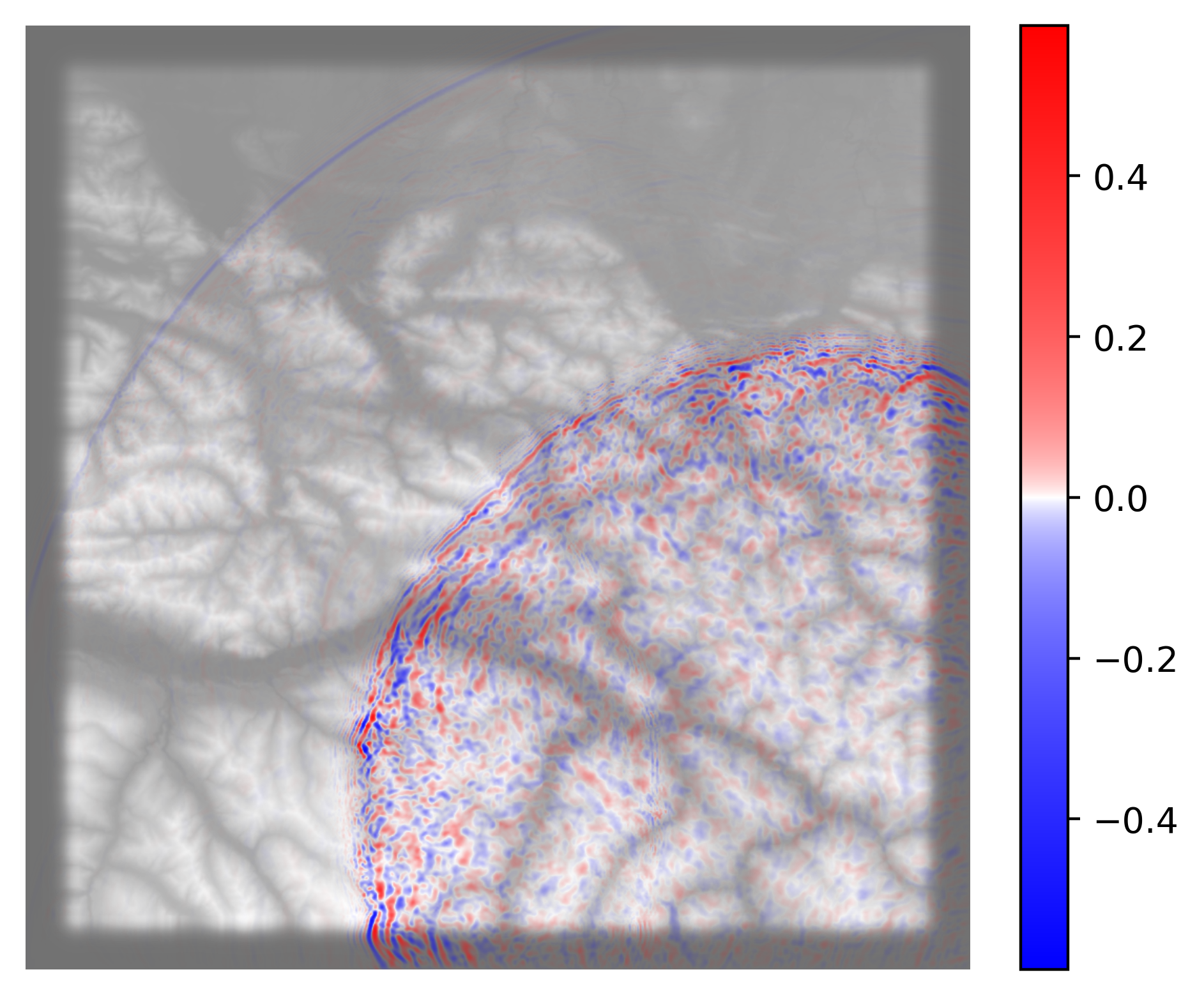}
    \caption{Simulations of high frequency  seismic wave scattering from complex free surface topography in a section of the European alpine region using the HPC software WaveQLab3D  \cite{DuruandDunham2016,DuruFungWilliams2020}. The background grey represents the altitude given by the underlying topography. We can see the scattering of high frequency waves conforming to the complex free-surface topography present and the efficient absorption properties  of the PML.  Note the absence of spurious reflected waves from the boundaries.
    } 
    \label{fig:topographydata}
\end{figure}
 The Zugspitze  simulations have been performed with both ExaHyPE and WaveQLab3D. We  run the simulation until the final time $t = 30$ s such that the elastic waves propagate through the media and leave the computational domain. 
As the waves propagate through the media, they interact with the complex topography and generate  high frequency scattered wave-modes. Because of the complex non-planar topography, the Zugspitze model has no analytical solution. 
 Snapshots of numerical seismic wave field on the free-surface topography are shown in Figure \ref{fig:topographydata}. Note the absence of spurious reflected waves from the boundaries. The seismograms are shown in Figure \ref{fig:zugs_r1}, for the 3 receiver stations.  We observe a near perfect agreement of the seismograms, for WaveQlab and ExaHyPE simulations, at sufficiently high frequencies. We refer the reader to \cite{ElasticDG_PML2019,Duru_exhype_2_2019,DuruFungWilliams2020} for more details and elaborate discussions.

We refer the reader to  \cite{DuruFungWilliams2020,ExaHyPE2019} for more details and results.

\section{Summary and conclusions}
The PML is now a standard  and flexible technique for constructing efficient and reliable domain truncation schemes for accurate numerical solution of  wave propagation problems.
However, the mathematical analysis and the development of provably stable and accurate numerical  approximation of the PML is a challenge in many applications. In this review article we have summarised the progress made, from mathematical, numerical and practical perspectives, point out open problems  and set the stage for future work.

The overarching results are are obtained from the use of mode analysis to prove well-posedness and stability of the PML IVPs and IBVPs. The result is that as long as the underlying hyperbolic IVPs do not violate the {\it geometric stability condition} and boundary conditions stable and well-posed without the PML,  the PML transformation will not move the wave modes into the unstable region in the complex plane. However, these results are   too technical to be extended to the analysis of numerical approximations.

Further, we review extensions of the results using the energy method in the Laplace space.  The energy estimates enable the development of stable and accurate numerical methods for the PML using SBP finite difference methods and DG methods. Numerical experiments in acoustic and elastic media verify the theoretical results.

Finally,  we perform numerical simulations of seismological application problems and demonstrate impact. We consider both a standard 3D seismological benchmark problem, the LOH1 \cite{Seismowine, Kristekova_etal2009, Kristekova_etal2006}, and a real-world wave propagation propagation problem which involves the simulation of 3D seismic waves in a section of European Alpine region, with strong non-planar free-surface topography. 
For the LOH1 benchmark problem, the PML allows us to sufficiently limit the domain to $ 3.8095 \%$ of the suggested large domain $\Omega_{L}$, thus saving as much as $96.19 \%  $ of  the required computational resources.
The algorithms and the PML have been implemented in two different freely open source  HPC software packages, WaveQLab3D \cite{DuruFungWilliams2020} and ExaHyPE \cite{ElasticDG_PML2019}, for large-scale simulation of seismic waves in geometrically complex 3D Earth models. The software package WaveQLab3D \cite{DuruFungWilliams2020} is  a  high order accurate SBP-SAT finite difference solver. ExaHyPE is a DG solver of arbitrary accuracy for large-scale numerical simulation of hyperbolic wave propagation problems on dynamically adaptive curvilinear meshes.

\section*{Acknowledgements}
Some of the works of the authors summarised in this review span over ten years.
The first author KD grateful acknowledges support from Uppsala University, Stanford University, Ludwig Maximilian University of Munich  and The Australian National University. This research was undertaken with the assistance of resources and services from the National Computational Infrastructure (NCI, Project fp92 on Gadi) which is supported by the Australian Government's National Collaborative Research Infrastructure Strategy (NCRIS). The authors also gratefully acknowledge the Gauss Centre for Supercomputing e.V. \footnote{www.gauss-centre.eu} for funding this project by providing computing time on the GCS Supercomputer SuperMUC-NG (project pr63qo) at Leibniz Supercomputing Centre \footnote{www.lrz.de}. 

\bibliographystyle{plain}
\bibliography{references}

\end{document}